\documentclass[12pt]{article}
\usepackage{amstex,amssymb}
\begin{document}
\newtheorem{thm}{Theorem}[subsection]
\newtheorem{lem}[thm]{Lemma}
\newtheorem{prop}[thm]{Proposition}
\newtheorem{cor}[thm]{Corollary}
\newenvironment{dfn}{\medskip\refstepcounter{thm}
\noindent{\bf Definition \thesubsection.\arabic{thm}\ }}{\medskip}
\newenvironment{ex}{\medskip\refstepcounter{thm}
\noindent{\bf Example \thesubsection.\arabic{thm}\ }}{\medskip}
\newenvironment{conj}{\medskip\refstepcounter{thm}
\noindent{\bf Conjecture \thesubsection.\arabic{thm}\ }}{\medskip}
\newenvironment{proof}{\medskip\noindent{\it Proof.}}{\hfill
$\square$\medskip}
\newenvironment{ax}[1]{\begin{axnumerate}
\item[{\bf Axiom #1. $(i)$}]
\addtocounter{enumi}{1}}{\end{axnumerate}}
\newenvironment{axw}[1]{\begin{axnumerate}\item[{\bf Axiom #1.}]}{\end{axnumerate}}
\newenvironment{axnumerate}{\begin{list}{$(\roman{enumi})$}
{\usecounter{enumi}
\setlength{\leftmargin}{22pt}
\setlength{\parsep}{1pt}
\setlength{\itemsep}{1pt}}}{\end{list}}
\def\bb{\mathbb}
\def\go{\mathfrak}
\def\md#1{\vert #1 \vert}
\def\eq#1{{\rm(\ref{#1})}}
\def\dim{\mathop{\rm dim}}
\def\mod{\mathop{\rm mod}}
\def\Im{\mathop{\rm Im}}
\def\Ker{\mathop{\rm Ker}}
\def\Aut{\mathop{\rm Aut}}
\def\Vect{\mathop{\rm Vect}}
\def\id{\mathop{\rm id}}
\def\ah{{\rm A}$\bb H$}
\def\H{{\mathbin{\bb H}}}
\def\R{{\mathbin{\bb R}}}
\def\Z{{\mathbin{\bb Z}}}
\def\C{{\mathbin{\bb C}}}
\def\I{{\mathbin{\bb I}}}
\def\d{\dagger}
\def\br{\buildrel}
\def\ov{\overline}
\def\ot{\otimes}
\def\ra{\rightarrow}
\def\longra{\longrightarrow}
\def\t{\times}
\def\ha{{\textstyle{1\over2}}}
\def\oth{{\textstyle\otimes_\H}}
\def\smalloth{{\scriptstyle\otimes_{\scriptscriptstyle\H}}}
\def\bigoth{{\textstyle\bigotimes}_\H}
\def\Lambh{{\textstyle\Lambda}_\H}
\def\Sh{S_\H}
\def\op{\oplus}
\def\g{{\go g}}

\title{A theory of quaternionic algebra, \\ with applications to 
hypercomplex geometry}
\author{Dominic Joyce, Lincoln College, Oxford}
\date{}
\maketitle

\section{Introduction}

\begin{quote}
{\it Mathematicians are like Frenchmen: whenever you say something
to them, they translate it into their own language, and at
once it is something entirely different.}

\hfill Goethe, Maxims and Reflections (1829)
\end{quote}

The subject of this paper is an algebraic device, a means
to construct algebraic structures over the quaternions $\H$
as though $\H$ were a commutative field. As far as the author can 
tell, the idea seems to be new. We shall provide the reader with 
a dictionary, giving the equivalents of simple concepts such as 
commutative field, vector space, tensor products of vector spaces, 
symmetric and antisymmetric products, dimensions of vector spaces, 
and so on. The reader will then be able to translate her own favourite
algebraic objects into this new quaternionic language. The results
often turn out to be surprisingly, entirely different.

The basic building blocks of the theory are the quaternionic analogues 
of vector space and the tensor product of vector spaces over a 
commutative field. A vector space is replaced by an {\it \ah-module
$(U,U')$}, which is a left $\H$-module $U$ with a given real 
vector subspace $U'$. The tensor product $\ot$ is replaced by 
the {\it quaternionic tensor product $\oth$}, which has a complex 
definition given in \S\ref{h11}. It shares some important 
properties of $\ot$ (e.g.~it is commutative and associative), 
but also has important differences (e.g.~the dimension of a quaternionic 
tensor product behaves strangely). 

The theory arose out of my attempts to understand
the algebraic structure of noncompact hypercomplex manifolds. Let $M$ be 
a $4n$-manifold. A {\it hypercomplex structure} on $M$ is a triple 
$(I_1,I_2,I_3)$ of complex structures satisfying $I_1I_2=I_3$. 
These induce an $\H$-action on the tangent bundle $TM$, and
$M$ is called a hypercomplex manifold. Hypercomplex manifolds are the subject of 
the second half of the paper.

The best quaternionic analogue of a holomorphic function on a complex 
manifold is a {\it q-holomorphic function} on a hypercomplex manifold 
$M$, defined in \S\ref{h31}. This is an $\H$-valued function 
on $M$ satisfying an equation analogous to the Cauchy-Riemann
equations, which was introduced by Fueter in 1935 for the case 
$M=\H$. Affine algebraic geometry is the study of complex 
manifolds using algebras of holomorphic functions upon them. 
Seeking to generalize this to hypercomplex manifolds using q-holomorphic 
functions, I was led to this quaternionic theory of algebra as the
best language to describe quaternionic algebraic geometry.

Although the applications are all to hypercomplex manifolds, I hope that 
much of the paper will be of interest to those who study algebra 
rather than geometry, and maybe to others who, like myself, are 
fascinated by the quaternions. The paper has been laid out with
this in mind. There are four chapters. Chapter 1 explains 
\ah-modules and the quaternionic tensor product. It is quite long and
wholly algebraic, without references to geometry. Chapter 2 
gives quaternionic analogues of various algebraic structures. It is short
and is mostly definitions. The most important idea is that of
an H-algebra, the quaternionic version of a commutative algebra.

Chapter 3 is about hypercomplex geometry. Q-holomorphic functions are
defined, their properties explored, and it is shown that the
q-holomorphic functions on a hypercomplex manifold form an H-algebra.
A similar result is proved for hyperk\"ahler manifolds. The problem of 
reconstructing a hypercomplex manifold from its H-algebra is considered, 
and H-algebras are used to study a special class of noncompact,
complete hyperk\"ahler manifolds, called asymptotically conical manifolds.

Chapter 4 is a collection of examples and applications of the 
theory. Two interesting topics covered here are an algebraic
treatment of the `coadjoint orbit' hyperk\"ahler manifolds, and some new
types of singularities of hypercomplex manifolds that have remarkable 
properties. To control the length of the paper I have kept the 
list of examples and applications short, and I have omitted
a number of proofs. However, I believe that there is much
interesting work still to be done on these ideas, and in
\S\ref{h47} I shall indicate some directions in which I would 
like to see the subject develop.

Much of the material in this paper has already been published by
the author in \cite{Joy2}; this paper is mostly an expanded version
of \cite{Joy2}, with rather more detail, and some new material on
the geometry of noncompact hypercomplex and hyperk\"ahler manifolds.
However, the results of \cite[\S 12]{Joy2} on coadjoint orbits are
not included in this paper, and instead we approach the subject
from a slightly different direction.

Two other references on the subject of this paper are Quillen
\cite{Qui}, who reinterprets the \ah-modules and the quaternionic 
tensor product $\oth$ in terms of sheaves on $\mathbb{CP}^1$, and 
Widdows \cite{Wid}, who explores a number of issues in hypercomplex
algebraic geometry, including the classification of finite-dimensional
\ah-modules up to isomorphism, and quaternionic analogues of the 
Dolbeault double complex on complex manifolds.

\subsection{Quaternionic tensor products}
\label{h11}

In this section we give some notation, and define
a quaternionic analogue of the tensor product of two vector spaces.
This idea is central to the whole paper. First, a remark
about the real tensor product. Let $U,V$ be real vector 
spaces. If $U,V$ are infinite-dimensional,
there is more than one possible definition for the
real tensor product $U\ot V$, if $U,V$ are equipped
with topologies. In this paper we choose the simplest 
definition: for us, every element of $U\ot V$ is a
{\it finite} sum $\Sigma_ju_j\ot v_j$ where $u_j\in U$,
$v_j\in V$. Also, in this paper the dual $U^*$ of $U$ 
means the vector space of {\it all} linear maps $U\ra\R$,
whether continuous in some topology or not. 

The quaternions are 
\begin{equation*}
\H=\{r_0+r_1i_1+r_2i_2+r_3i_3:r_0,\dots,r_3\in{\bb R}\},
\end{equation*}
and quaternion multiplication is given by
\begin{equation*}
i_1i_2=-i_2i_1=i_3,\;
i_2i_3=-i_3i_2=i_1,\;
i_3i_1=-i_1i_3=i_2,\;
i_1^2=i_2^2=i_3^2=-1.
\end{equation*}
The quaternions are an associative, noncommutative algebra.
The imaginary quaternions are $\I=\langle i_1,i_2,i_3\rangle$.
This notation $\I$ is not standard, but we will use it throughout
the paper. If $q=r_0+r_1i_1+r_2i_2+r_3i_3$ then we define the
{\it conjugate} $\ov q$ of $q$ by $\ov q=r_0-r_1i_1-r_2i_2-r_3i_3$. 
Then $\ov{(pq)}=\ov q\,\ov p$ for~$p,q\in\H$.

A {\it (left) $\H$-module} is a real vector space $U$ with an action
of $\H$ on the left. We write this action $(q,u)\mapsto q\cdot u$
or $qu$, for $q\in\H$ and $u\in U$. The action is a bilinear
map $\H\times U\ra U$, and satisfies 
$p\cdot(q\cdot u)=(pq)\cdot u$ for $p,q\in\H$ and $u\in U$.
In this paper, all $\H$-modules will be left $\H$-modules.

Let $U$ be an $\H$-module. We define the {\it dual $\H$-module $U^\t$}
to be the vector space of linear maps $\alpha:U\ra\H$ that
satisfy $\alpha(qu)=q\alpha(u)$ for all $q\in\H$ and $u\in U$.
If $q\in\H$ and $\alpha\in U^\t$ we may define $q\cdot\alpha$ by 
$(q\cdot\alpha)(u)=\alpha(u)\ov q$ for $u\in U$. Then 
$q\cdot\alpha\in U^\t$, and $U^\t$ is a (left) $\H$-module.
If $V$ is a real vector space, we write $V^*$ for the
dual of $V$ as a real vector space. It is important to 
distinguish between the dual vector spaces and dual $\H$-modules.
Dual $\H$-modules behave just like dual vector spaces.
In particular, there is a canonical map $U\ra (U^\t)^\t$,
that is an isomorphism when $U$ is finite-dimensional.

\begin{dfn} Let $U$ be an $\H$-module. Let $U'$ be a real 
vector subspace of $U$, that need not be closed under the 
$\H$-action. Define a real vector subspace $U^\d$ of $U^\t$ by
\begin{equation*}
U^\d=\bigl\{\alpha\in U^\t:\text{$\alpha(u)\in\I$ for all 
$u\in U'$}\bigr\}.
\end{equation*}
We define an {\it augmented $\H$-module}, or {\it \ah-module}, to 
be a pair $(U,U')$, such that if $u\in U$ and $\alpha(u)=0$ for 
all $\alpha\in U^\d$, then $u=0$. Usually we will refer to $U$ as 
an \ah-module, implicitly assuming that $U'$ is also given. 
We consider $\H$ to be an \ah-module, with~$\H'=\I$.
\label{ahmdef}
\end{dfn}

\ah-modules should be thought of as the quaternionic analogues of
real vector spaces. It is easy to define the dual of an 
\ah-module, but we are not going to do this, as it seems
not to be a fruitful idea. We can interpret $U^\t$ as the
dual of $U$ as a real vector space, and then $U^\d$ is the 
annihilator of $U'$. Thus if $U$ is finite-dimensional, $\dim 
U'+\dim U^\d=\dim U=\dim U^\t$. The letters $A,U,V,\dots,Z$ 
will usually denote \ah-modules. Here are the natural concepts 
of linear map between \ah-modules, and \ah-submodules.

\begin{dfn} Let $U,V$ be \ah-modules. Let $\phi:U\ra V$
be a linear map satisfying $\phi(qu)=q\phi(u)$ for each $q\in\H$
and $u\in U$. Such a map is called {\it quaternion linear},
or {\it $\H$-linear}. We say that $\phi$ is a {\it morphism of 
\ah-modules}, or {\it \ah-morphism}, if $\phi:U\ra V$ is $\H$-linear
and satisfies $\phi(U')\subset V'$. Define a linear map 
$\phi^\t:V^\t\ra U^\t$ by $\phi^\t(\beta)(u)=\beta(\phi(u))$ 
for $\beta\in V^\t$ and $u\in U$. Then $\phi(U')\subset V'$ 
implies that $\phi^\t(V^\d)\subset U^\d$. If $\phi$ is an isomorphism
of $\H$-modules and $\phi(U')=V'$, we say that $\phi$ is an 
{\it \ah-isomorphism}.

Clearly, if $\phi:U\ra V$ and $\psi:V\ra W$ are \ah-morphisms, 
then $\psi\circ\phi:U\ra W$ is an \ah-morphism. If $V$ is an 
\ah-module, we say that {\it $U$ is an \ah-submodule of\/ $V$} if 
$U$ is an $\H$-submodule of $V$ and $U'=U\cap V'$. This implies
that $U^\d$ is the restriction of $V^\d$ to $U$, so the condition 
that $u=0$ if $\alpha(u)=0$ for all $\alpha\in U^\d$ holds 
automatically, and $U$ is an \ah-module.
\end{dfn}

In this paper `$\Im$' always means the image of a map, and
we will write `$\id$' for the identity map on any vector 
space. If $U$ is an \ah-module, then $\id:U\ra U$ is an 
\ah-morphism. Next we will define a sort of tensor product 
of \ah-modules, which is the key to the whole paper.

\begin{dfn} Let $U$ be an \ah-module. Then $\H\ot(U^\d)^*$
is an $\H$-module, with $\H$-action $p\cdot(q\ot x)
=(pq)\ot x$. Define a map $\iota_U:U\ra\H\ot(U^\d)^*$ 
by $\iota_U(u)\cdot\alpha=\alpha(u)$, for $u\in U$ and 
$\alpha\in U^\d$. Then $\iota_U(q\cdot u)=q\cdot\iota_U(u)$
for $u\in U$. Thus $\iota_U(U)$ is an $\H$-submodule of
$\H\ot(U^\d)^*$. Suppose $u\in\Ker\iota_U$. 
Then $\alpha(u)=0$ for all $\alpha\in U^\d$, so that 
$u=0$ as $U$ is an \ah-module. Thus $\iota_U$ is injective,
and~$\iota_U(U)\cong U$.
\label{iotadef}
\end{dfn}

\begin{dfn} Let $U,V$ be \ah-modules. Then 
$\H\ot(U^\d)^*\ot(V^\d)^*$ is an $\H$-module, with
$\H$-action $p\cdot(q\ot x\ot y)=(pq)\ot x\ot y$.
Exchanging the factors of $\H$ and $(U^\d)^*$, we may regard
$(U^\d)^*\ot\iota_V(V)$ as a subspace of $\H\ot(U^\d)^*\ot(V^\d)^*$.
Thus $\iota_U(U)\ot(V^\d)^*$ and $(U^\d)^*\ot\iota_V(V)$ 
are $\H$-submodules of $\H\ot(U^\d)^*\ot(V^\d)^*$.
Define an $\H$-module $U\oth V$ by
\begin{equation}
\!\!\!\!\!\!\!\!\!\!\!\!\!\!\!\!
U\oth V=\bigl(\iota_U(U)\ot(V^\d)^*\bigr)
\cap\bigl((U^\d)^*\ot\iota_V(V)\bigr)
\subset\H\ot(U^\d)^*\ot(V^\d)^*.
\label{uothveq}
\end{equation}
Define a vector subspace $(U\oth V)'$ by $(U\oth V)'=(U\oth V)
\cap\bigl(\I\ot(U^\d)^*\ot(V^\d)^*\bigr)$. Define a linear map 
$\lambda_{U,V}:U^\d\ot V^\d\ra (U\oth V)^\t$ by 
$\lambda_{U,V}(x)(y)=y\cdot x\in\H$, for $x\in U^\d\ot V^\d$, 
$y\in U\oth V$, where `$\cdot$' contracts together the factors 
of $U^\d\ot V^\d$ and~$(U^\d)^*\ot(V^\d)^*$. 

Clearly, if $x\in U^\d\ot V^\d$ and $y\in(U\oth V)'$, then
$\lambda_{U,V}(x)(y)\in\I$. As this holds for all $y\in(U\oth V)'$,
$\lambda_{U,V}(x)\in(U\oth V)^\d$, so that $\lambda_{U,V}$ maps
$U^\d\ot V^\d\ra (U\oth V)^\d$. If $y\in U\oth V$, then 
$\lambda_{U,V}(x)(y)=0$ for all $x\in U^\d\ot V^\d$ if and only 
if $y=0$. Thus $U\oth V$ is an \ah-module, by Definition \ref{ahmdef}. 
This \ah-module will be called the {\it quaternionic tensor product of\/ $U$ 
and\/ $V$}, and the operation $\oth$ will be called the {\it quaternionic 
tensor product}. When $U,V$ are finite-dimensional, $\lambda_{U,V}$
is surjective, so that~$(U\oth V)^\d=\lambda_{U,V}(U^\d\ot V^\d)$.
\label{qctpdef}
\end{dfn}

Here are some basic properties of the operation~$\oth$.

\begin{lem} Let $U,V,W$ be \ah-modules. Then there are
canonical \ah-isomorphisms
\begin{equation*}
\H\oth U\cong U,\quad
U\oth V\cong V\oth U,\quad\text{and}\quad
(U\oth V)\oth W\cong U\oth (V\oth W).
\end{equation*}
\label{tenslem}
\end{lem}

\begin{proof} As $\H^\d\cong\R$, we may identify 
$\H\ot(\H^\d)^*\ot(U^\d)^*$ and $\H\ot(U^\d)^*$. Under 
this identification, it is easy to see that $\H\oth U$ and 
$\iota_U(U)$ are identified. Since $\iota_U(U)\cong U$ as 
in Definition \ref{iotadef}, this gives an isomorphism 
$\H\oth U\cong U$, which is an \ah-isomorphism. The 
\ah-isomorphism $U\oth V\cong V\oth U$ is trivial, because 
the definition of $U\oth V$ is symmetric in $U$ and~$V$. 

It remains to show that $(U\oth V)\oth W\cong U\oth (V\oth W)$.
The maps $\lambda_{U,V}:U^\d\ot V^\d\ra(U\oth V)^\d$ and 
$\lambda_{U\smalloth V,W}:(U\oth V)^\d\ot W^\d\ra
((U\oth V)\oth W)^\d$ compose to give a linear map 
$\lambda_{UV,W}:U^\d\ot V^\d\ot W^\d\ra((U\oth V)\oth W)^\d$, 
defined in the obvious way. Define a linear map $\iota_{UV,W}:
(U\oth V)\oth W\ra\H\ot(U^\d)^*\ot(V^\d)^*\ot(W^\d)^*$ by
$\iota_{UV,W}(y)\cdot x=\lambda_{UV,W}(x)(y)\in\H$, for
each~$x\in U^\d\ot V^\d\ot W^\d$.

Suppose that $y\in(U\oth V)\oth W$, and $\iota_{UV,W}(y)=0$.
Then $\lambda_{UV,W}(x)(y)=0$ for each $x\in U^\d\ot V^\d\ot W^\d$.
It can be shown that this implies that $y=0$. Thus $\iota_{UV,W}$ 
is injective. Now from \eq{uothveq} and the definitions, it is 
easy to show that
\begin{equation}
\begin{split}
\iota_{UV,W}\bigl(&(U\oth V)\oth W\bigr)
=\bigl(\iota_U(U)\ot(V^\d)^*\ot(W^\d)^*\bigr)\\
&\cap\bigl((U^\d)^*\ot\iota_V(V)\ot(W^\d)^*\bigr)
\cap\bigl((U^\d)^*\ot(V^\d)^*\ot\iota_W(W)\bigr),
\end{split}
\label{uvweq}
\end{equation}
interpreting this equation as we did \eq{uothveq}. As 
$\iota_{UV,W}$ is injective, $(U\oth V)\oth W$ is isomorphic to 
the r.h.s.~of \eq{uvweq}. By the same argument, $U\oth(V\oth W)$ 
is also isomorphic to the r.h.s.~of \eq{uvweq}. This gives
a canonical isomorphism $(U\oth V)\oth W\cong U\oth(V\oth W)$.
It turns out to be an \ah-isomorphism, and the lemma is complete.
\end{proof}

Lemma \ref{tenslem} tells us that $\oth$ is commutative and
associative, and that $\H$ acts as an identity element for $\oth$.
Since $\oth$ is associative, we shall not bother to put
brackets in multiple products such as $U\oth V\oth W$.
Also, the commutativity and associativity of $\oth$ enable
us to define symmetric and antisymmetric products, analogous
to $S^kV$ and $\Lambda^kV$ for $V$ a real vector space.

\begin{dfn} Let $U$ be an \ah-module. Write $\bigoth^kU$ for
the product $U\oth\cdots\oth U$ of $k$ copies of $U$. Then
the $k^{\rm th}$ symmetric group $S_k$ acts on 
$\bigoth^kU$ by permutation of the $U$ factors in the
obvious way. Define $\Sh^kU$ to be the \ah-submodule of
$\bigoth^kU$ that is symmetric under these permutations,
and $\Lambh^kU$ to be the \ah-submodule of $\bigoth^kU$ that is
antisymmetric under these permutations. Define $\bigoth^0U$,
$\Sh^0U$ and $\Lambh^0U$ to be the \ah-module $\H$. The 
{\it symmetrization operator $\sigma_\H$}, defined in the 
obvious way, is a projection $\sigma_\H:\bigoth^kU
\ra\Sh^kU$. Clearly, $\sigma_\H$ is an \ah-morphism. 
Similarly, there is an {\it antisymmetrization operator}, that 
is an \ah-morphism projection from $\bigoth^kU$ to~$\Lambh^kU$.
\label{symmantisymm}
\end{dfn}

Here is the definition of the tensor product of two \ah-morphisms.

\begin{dfn} Let $U,V,W,X$ be \ah-modules, and let 
$\phi:U\ra W$ and $\psi:V\ra X$ be \ah-morphisms. 
Then $\phi^\t(W^\d)\subset U^\d$ and $\psi^\t(X^\d)\subset V^\d$, 
by definition. Taking the duals gives maps $(\phi^\t)^*:(U^\d)^*
\ra(W^\d)^*$ and $(\psi^\t)^*:(V^\d)^*\ra(X^\d)^*$. 
Combining these, we have a map
\begin{equation}
\id\ot(\phi^\t)^*\ot(\psi^\t)^*:
\H\ot(U^\d)^*\ot(V^\d)^*\ra
\H\ot(W^\d)^*\ot(X^\d)^*.
\label{idphipsi}
\end{equation}
Now $U\oth V\subset\H\ot(U^\d)^*\ot(V^\d)^*$
and $W\oth X\subset\H\ot(W^\d)^*\ot(X^\d)^*$. It is easy to
show that $\bigl(\id\ot(\phi^\t)^*\ot(\psi^\t)^*\bigr)
(U\oth V)\subset W\oth X$. Define 
$\phi\oth\psi:U\oth V\ra W\oth X$ to
be the restriction of $\id\ot(\phi^\t)^*\ot(\psi^\t)^*$
to $U\oth V$. It follows trivially from the definitions
that $\phi\oth\psi$ is $\H$-linear and satisfies 
$(\phi\oth\psi)\bigl((U\oth V)'\bigr)\subset
(W\oth X)'$. Thus $\phi\oth\psi$ is an \ah-morphism
from $U\oth V$ to $W\oth X$. This is the {\it quaternionic tensor
product of\/ $\phi$ and\/~$\psi$}.
\end{dfn}

\begin{lem} Suppose that $\phi:U\ra W$ and
$\psi:V\ra X$ are injective \ah-morphisms.
Then $\phi\oth\psi:U\oth V\ra W\oth X$ is an
injective \ah-morphism.
\label{injlem}
\end{lem}

\begin{proof} Consider the map $\id\ot(\phi^\t)^*\ot(\psi^\t)^*$ 
of \eq{idphipsi}. Clearly this maps $\iota_U(U)\ot(V^\d)^*$ 
to $\iota_W(W)\ot(X^\d)^*$. As $\iota_U(U)\cong U$ and 
$\iota_W(W)\cong W$ and the map $\phi:U\ra W$ is injective, 
we see that the kernel of $\id\ot(\phi^\t)^*\ot(\psi^\t)^*$ on
$\iota_U(U)\ot(V^\d)^*$ is $\iota_U(U)\ot\Ker(\psi^\t)^*$.
Similarly, the kernel on $(U^\d)^*\ot\iota_V(V)$ is 
$\Ker(\phi^\t)^*\ot\iota_V(V)$. Thus the kernel of 
$\phi\oth\psi$ is
\begin{equation*}
\Ker(\phi\oth\psi)=
\bigl(\iota_U(U)\ot\Ker(\psi^\t)^*\bigr)\cap
\bigl(\Ker(\phi^\t)^*\ot\iota_V(V)\bigr). 
\end{equation*}
But this is contained in $\bigl(\iota_U(U)\cap(\H\ot
\Ker(\phi^\t)^*)\bigr)\ot(V^\d)^*$. Now
$\iota_U(U)\cap\bigl(\H\ot\Ker(\phi^\t)^*\bigr)=0$,
since if $\iota(u)$ lies in $\H\ot\Ker(\phi^\t)^*$ then
$\phi(u)=0$ in $W$, so $u=0$ as $\phi$ is injective.
Thus $\Ker(\phi\oth\psi)=0$, and $\phi\oth\psi$ is injective.
\end{proof}

The following lemma is trivial to prove.

\begin{lem} Let $U,V$ be \ah-modules, and let $u\in U$ and 
$v\in V$ be nonzero. Suppose that $\alpha(u)\beta(v)=
\beta(v)\alpha(u)\in\H$ for every $\alpha\in U^\d$ and 
$\beta\in V^\d$. Define an element $u\oth v$ of\/
$\H\ot(U^\d)^*\ot(V^\d)^*$ by $(u\oth v)\cdot(\alpha\ot\beta)
=\alpha(u)\beta(v)\in\H$. Then $u\oth v$ is a nonzero
element of\/~$U\oth V$.
\label{uveltlem}
\end{lem}

The philosophy of the algebraic side of this paper
is that much algebra that works over a commutative
field such as $\R$ or $\C$ also has a close analogue
over $\H$ (or some other noncommutative algebra), when
we replace vector spaces over $\R$ or $\C$ by \ah-modules,
and tensor products of vector spaces by the quaternionic tensor product 
$\oth$. Lemmas \ref{tenslem}, \ref{injlem} and \ref{uveltlem} 
are examples of this philosophy, as they show that essential
properties of the usual tensor product also hold for
the quaternionic tensor product.

However, the quaternionic tensor product also has properties that are
very unlike the usual tensor product, which come from the 
noncommutativity of the quaternions. For example, let 
$U,V$ be $\H$-modules, and let $U'=V'=\{0\}$. Then
$U^\d=U^\t$ and $V^\d=V^\t$. Suppose that $x\in U\oth V$, and 
let $p,q\in\H$, $\alpha\in U^\d$ and $\beta\in V^\d$. Then
\begin{equation}
\begin{split}
x(\alpha\ot\beta)\ov p\,\ov q
&=x\bigl((p\cdot\alpha)\ot\beta\bigr)\ov q
=x\bigl((p\cdot\alpha)\ot(q\cdot\beta)\bigr)\\
&=x\bigl(\alpha\ot(q\cdot\beta)\bigr)\ov p
=x(\alpha\ot\beta)\ov q\,\ov p,
\end{split}
\label{xpqeq}
\end{equation}
where we use the facts that $U^\d,V^\d$ are closed under the
$\H$-action on $U^\t$, and the definition \eq{uothveq} of
$U\oth V$. Choosing $p$ and $q$ such that $\ov p\,\ov q\ne
\ov q\,\ov p$, equation \eq{xpqeq} shows that 
$x(\alpha\ot\beta)=0$. Since this holds for all $\alpha,\beta$,
$x=0$. Thus $U\oth V=\{0\}$. We have shown that the quaternionic tensor 
product of two nonzero \ah-modules can be zero. 

Suppose that $U\cong\H^k$. Then the condition in Definition 
\ref{ahmdef} implies that $\dim(U^\d)\ge k$. But 
$\dim(U')+\dim(U^\d)=4k$, so $\dim(U')\le 3k$. The
example above illustrates the general principle that if 
$\dim(U')$ is small, then quaternionic tensor products involving 
$U$ tend to be small or zero. A good rule is that the
most interesting \ah-modules $U$ are those in the 
range~$2k\le\dim(U')\le 3k$.

Here are some differences between the quaternionic and ordinary
tensor products.

\begin{itemize}
\item In contrast to Lemma \ref{uveltlem}, if $u\in U$ and 
$v\in V$, there is in general no element `$u\oth v$' in 
$U\oth V$. At best, there is a linear map from some vector 
subspace of $U\ot V$ to~$U\oth V$.
\item Suppose that $U,V$ are finite-dimensional \ah-modules,
with $U\cong\H^k$, $V\cong\H^l$. It is easy to 
show that $U\oth V\cong\H^n$, for some integer $n$ with
$0\le n\le kl$. However, $n$ can vary discontinuously
under smooth variations of $U^\d,V^\d$.
\item If we wish, we can make $U^\t,V^\t$ into \ah-modules.
However, it is not in general true that 
$U^\t\oth V^\t\cong(U\oth V)^\t$, and there is no reason for
$U^\t\oth V^\t$ and $(U\oth V)^\t$ even to have the same
dimension. For this reason, dual \ah-modules seem not to be 
a very powerful tool.
\item In contrast to Lemma \ref{injlem}, if 
$\phi:U\ra W$ and $\psi:V\ra X$ are 
{\it surjective} \ah-morphisms, then
$\phi\oth\psi:U\oth V\ra W\oth X$ does not
have to be surjective. In particular, $U\oth V$ may be
zero, but $W\oth X$ nonzero.
\end{itemize}

\subsection{Stable and semistable \ah-modules}
\label{h12}

Now two special sorts of \ah-modules will be defined, called
{\it stable} and {\it semistable} \ah-modules. Our aim in
this paper has been to develop a strong analogy between the 
theories of \ah-modules and vector spaces over a field. For
stable \ah-modules it turns out that this analogy is more 
complete than in the general case, because various important 
properties of the vector space theory hold for stable but
not for general \ah-modules. Therefore, in applications of the
theory it will often be useful to restrict to stable \ah-modules,
to exploit their better behaviour. We begin with a definition.

\begin{dfn} Let $q\in\I$ be nonzero. Define an \ah-module
$X_q$ by $X_q=\H$, and $X_q'=\{p\in\H:pq=-qp\}$. Then
$X_q'\subset\I$ and $\dim X_q'=2$. Let $\chi_q:X_q\ra\H$ 
be the identity. Then $\chi_q(X_q')\subset\I=\H'$, so
$\chi_q$ is an \ah-morphism. Suppose $U$ is any
finite-dimensional \ah-module. Then $U\oth X_q$ is an 
\ah-module, and $U\oth\H\cong U$, so there is a canonical 
\ah-morphism $\id\oth\chi_q:U\oth X_q\ra U$. 
We say that {\it $U$ is a semistable \ah-module} if

\begin{equation*}
U=\bigl\langle(\id\oth\chi_q)(U\oth X_q):
0\ne q\in\I\bigr\rangle,
\end{equation*}
that is, $U$ is generated as an $\H$-module by the images
$\Im(\id\oth\chi_q)$ for nonzero $q\in\I$. Clearly,
every \ah-module contains a unique maximal semistable \ah-submodule,
the submodule generated by the images~$\Im(\id\oth\chi_q)$. 
\label{submindef}
\end{dfn}

\begin{prop} Let\/ $q\in\I$ be nonzero. Then $X_q\oth X_q$ 
is \ah-isomorphic to $X_q$. Let\/ $U$ be a finite-dimensional 
\ah-module. Then $U\oth X_q$ is \ah-isomorphic to $nX_q$, the 
direct sum of\/ $n$ copies of\/ $X_q$, for some $n\ge 0$.
Now suppose that\/ $U$ is semistable, with\/ $\dim U=4j$ 
and\/ $\dim U'=2j+r$, for integers $j,r$. Then $U\oth X_q
\cong nX_q$, where $n\ge r$ for all nonzero $q\in\I$, and 
$n=r$ for generic $q\in\I$. Thus~$r\ge 0$.
\label{xqprop}
\end{prop}

\begin{proof} A short calculation shows that $X_q\oth X_q
\cong X_q$ as \ah-modules. Let $U$ be a finite-dimensional
\ah-module. Then $\id\oth\chi_q:U\oth X_q\ra U$ is an
injective \ah-morphism, by Lemma \ref{injlem}. It is easy
to show that $(U\oth X_q)'$ is identified with
$U'\cap(q\cdot U')$ by $\id\oth\chi_q$. Now 
$\langle 1,q\rangle\subset\H$ is a subalgebra $\C_q$ 
of $\H$ isomorphic to $\C$, and clearly $(U\oth X_q)'$
is closed under $\C_q$. Choose a basis $u_1,\dots,u_n$ of 
$(U\oth X_q)'$ over the field $\C_q$. Suppose that 
$\Sigma q_ju_j=0$ in $U\oth X_q$, for $q_1,\dots,q_n\in\H$. 
Let $p\in\H$ be nonzero, and such that $pq=-qp$. Then $\H$ 
splits as $\H=\C_q\op p\C_q$. Using this splitting, write 
$q_j=a_j+pb_j$, with $a_j,b_j\in\C_q$. 

If $\alpha\in(U\oth X_q)^\d$, then $\alpha(u_j)\in p\C_q$, 
so $\alpha(\Sigma q_ju_j)=0$ implies that
$\alpha(\Sigma a_ju_j)=\alpha(\Sigma b_ju_j)=0$.
As this hold for all $\alpha\in(U\oth X_q)^\d$, we have
$\Sigma a_ju_j=0$ and $\Sigma b_ju_j=0$. But $\{u_j\}$
is a basis over $\C_q$ and $a_j,b_j\in\C_q$. Thus 
$a_j=b_j=0$, and $q_j=0$. We have proved that 
$u_1,\dots,u_n$ are linearly independent over~$\H$.
It is now easily shown that $\H\cdot u_j\cong X_q$,
and that $U\oth X_q=\H\cdot u_1\op\cdots\op\H\cdot u_n
\cong nX_q$, as we have to prove.

Now let $U$ be semistable, with $\dim U=4j$ and 
$\dim U'=2j+r$. From above, $U'\cap(q\cdot U')=\C_q^n$.
But we have
\begin{equation}
\begin{split}
\dim\bigl(&U'\cap(q\cdot U')\bigr)
+\dim\bigl(U'+q\cdot U'\bigr)\\
&=\dim U'+\dim(q\cdot U')=4j+2r.
\end{split}
\label{uintdimeq}
\end{equation}
Since $U'+q\cdot U'\subset U$, $\dim(U'+q\cdot U')\le 4j$, and so 
\eq{uintdimeq} shows that $2n\ge 2r$, with equality if and only 
if $U=U'+q\cdot U'$. Thus $n\ge r$ for all nonzero $q$, as we have 
to prove. To complete the proposition, it is enough to show that
$U=U'+q\cdot U'$ for generic~$q\in\I$.

As $U$ is semistable, it is generated by the images
$\Im(\id\oth\chi_q)$. So suppose $U$ is generated by 
$\Im(\id\oth\chi_{q_j})$ for $j=1,\dots,k$, where 
$0\ne q_j\in\I$. Let $q\in\I$, and suppose that 
$qq_j\ne q_jq$ for $j=1,\dots,k$. This is true for generic 
$q$. Clearly $X_{q_j}'+q\cdot X_{q_j}'=X_{q_j}$, as
$qq_j\ne q_jq$. Thus $(U\oth X_{q_j})'+q\cdot(U\oth X_{q_j})'
=U\oth X_{q_j}$, as $U\oth X_{q_j}\cong nX_{q_j}$.
We deduce that $\Im(\id\oth\chi_{q_j})$ is contained in
$U'+q\cdot U'$. But $U$ is generated by the spaces
$\Im(\id\oth\chi_{q_j})$, so $U=U'+q\cdot U'$.
This completes the proof.
\end{proof}

Now we can define stable \ah-modules.

\begin{dfn} Let $U$ be a finite-dimensional \ah-module.
Then $\dim U=4j$ and $\dim U'=2j+r$, for some integers $j,r$.
Define the {\it virtual dimension of\/ $U$} to be $r$. We say that 
{\it $U$ is a stable \ah-module} if $U$ is a semistable 
\ah-module, so that $r\ge 0$, and $U\oth X_q\cong rX_q$ for each 
nonzero~$q\in\I$.
\label{mindef}
\end{dfn}

The point of this definition will become clear soon. Here are 
two propositions about stable and semistable \ah-modules.

\begin{prop} Let\/ $U$ be a stable \ah-module with
$\dim U=4j$, and $\dim U'=2j+r$ for integers $j,r$.
Let\/ $V$ be a semistable \ah-module with\/ $\dim V=4k$ 
and\/ $\dim V'=2k+s$ for integers $k,s$. Then 
$\dim(U\oth V)=4l$, where~$l=js+rk-rs$. 
\label{mstprop}
\end{prop}

\begin{proof} Regard $\H\ot U^\d\ot V^\d$ and $\H\ot(U^\d)^*
\ot(V^\d)^*$ as $\H$-modules in the obvious way. Define a bilinear 
map $\Theta:\H\ot(U^\d)^*\ot(V^\d)^*\t\H\ot U^\d\ot V^\d\ra\H$ by 
$\Theta(p\ot\alpha\ot\beta,q\ot x\ot y)=\alpha(x)\beta(y)p\ov q$
for $p,q\in\H$, $x\in U^\d$, $y\in V^\d$, $\alpha\in(U^\d)^*$
and $\beta\in(V^\d)^*$. Recall that $\iota_U(U)\ot(V^\d)^*$ and
$(U^\d)^*\ot\iota_V(V)$ are $\H$-submodules of $\H\ot(U^\d)^*\ot
(V^\d)^*$. Define a subspace $K_{U,V}$ of $\H\ot U^\d\ot V^\d$ 
by $z\in K_{U,V}$ if $\Theta(\zeta,z)=0$ whenever $\zeta\in
\iota_U(U)\ot(V^\d)^*$ or $\zeta\in(U^\d)^*\ot\iota_V(V)$.
Then $K_{U,V}$ is an $\H$-submodule of~$\H\ot U^\d\ot V^\d$.

Now $\iota_U(U)\ot(V^\d)^*+(U^\d)^*\ot\iota_V(V)$ is an 
$\H$-submodule of $\H\ot(U^\d)^*\ot(V^\d)^*$, and clearly
\begin{align*}
\dim K_{U,V}&+\dim\bigl(\iota_U(U)\ot(V^\d)^*+(U^\d)^*\ot\iota_V(V)\bigr)\\
&=\dim\H\ot(U^\d)^*\ot(V^\d)^*=4(2j-r)(2k-s).
\end{align*}
But $U\oth V=\bigl(\iota_U(U)\ot(V^\d)^*\bigr)
\cap\bigl((U^\d)^*\ot\iota_V(V)\bigr)$, and thus
\begin{align*}
\dim(U\oth V)=
&\dim\bigl(\iota_U(U)\ot(V^\d)^*\bigr)
+\dim\bigl((U^\d)^*\ot\iota_V(V)\bigr)\\
&-\dim\bigl(\iota_U(U)\ot(V^\d)^*+(U^\d)^*\ot\iota_V(V)\bigr),
\end{align*}
so that $\dim(U\oth V)=4j(2k\!-\!s)\!+\!(2j\!-\!r)4k\!-
\!\bigl\{4(2j\!-\!r)(2k\!-\!s)\!-\!\dim K_{U,V}\bigr\}
=4l+\dim K_{U,V}$, where $l=js+rk-rs$. Therefore 
$\dim(U\oth V)=4l$ if and only if $K_{U,V}=\{0\}$.

Suppose that $W$ is an \ah-module, and $\phi:W\ra V$ is 
an \ah-morphism. Then $\phi^\t:V^\d\ra W^\d$, so that
$\id\ot\phi^\t:\H\ot U^\d\ot V^\d\ra\H\ot U^\d\ot W^\d$.
We have $K_{U,V}\subset\H\ot U^\d\ot V^\d$ and $K_{U,W}
\subset\H\ot U^\d\ot W^\d$. It is easy to show that
$(\id\ot\phi^\t)(K_{U,V})\subset K_{U,W}$. Let
$0\ne q\in\I$, and put $W=V\oth X_q$, and $\phi=\id\oth\chi_q$. 
In this case $W\cong nX_q$. The argument above shows that
$K_{U,X_q}=\{0\}$ if and only if $\dim(U\oth X_q)=4r$.
But by this holds by Definition \ref{mindef}, as $U$ is stable. 

Thus $K_{U,X_q}=\{0\}$, and $K_{U,W}=\{0\}$ as
$W\cong nX_q$. It follows that $(\id\ot\phi^\t)(K_{U,V})=
\{0\}$, so $K_{U,V}\subset\H\ot U^\d\ot\Ker\phi^\t$.
Now $V$ is semistable. Therefore $V$ is generated by
submodules $\phi(W)$ of the above type, and the 
intersection of the subspaces $\Ker\phi^\t\subset V^\d$
for all nonzero $q$, must be zero. So $K_{U,V}\subset
\H\ot U^\d\ot\{0\}$, giving $K_{U,V}=\{0\}$, and 
$\dim(U\oth V)=4l$ from above, which completes the proof.
\end{proof}

\begin{prop} Let\/ $U$ be a stable \ah-module, and\/
$V$ a semistable \ah-module. Then $U\oth V$ is semistable.
\label{sstprop}
\end{prop}

\begin{proof} Let $\dim U=4j$, $\dim V=4k$, $\dim U'=2j+r$ 
and $\dim V'=2k+s$. Then Proposition \ref{mstprop} shows
that $\dim(U\oth V)=4l$, where $l=js+rk-rs$. Let $W\subset 
U\oth V$ be the \ah-submodule of $U\oth V$ generated by
the images $\Im(\id\oth\chi_q)$, where $\id\oth\chi_q:
U\oth V\oth X_q\ra U\oth V$ and $0\ne q\in\I$. Then $W$ is the 
maximal semistable \ah-submodule of $U\oth V$. We shall prove 
the proposition by explicitly constructing $l$ elements of $W$, 
that are linearly independent over $\H$. This will imply that 
$\dim W\ge 4l$. Since $W\subset U\oth V$ and $\dim(U\oth V)=4l$, 
we see that $W=U\oth V$, so $U\oth V$ is semistable.

Here is some new notation. Let $0\ne q\in\I$, and define 
$U_q=\bigl\{u\in U:\alpha(u)\in\langle 1,q\rangle$ for all
$\alpha\in U^\d\bigr\}$. Similarly define $V_q,(U\oth V)_q$. 
It can be shown that $\H\cdot U_q=\Im(\id\oth\chi_q)$, and
similarly for $V_q,(U\oth V)_q$. Thus $(U\oth V)_q\subset W$.
Since $V$ is semistable, we can choose nonzero elements 
$q_1,\dots,q_k\in\I$ and $v_1,\dots,v_k\in V$, such that 
$v_a\in V_{q_a}$ and $(v_1,\dots,v_k)$ is a basis for $V$ over 
$\H$. Since $U$ is stable, $U_q\cong\C^r$ for each $0\ne q\in\I$. 
Therefore for each $a=1,\dots,k$ we may choose elements 
$u_{a1},\dots,u_{ar}$ of $U$, such that $u_{ab}\in U_{q_a}$ and 
$u_{a1},\dots,u_{ar}$ are linearly independent over $\H$ in~$U$. 

As $U$ is stable, it is not difficult to see that there
are nonzero elements $p_1,\dots,p_{j-r}\in\I$ and 
$u_1,\dots,u_{j-r}\in U$, such that $u_c\in U_{p_c}$, and
for each $a=1,\dots,k$, the set $u_{a1},\dots,u_{ar},u_1,
\dots,u_{j-r}$ is linearly independent over $\H$. This is
just a matter of picking {\it generic} elements $p_c$ and
$u_c$, and showing that generically, linear independence
holds. 

We shall also need another property. Define $F\subset U^\d$
by $F=\{\alpha\in U^\d:\alpha(u_c)=0$ for $c=1,\dots,j-r\}$.
Now as $\alpha(u_c)\in\langle 1,q_c\rangle$ for each
$\alpha\in U^\d$, the codimension of $F$ in $U^\d$ is
at most $2(j-r)$. Since $\dim U^\d=2j-r$, this gives 
$\dim F\ge r$. The second property we need is that for
each $a=1,\dots,k$, if $u\in\langle u_{a1},\dots,u_{ar}
\rangle_\H$ and $\alpha(u)=0$ for all $\alpha\in F$, then
$u=0$. Here $\langle\,,\,\rangle_\H$ means the linear span 
over $\H$. Again, it can be shown that for generic choice
of $p_c,v_c$, this property holds.

Now for $a=1,\dots,k$ and $b=1,\dots,r$, $u_{ab}\in U_{q_a}$
and $v_a\in V_{q_a}$. Lemma \ref{uveltlem} gives an element 
$u_{ab}\oth v_a$ in $U\oth V$. Moreover $u_{ab}\oth v_a\in
(U\oth V)_{q_a}$. Thus $u_{ab}\oth v_a\in W$. Therefore, we 
have made $kr$ elements $u_{ab}\oth v_a$ of $W$. Similarly, 
for $c=1,\dots,j\!-\!r$ and $d=1,\dots,s$, the elements 
$u_c\oth v_{cd}$ exist in $W$, which gives a further $(j\!-\!r)s$ 
elements of $W$. Since $kr+(j\!-\!r)s=l$, we have constructed $l$ 
explicit elements $u_{ab}\oth v_a$ and $u_c\oth v_{cd}$ of~$W$.

Suppose that $\Sigma_{a,b}\,x_{ab}u_{ab}\oth v_a
+\Sigma_{c,d}\,y_{cd}u_c\oth v_{cd}=0$ in $U\oth V$, where
$x_{ab}\in\H$ for $a=1,\dots,k$, $b=1,\dots,r$ and $y_{cd}\in\H$ 
for $c=1,\dots,j-r$ and $d=1,\dots,s$. We shall show that
$x_{ab}=y_{cd}=0$. Now this equation implies that
\begin{equation}
\sum_{a,b}x_{ab}\alpha(u_{ab})\beta(v_a)
+\sum_{c,d}y_{cd}\alpha(u_c)\beta(v_{cd})=0
\label{xyuvsumeq}
\end{equation}
for all $\alpha\in U^\d$ and $\beta\in V^\d$, by Lemma
\ref{uveltlem}. Let $\alpha\in F$. Then $\alpha(u_c)=0$,
so $\Sigma_{a,b}\,x_{ab}\alpha(u_{ab})\beta(v_a)=0$. As this
holds for all $\beta\in V^\d$ and the $v_a$ are linearly 
independent over $\H$, it follows that $\Sigma_b\,x_{ab}
\alpha(u_{ab})=0$ for all $a$ and all~$\alpha\in F$.

By the property of $F$ assumed above, this implies that 
$x_{ab}=0$ for all $a,b$. It is now easy to show that $y_{cd}=0$ 
for all $c,d$. Thus the $l$ elements $u_{ab}\oth v_a$ and 
$u_c\oth v_{cd}$ of $W$ are linearly independent over $\H$, 
so $\dim W\ge 4l$. But $\dim(U\oth V)=4l$ by Proposition
\ref{mstprop}. So $U\oth V=W$, and $U\oth V$ is 
semistable. This finishes the proof.
\end{proof}

We can now prove the main result of this section.

\begin{thm} Let\/ $U$ be a stable \ah-module and\/ $V$ be a 
semistable \ah-module with
\begin{equation*}
\dim U=4j, \; \dim U'=2j+r, \;
\dim V=4k\;\text{and}\; \dim V'=2k+s 
\end{equation*}
for integers $j,k,r$ and $s$. Then $U\oth V$ is a semistable \ah-module 
with $\dim(U\oth V)=4l$ and\/ $\dim(U\oth V)'=2l+t$, where $l=js+rk-rs$ 
and\/ $t=rs$. If\/ $V$ is stable, then $U\oth V$ is stable.
\label{minsubtensthm}
\end{thm}

\begin{proof} Let $l=js+rk-rs$ and $t=rs$. Proposition \ref{mstprop} 
shows that $\dim(U\oth V)=4l$. As $U$ is stable, $U\oth X_q\cong rX_q$ 
for nonzero $q\in\I$. As $V$ is semistable, Proposition \ref{xqprop} 
shows that $V\oth X_q\cong sX_q$ for generic $q\in\I$. Thus 
$U\oth V\oth X_q\cong rsX_q=tX_q$ for generic $q\in\I$. Also,
Proposition \ref{sstprop} shows that $U\oth V$ is semistable.
Combining these two facts with Proposition \ref{xqprop}, we
see that $\dim(U\oth V)'=2l+t$, as we have to prove.

It remains to show that if $V$ is stable, then $U\oth V$ is stable.
Suppose $V$ is stable. Then $V\oth X_q\cong sX_q$ for all nonzero
$q\in\I$, so $U\oth V\oth X_q=tX_q$ for all nonzero $q\in\I$.
As $U\oth V$ is semistable, it is stable, by Definition \ref{mindef}.
This completes the proof.
\end{proof}

Theorem \ref{minsubtensthm} shows that if $U$ is stable and 
$V$ semistable, then the virtual dimension of $U\oth V$ is the 
product of the virtual dimensions of $U$ and $V$. Thus the 
virtual dimension is a good analogue of the dimension of a vector 
space, as it multiplies under $\oth$. Note also that 
$(j-r)/r+(k-s)/s=(l-t)/t$, so that the nonnegative function
$U\mapsto(j-r)/r$ behaves additively under~$\oth$.

Next we will show that generic \ah-modules $(U,U')$ with
positive virtual dimension are stable. Thus there are many 
stable \ah-modules.

\begin{lem} Let\/ $j,r$ be integers with\/ $0<r\le j$. Let\/ 
$U=\H^j$, and let\/ $U'$ be a real vector subspace of\/ $U$ 
with\/ $\dim U'=2j+r$. Then for generic subspaces $U'$,
$(U,U')$ is a stable \ah-module.
\end{lem}

\begin{proof} Let $G$ be the Grassmannian of real $(2j\!+\!r)$-planes
in $U\cong\R^{4j}$. Then $U'\in G$, and $\dim G=4j^2-r^2$. 
The condition for $(U,U')$ to be an \ah-module is that
$\H\cdot U^\d=U^\t$. A calculation shows that this fails
for a subset of $G$ of codimension $4(j-r+1)$. Thus for generic
$U'\in G$, $(U,U')$ is an \ah-module. Suppose $(U,U')$
is an \ah-module. Let $W$ be the maximal semistable \ah-submodule 
in $U$. Then $W\cong\H^k$, for some $k$ with $r\le k\le j$. A 
calculation shows that for given $k$, the subset of $G$ with
$W\cong\H^k$ is of codimension $2(j-k)r$. Thus for generic
$U'\in G$, $W=\H^j=U$, and $U$ is semistable.

Suppose $(U,U')$ is semistable. Let $0\ne q\in\I$. Then 
$U\oth X_q\cong rX_q$ if and only if $U'\cap q\cdot U'=\R^{2r}$. 
A computation shows that this fails for a subset of $G$ of 
codimension $2r+2$. Thus the condition $U'\cap q\cdot U'=
\R^{2r}$ for all nonzero $q\in\I$ fails for a subset of $G$ of 
codimension at most $2r$, since this subset is the union of a 
2-dimensional family of $2r+2$-codimensional subsets, the 
2-dimensional family being ${\cal S}^2$, the unit sphere in $\I$. 
Therefore for generic $U'\in G$, $U\oth X_q\cong rX_q$ for all 
nonzero $q\in\I$, and $U$ is stable.
\end{proof}

We leave the proof of this proposition to the reader, as a 
(difficult) exercise.

\begin{prop} Let\/ $U$ be a stable \ah-module, with\/ 
$\dim U=4j$ and\/ $\dim U'=2j+r$. Let\/ $n$ be a positive
integer. Then $\Sh^nU$ and $\Lambh^nU$ are stable
\ah-modules, with\/ $\dim(\Sh^nU)=4k$, $\dim(\Sh^nU)'=2k+s$,
$\dim(\Lambh^nU)=4l$ and\/ $\dim(\Lambh^nU)'=2l+t$, where
\begin{alignat*}{2}
k&=(j-r)\binom{r+n-1}{n-1}+\binom{r+n-1}{n},&\quad
s&=\binom{r+n-1}{n},\\
l&=(j-r)\binom{r-1}{n-1}+\binom{r}{n}\qquad\text{and}&\quad
t&=\binom{r}{n}.
\end{alignat*}
\label{symmantiprop}
\end{prop}

\subsection{Stable \ah-modules and exact sequences}
\label{h13}

Recall that if $U,V,W$ are vector spaces and $\phi:U\ra V$,
$\psi:V\ra W$ are linear maps, then we say that the sequence
$U{\br\phi\over\ra}V{\br\psi\over\ra}W$ is {\it exact at $V$} 
if $\Im\phi=\Ker\psi$. Here is the analogue of this for \ah-modules.

\begin{dfn} Let $U,V,W$ be \ah-modules, and let $\phi:U\ra V$ 
and $\psi:V\ra W$ be \ah-morphisms. We say that 
{\it the sequence $U{\br\phi\over\ra}V{\br\psi\over\ra}W$ is 
\ah-exact at\/ $V$} if the sequence $U{\br\phi\over\ra}V
{\br\psi\over\ra}W$ is exact at $V$, and the sequence 
$U'{\br\phi\over\ra}V'{\br\psi\over\ra}W'$ is exact at 
$V'$. We say that a sequence of \ah-morphisms is \ah-exact 
if it is \ah-exact at every term.
\end{dfn}

Here is an example of some bad behaviour of the theory.

\begin{ex} Define $U=\H$ and $U'=\{0\}$. Define 
$V=\H^2$ and $V'=\bigl\langle(1,i_1),(1,i_2)\bigr\rangle$.
Define $W=\H$ and $W'=\langle i_1,i_2\rangle$. Then
$U,V,W$ are \ah-modules. Define linear maps 
$\phi:U\ra V$ by $\phi(q)=(q,0)$, and $\psi:V\ra W$ by
$\psi\bigl((p,q)\bigr)=q$. Then $\phi,\psi$ are \ah-morphisms,
and the sequence $0\ra U{\br\phi\over\ra}V
{\br\psi\over\ra} W\ra 0$ is \ah-exact.

Now set $Z$ to be the \ah-module $W$. A short calculation 
shows that $U\oth Z=\{0\}$, $V\oth Z=\{0\}$, but
$W\oth Z\cong W$. It follows that the sequence
\begin{equation*}
0\ra 
U\oth Z\;{\br\phi\smalloth\id\over\longra}\;
V\oth Z\;{\br\psi\smalloth\id\over\longra}\;
W\oth Z\ra 0
\end{equation*}
is {\it not} \ah-exact at $W\oth Z$. This contrasts with 
the behaviour of exact sequences of real vector spaces 
under the tensor product.
\end{ex}

The following proposition gives a clearer idea of what is happening.

\begin{prop} Suppose that\/ $U,V,W$ and\/ $Z$ are finite-dimensional
\ah-modules, that\/ $\phi:U\ra V$ and\/ $\psi:V\ra W$ are 
\ah-morphisms, and that the sequence $0\ra U{\br\phi\over\ra}V
{\br\psi\over\ra}W\ra 0$ is \ah-exact. Then the sequence 
\begin{equation}
0\ra W^\d\,{\br\psi^\t\over\longra}\,V^\d\,
{\br\phi^\t\over\longra}\,U^\d\ra 0
\label{dagseq}
\end{equation} 
is exact, and the sequence
\begin{equation}
0\ra U\oth Z\;{\br\phi\smalloth\id\over\longra}\;
V\oth Z\;{\br\psi\smalloth\id\over\longra}\;
W\oth Z\ra 0
\label{tensseq}
\end{equation}
is \ah-exact at\/ $U\oth Z$ and\/ $V\oth Z$, but it need not
be \ah-exact at\/~$W\oth Z$.
\label{exactprop}
\end{prop}

\begin{proof} Since $\psi:V\ra W$ is surjective, $\psi^\t:
W^\t\ra V^\t$ is injective, and thus $\psi^\t:W^\d\ra V^\d$ 
is injective. Now $\phi$ induces a map $\phi:U/U'\ra V/V'$.
Suppose that $u+U'$ lies in the kernel of this map. Then
$\phi(u)\in V'$. By exactness, $\psi(\phi(u))=0$. As
the sequence $U'\ra V'\ra W'$ is also exact,
$\phi(u)=\phi(u')$ for some $u'\in U'$. But $\phi$
is injective, so $u\in U'$. Thus the map
$\phi:U/U'\ra V/V'$ is injective. It follows that
the map $\phi^*:(V/V')^*\ra(U/U')^*$ is surjective.

Since $U,V$ are finite-dimensional, it follows that 
$U^\d\cong(U/U')^*$ and $V^\d\cong(V/V')^*$. Under
these isomorphisms $\phi^*$ is identified with $\phi^\t$.
Thus $\phi^\t:V^\d\ra U^\d$ is surjective. Using exactness
we see that $\dim V=\dim U+\dim W$ and $\dim V'=\dim U'
+\dim W'$. By subtraction we find that $\dim V^\d=
\dim U^\d+\dim W^\d$. But we have already shown that
$\psi^\t:W^\d\ra V^\d$ is injective, and $\phi^\t:V^\d\ra U^\d$
is surjective. Using these facts, we see that the sequence
\eq{dagseq} is exact, as we have to prove.

Now $\phi:U\ra V$ is injective, and clearly $\id:Z\ra Z$ is 
injective. Thus Lemma \ref{injlem} shows that 
$\phi\oth\id:U\oth Z\ra V\oth Z$ is injective. So the sequence 
\eq{tensseq} is \ah-exact at $U\oth Z$, as we have to prove. 
Suppose that $x\in\Ker(\psi\oth\id)$. It is easy to show,
using exactness, injectivity, and the definition of $\oth$,
that $x\in\Im(\phi\oth\id)$. Thus $\Ker(\psi\oth\id)=
\Im(\phi\oth\id)$, and the sequence \eq{tensseq} is 
exact at $V\oth Z$, as we have to prove. 

Recall from Definition \ref{qctpdef} that as $U,Z$ are 
finite-dimensional, $\lambda_{U,Z}$ is surjective. Suppose 
that $a\in(U\oth Z)^\d$. Then $a=\lambda_{U,Z}(b)$ for some 
$b\in U^\d\ot Z^\d$. The map $\phi^\t:V^\d\ra U^\d$ is 
surjective from above, and thus $b=(\phi^\t\ot\id)(c)$ for 
some $c\in V^\d\ot Z^\d$. Thus $a=\lambda_{U,Z}\circ
(\phi^\t\ot\id)(c)$. But it is easy to see that $\lambda_{U,Z}
\circ(\phi^\t\ot\id)=(\phi\oth\id)^\t\circ\lambda_{V,Z}$, as 
maps $V^\d\ot Z^\d\ra(U\oth Z)^\d$. Thus if $a\in(U\oth Z)^\d$, 
then $a=(\phi\oth\id)^\t(d)$, for $d=\lambda_{V,Z}(c)\in
(V\oth Z)^\d$. Therefore $(\phi\oth\id)^\t:(V\oth Z)^\d
\ra(U\oth Z)^\d$ is surjective.

Let $x\in(V\oth Z)'$, and suppose $(\psi\oth\id)(x)=0$.
As \eq{tensseq} is exact, $x=(\phi\oth\id)(y)$, and as
$\phi\oth\id$ is injective, $y$ is unique. We shall show
that $y\in(U\oth Z)'$. It is enough to show that for
each $\alpha\in(U\oth Z)^\d$, $\alpha(y)\in\I$. Since
$(\phi\oth\id)^\t:(V\oth Z)^\d\ra(U\oth Z)^\d$ is 
surjective, $\alpha=(\phi\oth\id)^\t(\beta)$ for
$\beta\in(V\oth Z)^\d$. Then $\alpha(y)=\beta(x)$.
But $x\in(V\oth Z)'$ and $\beta\in(V\oth Z)^\d$,
so $\beta(x)\in\I$. Thus $\alpha(y)\in\I$, and
$y\in(U\oth Z)'$. It follows that the sequence
$(U\oth Z)'{\br\phi\smalloth\id\over\longra}
(V\oth Z)'{\br\phi\smalloth\id\over\longra}
(W\oth Z)'$ is exact at $(V\oth Z)'$. Therefore
\eq{tensseq} is \ah-exact at $V\oth Z$, as we have 
to prove. The example above shows that the sequence need 
not be \ah-exact at $W\oth Z$, and the proposition is finished.
\end{proof}

In the language of category theory, Proposition \ref{exactprop}
shows that when $Z$ is a finite-dimensional \ah-module, 
the operation $\oth Z$ is a {\it left-exact functor}, but 
may not be a {\it right-exact functor}. However, the next 
proposition shows that right-exactness does hold for stable and 
semistable \ah-modules.

\begin{prop} Suppose that\/ $U,V,W$ are \ah-modules, that\/ $U$
and\/ $W$ are stable, that\/ $\phi:U\ra V$ and\/ $\psi:V\ra W$ 
are \ah-morphisms, and that the sequence $0\ra U{\br\phi\over\ra}V
{\br\psi\over\ra}W\ra 0$ is \ah-exact. Let $Z$ be a semistable 
\ah-module. Then the following sequence is \ah-exact:
\begin{equation}
0\ra U\oth Z\;{\br\phi\smalloth\id\over\longra}\;
V\oth Z\;{\br\psi\smalloth\id\over\longra}\;
W\oth Z\ra 0.
\label{newtensseq}
\end{equation}
\label{minexactprop}
\end{prop} 

\begin{proof} Let $\dim U=4j$, $\dim V=4k$, $\dim W=4l$,
$\dim U'=2j+r$, $\dim V'=2k+s$, $\dim W'=2l+t$,
$\dim Z=4a$ and $\dim Z'=2a+b$. Then by Theorem
\ref{minsubtensthm} we have
\begin{equation}
\begin{alignedat}{2}
\dim (U\oth Z)&=4(jb+ra-rb),&\;
\dim (W\oth Z)&=4(lb+ta-tb),\\
\dim (U\oth Z)'&=2jb+2ra-rb,&\;
\dim (W\oth Z)'&=2lb+2ta-tb.
\end{alignedat}
\label{uwdimeq}
\end{equation}
Theorem \ref{minsubtensthm} calculates the dimensions of a quaternionic
tensor product of stable and semistable \ah-modules. Examining the 
proof, it is easy to see that these dimensions are actually {\it 
lower bounds} for the dimensions of a quaternionic tensor product of general
\ah-modules. Therefore 
\begin{equation}
\label{vdimeq}
\!\!\!\!\!\!\!
\dim (V\oth Z)\ge 4(kb+sa-sb)\;\text{and}\;
\dim (V\oth Z)'\ge 2kb+2sa-sb.
\end{equation}

By Proposition \ref{exactprop}, the sequence \eq{newtensseq} is
\ah-exact at $U\oth Z$ and $V\oth Z$. The only way \ah-exactness at
$W\oth Z$ can fail is for $\psi\oth\id:V\oth Z\ra W\oth Z$ or 
$\psi\oth\id:(V\oth Z)'\ra(W\oth Z)'$ not to be surjective.
We deduce that
\begin{equation}
\begin{split}
\dim(U\oth Z)+\dim(W\oth Z)&\ge\dim(V\oth Z),\\
\dim(U\oth Z)'+\dim(W\oth Z)'&\ge\dim(V\oth Z)'.
\end{split}
\label{exdimeq}
\end{equation}
But $0\ra U{\br\phi\over\ra}V{\br\psi\over\ra}W\ra 0$ is \ah-exact, 
so that $k=j+l$ and $s=r+t$. Combining \eq{uwdimeq}, \eq{vdimeq}
and \eq{exdimeq}, we see that equality holds in \eq{vdimeq}
and \eq{exdimeq}, because the inequalities go opposite ways.
Counting dimensions, $\psi\oth\id:V\oth Z\ra W\oth Z$ and 
$\psi\oth\id:(V\oth Z)'\ra(W\oth Z)'$ must be surjective.
Thus by definition, \eq{newtensseq} is \ah-exact at $W\oth Z$,
and the proposition is proved.
\end{proof}

The proposition gives one reason why it is convenient, in many
situations, to work with stable \ah-modules rather than
general \ah-modules.

\begin{prop} Suppose $U,V,W$ are \ah-modules, $\phi:U\ra V$ and\/ 
$\psi:V\ra W$ are \ah-morphisms, and that the sequence 
$0\ra U{\br\phi\over\ra}V{\br\psi\over\ra}W\ra 0$ is \ah-exact. 
If\/ $U$ and\/ $W$ are stable \ah-modules, then $V$ is a stable 
\ah-module.
\label{midminprop}
\end{prop} 

\begin{proof} Let $q\in\I$ be nonzero, and apply Proposition 
\ref{minexactprop} with $Z=X_q$, the semistable \ah-module defined
in \S\ref{h12}. Let $\dim V=4k$ and $\dim V'=2k+s$ as in the
proof of the proposition. Since equality holds in \eq{vdimeq},
$\dim(V\oth X_q)=4s$, so that $V\oth X_q\cong sX_q$. But this is
the main condition for $V$ to be stable. Thus, it remains only to
show that $V$ is semistable.

Let $S$ be the maximal semistable \ah-submodule of $V$. As
$U$ is semistable, $\phi(U)\subset S$. Also, from above the map
$\psi\oth\id:V\oth X_q\ra W\oth X_q$ is surjective. As $W$ is 
semistable, we deduce that $\psi(S)=W$. But $0\ra U{\br\phi\over\ra}
V{\br\psi\over\ra}W\ra 0$ is \ah-exact, so $\phi(U)\subset S$ and
$\psi(S)=W$ imply that $S=V$. Therefore $V$ is semistable, so $V$ 
is stable.
\end{proof}

\section{Algebraic structures over the quaternions}

In this chapter, the machinery of Chapter 1 will be
used to define quaternionic analogues of various algebraic concepts.
We shall only discuss those structures we shall need for
our study of hypercomplex geometry, but the reader will soon see how
the process works in general. First, in \S\ref{h21} we define
H-algebras, the quaternionic version of commutative algebras, and
modules over H-algebras. Then in \S\ref{h22} we define HL-algebras 
and HP-algebras, the analogues of Lie algebras and Poisson algebras. 
Section \ref{h23} is about filtered and graded H-algebras, and
\S\ref{h24} considers free and finitely-generated H-algebras.

\subsection{H-algebras and modules}
\label{h21}

Now we will define the quaternionic version of a commutative algebra,
which we shall call an H-algebra, and also modules over H-algebras.
In Chapter 3 we shall see that the q-holomorphic functions
on a hypercomplex manifold form an H-algebra. Here are two axioms.

\begin{ax}{A} $A$ is an \ah-module.
\item There is an \ah-morphism $\mu_A:A\oth A\ra A$,
called the {\it multiplication map}.
\item $\Lambh^2A\subset\Ker\mu_A$. Thus $\mu_A$ is {\it commutative}.
\item The \ah-morphisms $\mu_A:A\oth A\ra A$ and
$\id:A\ra A$ combine to give \ah-morphisms
$\mu_A\oth\id$ and $\id\oth\mu_A:A\oth A\oth A\ra A\oth A$.
Composing with $\mu_A$ gives \ah-morphisms 
$\mu_A\circ(\mu_A\oth\id)$ and $\mu_A\circ(\id\oth\mu_A):A\oth 
A\oth A\ra A$. Then $\mu_A\circ(\mu_A\oth\id)=
\mu_A\circ(\id\oth\mu_A)$. This is {\it associativity of
multiplication}.
\item An element $1\in A$ called the {\it identity} is given, 
with $1\notin A'$ and~$\I\cdot 1\subset A'$. 
\item Part $(v)$ implies that if $\alpha\in A^\d$ then 
$\alpha(1)\in\R$. Thus for each $a\in A$, $1\oth a$ and 
$a\oth 1\in A\oth A$ by Lemma \ref{uveltlem}. Then 
$\mu_A(1\oth a)=\mu_A(a\oth 1)=a$ for each $a\in A$.
Thus {\it $1$ is a multiplicative identity}.
\end{ax}

\begin{ax}M Let $U$ be an \ah-module.
\item There is an \ah-morphism $\mu_U:A\oth U\ra U$ 
called the {\it module multiplication map}.
\item The maps $\mu_A$ and $\mu_U$ define \ah-morphisms
$\mu_A\oth\id$ and $\id\oth\mu_U:A\oth A\oth U\ra A\oth U$.
Composing with $\mu_U$ gives \ah-morphisms
$\mu_U\circ(\mu_A\oth\id)$ and $\mu_U\circ(\id\oth\mu_U): 
A\oth A\oth U\ra U$. Then $\mu_U\circ(\mu_A\oth\id)=
\mu_U\circ(\id\oth\mu_U)$. This is {\it associativity of
module multiplication}.
\item For $u\in U$, $1\oth u\in A\oth U$ by Lemma
\ref{uveltlem}. Then $\mu_U(1\oth u)=u$ for all $u\in U$.
Thus 1 acts as an identity on~$U$.
\end{ax}

Now we can define H-algebras and modules over them. Here H-algebra 
stands for {\it Hamilton algebra} (but my wife calls them Happy algebras). 

\begin{dfn} 
{\setlength{\parskip}{-8pt}
\begin{itemize}
\setlength{\itemsep}{-4pt}
\setlength{\parsep}{0pt}
\item An {\it H-algebra} satisfies Axiom A. 
\item A {\it noncommutative H-algebra} satisfies Axiom A, 
except part~$(iii)$.
\item Let $A$ be an H-algebra. A {\it module $U$ over $A$}
satisfies Axiom M.
\end{itemize}}
\label{halgdef}
\end{dfn}

An H-algebra is basically a commutative algebra over
the skew field $\H$. This is a strange idea:
how can the algebra commute when the field does not?
The obvious answer is that the algebra is only a partial
algebra, and multiplication is only allowed when the elements
commute. I'm not sure if this is the full story, though. In 
this paper our principal interest is in commutative H-algebras,
but noncommutative H-algebras also exist.

The associative axiom A$(iv)$ gives a good example of the
issues involved in finding quaternionic analogues of algebraic structures.
The usual formulation is that $(ab)c=a(bc)$ for all $a,b,c\in A$. 
This is not suitable for the quaternionic case, as not all elements in $A$ 
can be multiplied, so we rewrite the axiom in terms of linear maps 
of tensor products, and the quaternionic analogue becomes clear.
Finally we define morphisms of H-algebras.

\begin{dfn} Let $A,B$ be H-algebras, and let $\phi:A\ra B$ be 
an \ah-morphism. Write $1_A,1_B$ for the identities in $A,B$ 
respectively. We say $\phi$ is a {\it morphism of H-algebras} or
an {\it H-algebra morphism} if $\phi(1_A)=1_B$ and 
$\mu_B\circ(\phi\oth\phi)=\phi\circ\mu_A$ as 
\ah-morphisms~$A\oth A\ra B$. 
\label{halgmordef}
\end{dfn}

\subsection{H-algebras and Poisson brackets}
\label{h22}

Recall the idea of a {\it Lie bracket} on a vector space. 
A real algebra may be equipped with a Lie bracket satisfying 
certain conditions, and in this case the Lie bracket is 
called a {\it Poisson bracket}, and the algebra is called a 
{\it Poisson algebra}. Poisson algebras are studied in \cite{BV}. 
In this section we define one possible analogue of these 
concepts in our theory of quaternionic algebra. The analogue of a 
Poisson algebra will be called an {\it HP-algebra}. We begin by 
defining a special \ah-module.

\begin{dfn} Define $Y\subset\H^3$ by $Y=\bigl\{(q_1,q_2,q_3):
q_1i_1+q_2i_2+q_3i_3=0\bigr\}$. Then $Y\cong\H^2$ is an 
$\H$-module. Define $Y'\subset Y$ by $Y'=\bigl\{(q_1,q_2,q_3)
\in Y:q_j\in\I\bigr\}$. Then $\dim Y'=5$ and $\dim Y^\d=3$. 
Thus $\dim Y=4j$ and $\dim Y'=2j+r$, where $j=2$ and $r=1$. 
Define a map $\nu:Y\ra\H$ by $\nu\bigl((q_1,q_2,q_3)\bigr)
=i_1q_1+i_2q_2+i_3q_3$. Then $\Im\nu=\I$, and $\Ker\nu=Y'$.
But $Y/Y'\cong(Y^\d)^*$, so that $\nu$ induces an isomorphism 
$\nu:(Y^\d)^*\ra\I$. Since $\I\cong\I^*$, we have $(Y^\d)^*\cong
\I\cong Y^\d$. It is easy to see that $Y$ is a stable 
\ah-module. Thus $Y\oth Y$ satisfies $\dim(Y\oth Y)=12$ and 
$\dim(Y\oth Y)'=7$, by Theorem \ref{minsubtensthm}. Proposition 
\ref{symmantiprop} then shows that $\Sh^2Y=Y\oth Y$ and~$\Lambh^2Y=\{0\}$.
\label{ymoddef}
\end{dfn}

Let $A$ be an \ah-module. Here is the axiom for a Lie bracket on~$A$.

\begin{ax}{P1} There is an \ah-morphism $\xi_A:A\oth A\ra
A\oth Y$ called the {\it Lie bracket} or {\it Poisson bracket}, 
where $Y$ is the \ah-module of Definition~\ref{ymoddef}.
\item $\Sh^2A\subset\Ker\xi_A$. Thus $\xi_A$ is 
{\it antisymmetric}.
\item There are \ah-morphisms $\id\oth\xi_A:A\oth A\oth A\ra
A\oth A\oth Y$ and $\xi_A\oth\id:A\oth A\oth Y\!\ra\!A\oth Y\oth Y$.
Composing gives an \ah-morphism $(\xi_A\oth\id)\circ(\id\oth\xi_A):
A\oth A\oth A\ra A\oth Y\oth Y$. Then \newline
$\Lambh^3A\subset\Ker\bigl((\xi_A\oth\id)\circ(\id\oth\xi_A)\bigr)$. 
This is {\it the Jacobi identity for\/~$\xi_A$}.
\end{ax}

For the next axiom, let $A$ be an H-algebra.

\begin{ax}{P2} If $a\in A$, we have $1\oth a\in A\oth A$.
Then~$\xi_A(1\oth a)=0$.
\item There are \ah-morphisms $\id\oth\xi_A:A\oth A\oth A\ra
A\oth A\oth Y$ and $\mu_A\oth\id:A\oth A\oth Y\ra A\oth Y$.
Composing gives an \ah-morphism $(\mu_A\oth\id)\circ(\id\oth\xi_A):
A\oth A\oth A\ra A\oth Y$. Similarly, there are \ah-morphisms 
$\mu_A\oth\id:A\oth A\oth A\ra A\oth A$ and $\xi_A:A\oth A\ra 
A\oth Y$. Composing gives an \ah-morphism $\xi_A\circ(\mu_A\oth\id):
A\oth A\oth A\!\ra\!A\oth Y$. Then $\xi_A\circ(\mu_A\oth\id)=
2(\mu_A\oth\id)\circ(\id\oth\xi_A)$ on $\Sh^2A\oth A$. This is
{\it the derivation property.}
\end{ax}

Now we can define HL-algebras and HP-algebras. Here HL-algebra
stands for {\it Hamilton-Lie algebra}, and HP-algebra stands
for {\it Hamilton-Poisson algebra} (but my wife calls these 
Happy Fish algebras).

\begin{dfn} 
{\setlength{\parskip}{-8pt}
\begin{itemize}
\setlength{\itemsep}{-4pt}
\setlength{\parsep}{0pt}
\item an {\it HL-algebra} is an \ah-module $A$ satisfying 
Axiom~P1.
\item an {\it HP-algebra} satisfies Axioms A, P1 and P2.
\end{itemize}}
\label{hpldef}
\end{dfn}

Here is a little motivation for the definitions above. 
If $M$ is a symplectic manifold, then the algebra of
smooth functions on $M$ acquires a Poisson bracket. Since
a hyperk\"ahler manifold $M$ has 3 symplectic structures, the algebra
of smooth functions on $M$ has 3 Poisson brackets, and these
interact with the H-algebra $A$ of q-holomorphic functions on $M$,
that will be defined in Chapter~3. 

Our definition of HP-algebra is an attempt to capture the essential 
algebraic properties of this interaction. We may regard $A\oth Y$ as 
a subspace of $A\ot(Y^\d)^*$, and $(Y^\d)^*\cong\I$ by Definition 
\ref{ymoddef}. Thus $\xi_A$ is an antisymmetric map from $A\oth A$ 
to $A\ot\I$, i.e.~a triple of antisymmetric maps from $A\oth A$ to $A$. 
These 3 antisymmetric maps should be interpreted as the 3 Poisson 
brackets on the hyperk\"ahler manifold.

\subsection{Filtered and graded H-algebras}
\label{h23}

We begin by defining filtered and graded \ah-modules.

\begin{dfn} Let $U$ be an \ah-module. A {\it filtration
of\/ $U$} is a sequence $U_0,U_1,\dots$ of \ah-submodules
of $U$, such that $U_j\subset U_k$ whenever $j\le k$, and 
$U=\bigcup_{k=0}^\infty U_k$. We call $U$ a {\it filtered\/
\ah-module} if it has a filtration~$U_0,U_1,\dots$. 

Let $U$ be an \ah-module. A {\it grading of\/ $U$} is a 
sequence $U^0,U^1,\dots$ of \ah-submodules of $U$, such that 
$U=\bigoplus_{k=0}^\infty U^k$. We call $U$ a {\it graded\/
\ah-module} if it has a grading $U^0,U^1,\dots$. If $U$
is a graded \ah-module, define $U_k=\bigoplus_{j=0}^kU^j$.
Then $U_0,U_1,\dots$ is a filtration of $U$, so every
graded \ah-module is also a filtered \ah-module.

Let $U,V$ be filtered \ah-modules, and $\phi:U\ra V$ be an
\ah-morphism. We say that {\it $\phi$ is a filtered\/
\ah-morphism} if $\phi(U_k)\subset V_k$ for each $k\ge 0$.
Graded \ah-morphisms are also defined in the obvious way.
\end{dfn}

Here are axioms for filtered and graded H- and HP-algebras.

\begin{ax}{AF} $A$ is a filtered \ah-module.
\item $\H\cdot 1\subset A_0$.
\item For each $j,k$, $\mu_A(A_j\oth A_k)\subset A_{j+k}$.
\end{ax}

\begin{ax}{AG} $A$ is a graded \ah-module.
\item $\H\cdot 1\subset A^0$.
\item For each $j,k$, $\mu_A(A^j\oth A^k)\subset A^{j+k}$.
\end{ax}

\begin{axw}{PF} For each $j,k$, 
$\xi_A(A_j\oth A_k)\subset A_{j+k-1}\oth Y$.
\end{axw}

\begin{axw}{PG} For each $j,k$, 
$\xi_A(A^j\oth A^k)\subset A^{j+k-1}\oth Y$.
\end{axw}

\begin{dfn} 
{\setlength{\parskip}{-8pt}
\begin{itemize}
\setlength{\itemsep}{-4pt}
\setlength{\parsep}{0pt}
\item A {\it filtered H-algebra} satisfies Axioms A and AF.
\item A {\it graded H-algebra} satisfies Axioms A and AG.
\item a {\it filtered HP-algebra} satisfies Axioms A, AF, 
P1, P2, and~PF.
\item a {\it graded HP-algebra} satisfies Axioms A, AG, 
P1, P2, and~PG.
\end{itemize}}
\label{filterdef}
\end{dfn}

Morphisms of filtered and graded H-algebras are defined in the
obvious way, following Definition \ref{halgmordef}. The choice of 
the grading $j\!+\!k\!-\!1$ in Axioms PF and PG is not always 
appropriate, but depends on the situation. For some of our 
applications, the grading $j\!+\!k\!-\!2$ is better. Now let
$U$ be an \ah-module, and $V$ an \ah-submodule of $V$.
Then $U/V$ is naturally an $\H$-module. As $V'=U'\cap V$,
we may interpret $U'/V'$ as a real vector subspace of
$U/V$. Put $(U/V)'=U'/V'$. Then $U/V$ is an $\H$-module
with a real vector subspace $(U/V)'$. Note that $U/V$ may or 
may not be an \ah-module, because it may not satisfy the
condition of Definition~\ref{ahmdef}.

\begin{dfn} Let $A$ be a filtered H-algebra. Define $A_k=\{0\}$
for $k<0$ in $\Z$. We say that $A$ is a {\it stable filtered 
H-algebra}, or {\it SFH-algebra}, if for each $k\ge 0$, $A_k/A_{k-1}$ 
is a stable \ah-module. Let $B$ be a graded H-algebra. We say $B$ is 
a {\it stable graded H-algebra}, or {\it SGH-algebra}, if $B^k$ is 
stable for each~$k\ge 0$.
\end{dfn}

\begin{lem} Let $A$ be an SFH-algebra. Then for each\/ $j,k\ge 0$, 
$A_k$ and $A_j/A_{j-k}$ are stable \ah-modules. Let $j,k\ge 0$ and
$l>0$ be integers. Then the multiplication map $\mu_A:A_j\oth A_k\ra 
A_{j+k}$ induces a natural \ah-morphism 
$\mu^A_{jkl}:(A_j/A_{j-l})\oth (A_k/A_{k-l})\ra A_{j+k}/A_{j+k-l}$.
\label{gradelem}
\end{lem}

\begin{proof} We shall prove that $A_k$ is stable, by induction on 
$k$. Firstly, $A_0=A_0/A_{-1}$ is stable, by definition. Suppose by 
induction that $A_{k-1}$ is stable. The sequence
$0\ra A_{k-1}\ra A_k{\buildrel\pi_k\over\longra} A_k/A_{k-1}\ra 0$ is 
\ah-exact, and $A_{k-1}$ and $A_k/A_{k-1}$ are stable. Therefore, 
Proposition \ref{midminprop} shows that $A_k$ is stable, so all $A_k$ 
are stable, by induction. By a similar argument involving induction on 
$k$, $A_j/A_{j-k}$ is stable.

Now let $j,k,l$ be as given, and let $\pi_m:A_m\ra A_m/A_{m-l}$
be the natural projection, for $m\ge 0$. Because
$0\ra A_{j-l}\ra A_j{\buildrel\pi_j\over\longra} A_j/A_{j-l}\ra 0$ and 
$0\ra A_{k-l}\ra A_k{\buildrel\pi_k\over\longra} A_k/A_{k-l}\ra 0$ are 
\ah-exact sequences of stable \ah-modules, two applications of 
Proposition \ref{minexactprop} show that the sequence 
\begin{equation}
A_j\oth A_k{\buildrel\pi_j\smalloth\pi_k\over\longra}
(A_j/A_{j-l})\oth(A_k/A_{k-l})\ra 0
\label{jklseqeq}
\end{equation}
is \ah-exact at the middle term. Now $\mu_A$ maps $A_j\oth A_k$ to 
$A_{j+k}$. By identifying the kernels of $\pi_j\oth\pi_k$ and 
$\pi_{j+k}\circ\mu_A$, it can be seen that there exists a linear map 
$\mu^A_{jkl}$ as in the lemma, such that $\mu^A_{jkl}\circ(\pi_j\oth\pi_k)
=\pi_{j+k}\circ\mu_A$ as maps $A_j\oth A_k\ra A_{j+k}/A_{j+k-l}$. Using 
the \ah-exactness of \eq{jklseqeq} and that $\pi_{j+k}\circ\mu_A$ 
is an \ah-morphism, we deduce that $\mu^A_{jkl}$ is unique and is an 
\ah-morphism.
\end{proof}

Using this lemma, we make a definition.

\begin{dfn} 
\label{asympalgdef}
Let $A,B$ be SFH-algebras and let $l>0$ be an integer. By Lemma 
\ref{gradelem}, $A_j/A_{j-l}$ and\/ $B_j/B_{j-l}$ are stable 
\ah-modules for each\/ $j\ge 0$. We say that {\it $A$ and\/ $B$ are 
isomorphic to order $l$} if the following holds. For $j\ge 0$ there 
are \ah-isomorphisms $\phi_j:A_j/A_{j-l}\ra B_j/B_{j-l}$. These satisfy 
$\phi_j(A_k/A_{j-l})\subset B_k/B_{j-l}$ when $j-l\le k\le j$. Therefore 
$\phi_j$ projects to a map $A_j/A_{j-l+1}\ra B_j/B_{j-l+1}$, and 
$\phi_{j+1}$ restricts to a map $A_j/A_{j-l+1}\ra B_j/B_{j-l+1}$. Then 
$\phi_j=\phi_{j+1}$ on $A_j/A_{j-l+1}$. Also, for all\/ $j,k\ge 0$, 
$\mu^B_{jkl}\circ(\phi_j\oth\phi_k)=\phi_{j+k}\circ\mu^A_{jkl}$ as 
\ah-morphisms from $(A_j/A_{j-l})\oth(A_k/A_{k-l})$ to 
$B_{j+k}/B_{j+k-l}$. Here $\mu^A_{jkl}$ and $\mu^B_{jkl}$ are defined 
in Lemma~\ref{gradelem}.
\end{dfn}

There is a well-known way to construct a graded algebra from a 
filtered algebra (for instance \cite[p.~35-37]{BV} gives the 
associated graded Poisson algebra of a filtered Poisson algebra). 
Here is the analogue of this for H-algebras, which will be
applied in~\S\ref{h35}.

\begin{prop} Let\/ $A$ be an SFH-algebra. Define $B^k=A_k/A_{k-1}$ 
for $k\ge 0$. Define $B=\bigoplus_{k=0}^\infty B^k$. Then $B$ has 
the structure of an SGH-algebra, in a natural way.
\label{asgrhaprop}
\end{prop}

\begin{proof} As $A$ is an SFH-algebra, $B^k$ is a stable 
\ah-module, by definition. Let $j,k\ge 0$ be integers. Putting $l=1$, 
Lemma \ref{gradelem} defines an \ah-morphism $\mu^A_{jk1}:B^j\oth B^k
\ra B^{j+k}$. Let $\mu_B:B\oth B\ra B$ be the unique \ah-morphism, 
such that the restriction of $\mu_B$ to $B^j\oth B^k$ is $\mu^A_{jk1}$.
It is elementary to show that because $A$ is an H-algebra, $\mu_B$
makes $B$ into an H-algebra, and we leave this to the reader.
As $B$ satisfies Axiom AG, $B$ is graded, so $B$ is an SGH-algebra,
and the proposition is complete.
\end{proof}

We call the SGH-algebra $B$ defined in Proposition \ref{asgrhaprop} 
the {\it associated graded H-algebra of\/ $A$}. Note that $A$
is isomorphic to $B$ to order 1, in the sense of Lemma 
\ref{asympalgdef}. One might ask if the construction would work even 
if $A$ were only a filtered H-algebra. There are two problems here: 
firstly, $B^k$ might not be an \ah-module, and secondly, even if $B^j$ 
and $B^k$ were \ah-modules, the map $\pi_j\oth\pi_k:A_j\oth A_k\ra 
B^j\oth B^k$ might not be surjective. If it were not, we could only 
define $\mu^A_{jk1}$ uniquely on part of $B^j\oth B^k$. For these reasons 
we prefer SFH-algebras.

\subsection{Free and finitely-generated H-algebras}
\label{h24}

First we define ideals in SFH-algebras.

\begin{dfn} Let $A$ be an SFH-algebra, and $I$ an \ah-submodule 
of $A$. Set $I_k=I\cap A_k$ for $k\ge 0$. We say that {\it $I$ 
is a stable filtered ideal in $A$} if $1\notin I$, $\mu_A(I\oth A)
\subset I$, and for each $k\ge 0$, $I_k$ and $A_k/I_k$ are stable 
\ah-modules. Suppose that $J$ is an \ah-submodule of $A$, and
$I=\mu_A(J\oth A)$. Then we say that {\it $I$ is generated by~$J$}.
\label{idealdef}
\end{dfn}

The proof of the next lemma is similar to that of Proposition
\ref{asgrhaprop}, so we omit it.

\begin{lem}
Let\/ $A$ be an SFH-algebra, and\/ $I$ a stable filtered ideal
in $A$. Then there exists a unique SFH-algebra $B$, with a
filtered H-algebra morphism $\pi:A\ra B$, such that\/
$0\ra I_k{\buildrel\iota\over\ra}A_k{\buildrel\pi\over\ra}B_k\ra 0$
is an \ah-exact sequence for each $k\ge 0$. Here 
$\iota:I\ra A$ is the inclusion map.
\label{ideallem}
\end{lem}

We shall call this new H-algebra $B$ the {\it quotient of\/ $A$ by 
$I$}. Next, here is the definition of a free H-algebra.

\begin{dfn} Let $Q$ be an \ah-module. Define the {\it free
H-algebra $F^Q$ generated by $Q$} as follows. Put
$F^Q=\bigoplus_{k=0}^\infty\Sh^kQ$. Then $F^Q$ is an 
\ah-module. Now $(\Sh^kQ)\oth(\Sh^lQ)\subset\bigoth^{k+l}Q$,
and from Definition \ref{symmantisymm} there is an
\ah-module projection $\sigma_\H:\bigoth^{k+l}Q\ra
\Sh^{k+l}Q$. Define $\mu_{k,l}:(\Sh^kQ)\oth(\Sh^lQ)
\ra\Sh^{k+l}Q$ to be the restriction of 
$\sigma_\H$ to~$(\Sh^kQ)\oth(\Sh^lQ)$. 

Define $\mu_{F^Q}:F^Q\oth F^Q\ra F^Q$ to be the unique 
linear map such that the restriction of $\mu_{F^Q}$ to
$(\Sh^kQ)\oth(\Sh^lQ)$ is $\mu_{k,l}$. Recall that $\Sh^0Q=\H$, 
and define $1\in F^Q$ to be $1\in\H=\Sh^0Q$. It is easy to show 
that with these definitions, $F^Q$ is an H-algebra. The
{\it natural grading on $F^Q$} is $(F^Q)^k=\Sh^kQ$ for $k\ge 0$.
The {\it natural filtration on $F^Q$} is $F^Q_k=\bigoplus_{j=0}^k
(F^Q)^j$ for~$k\ge 0$.
\label{freealgdef}
\end{dfn}

If $Q$ is stable, then $\Sh^kQ$ is stable by Proposition 
\ref{symmantiprop}. Thus, if $Q$ is stable then $F^Q$ is an 
SFH-algebra (SGH-algebra) with the natural filtration (grading). 
The proof of the next lemma is trivial, and we omit it.

\begin{lem} Let\/ $A$ be an H-algebra, and let\/ $Q\subset A$ be
an \ah-submodule. Let\/ $\iota_Q:Q\ra F^Q$ be the inclusion map.
Then there is a unique H-algebra morphism $\phi_Q:F^Q\ra A$ such 
that $\phi_Q\circ\iota_Q:Q\ra A$ is the identity on~$Q$.
\label{genlem}
\end{lem}

\begin{dfn}
A {\it free SFH-algebra} is an SFH-algebra $A$ that is isomorphic,
as an H-algebra, to some $F^Q$, for finite-dimensional $Q$. Note
that the filtration on $A$ need not be the natural filtration on
$F^Q$. A {\it finitely-generated SFH-algebra}, or {\it FGH-algebra} 
is the quotient $B$ of a free SFH-algebra $A$ by a stable filtered 
ideal $I$ in $A$. Suppose that $Q$ is an \ah-submodule of $B$, and
that there exists an H-algebra isomorphism $A\cong F^Q$ identifying 
the maps $\pi:A\ra B$ and $\phi_Q:F^Q\ra B$. Then we say that 
{\it $Q$ generates the FGH-algebra~$B$}.
\label{fghdef}
\end{dfn}

The purpose of this definition is as follows. The polynomials
on an affine algebraic variety form a finitely-generated, filtered
algebra, and in algebraic geometry one studies this algebra to
learn about the variety. In the opinion of the author, FGH-algebras
are the best quaternionic analogue of finitely-generated, filtered algebras.
Clearly, they are finitely-generated, filtered H-algebras, and the
extra stability conditions we impose enable us to exploit the
`right-exactness' results of~\S\ref{h13}.

Moreover, the author believes that there is a wide class 
of noncompact hypercomplex manifolds, to which one can naturally 
associate an FGH-algebra. The study of FGH-algebras should be 
interpreted as the `quaternionic algebraic geometry' of these 
hypercomplex manifolds. Therefore, the author proposes that the 
study of FGH-algebras, from the algebraic point of view, may be 
interesting and worthwhile. More will be said on these ideas in 
Chapter~4.

We leave the proof of this final result as an exercise. It
will be useful later.

\begin{prop}
Let\/ $A$ be an SFH-algebra, and\/ $B$ the associated graded 
algebra, as in \S\ref{h23}. Then $B$ is also an SFH-algebra. 
Suppose that\/ $B$ is an FGH-algebra, generated by $B_k$. Then 
$A$ is an FGH-algebra, generated by~$A_k$.
\label{asgrgenprop}
\end{prop}

\section{Hypercomplex geometry}

We begin in \S\ref{h31} by defining hypercomplex manifolds, hyperk\"ahler
manifolds and q-holomorphic functions on hypercomplex manifolds,
and some elementary properties of q-holomorphic functions
are given. Section \ref{h32} proves that the vector space
$A$ of q-holomorphic functions on a hypercomplex manifold $M$ forms an 
H-algebra, and \S\ref{h33} shows that if $M$ is hyperk\"ahler, then 
$A$ is an HP-algebra. Section \ref{h34} discusses the
possibility of reconstructing a hypercomplex manifold from an
H-algebra of q-holomorphic functions upon it. Finally, 
\S\ref{h35} discusses hyperk\"ahler manifolds that are asymptotic
to a conical metric.

\subsection{Q-holomorphic functions on hypercomplex manifolds} 
\label{h31}

We begin by defining hypercomplex manifolds (\cite[p.~137-139]{Sal}) 
and hyperk\"ahler manifolds (\cite[p.~114-123]{Sal}). Let $M$ be a 
manifold of dimension $4n$. A {\it hypercomplex structure} on $M$ 
is a triple $(I_1,I_2,I_3)$ on $M$, where $I_j$ is a complex structure 
on $M$, and $I_1I_2=I_3$. A {\it hyperk\"ahler structure} on $M$ is a 
quadruple $(g,I_1,I_2,I_3)$, where $g$ is a Riemannian metric 
on $M$, $(I_1,I_2,I_3)$ is a hypercomplex structure on $M$, and 
$g$ is K\"ahler w.r.t.~each $I_j$. If $M$ has a hypercomplex 
(hyperk\"ahler) structure, then $M$ is called a {\it hypercomplex 
(hyperk\"ahler) manifold}.

If $M$ is a hypercomplex manifold, then $I_1,I_2,I_3$ satisfy the 
quaternion relations, so that each tangent space $T_mM$ is an
$\H$-module isomorphic to $\H^n$. Also, if $r_1,r_2,r_3\in{\bb R}$ 
with $r_1^2+r_2^2+r_3^2=1$, then $r_1I_1+r_2I_2+r_3I_3$ is a complex 
structure. Thus a hypercomplex manifold possesses a 2-dimensional family of 
integrable complex structures, parametrized by ${\cal S}^2$. We will
often use ${\cal S}^2$ to denote this family of complex structures.

Let $M$ be a hypercomplex manifold. For $k\ge 0$, define 
$\Omega^k=C^\infty(\Lambda^kT^*M)$, and
$\Omega^k(\H)=C^\infty(\H\otimes\Lambda^kT^*M)$.
Then $\Omega^1$ is the vector space of smooth 1-forms on $M$,
and $\Omega^0(\H)$ is the vector space of smooth, 
quaternion-valued functions on $M$. Define an operator
$D:\Omega^0(\H)\rightarrow\Omega^1$ by
\begin{equation}
D(a_0+a_1i_1+a_2i_2+a_3i_3)
=da_0+I_1(da_1)+I_2(da_2)+I_3(da_3),
\label{qdirac}
\end{equation}
where $a_0,\dots,a_3$ are smooth real functions on $M$.

We define a {\it q-holomorphic function on $M$} to
be an element $a=a_0+a_1i_1+a_2i_2+a_3i_3$ of 
$\Omega^0(\H)$ for which $D(a)=0$. The term
q-holomorphic is short for {\it quaternion-holomorphic},
and it is intended to indicate that a q-holomorphic
function on a hypercomplex manifold is the appropriate
quaternionic analogue of a holomorphic function on a 
complex manifold. The operator $D$ of \eq{qdirac} should 
be thought of as the quaternionic analogue of the 
$\overline\partial$ operator on a complex manifold.
It can be seen that $D(a)=0$ is equivalent to the equation
\begin{equation}
da-I_1(da)i_1-I_2(da)i_2-I_3(da)i_3=0,
\label{altqdirac}
\end{equation}
where each term is an $\H$-valued 1-form, $I_j$ acts on 1-forms
and $i_j$ acts on $\H$ by multiplication.

Now in 1935, Fueter defined a class of `regular' $\H$-valued
functions on $\H$, using an analogue of the Cauchy-Riemann
equations, and Fueter and his co-workers went on to develop the 
theory of {\it quaternionic analysis}, by analogy with complex analysis. 
An account of this theory, with references, is given by Sudbery 
in \cite{Sud}. On the hypercomplex manifold $\H$, Fueter's definition
of regular function coincides with that of q-holomorphic 
function, given above. We shall make little reference to
Fueter's theory, because we are interested in rather different
questions. However, in Chapter 4 we will use our theory to give 
an elegant construction of the spaces of homogeneous q-holomorphic 
functions on $\H$, which are important in quaternionic analysis.

Suppose that $M$ is hyperk\"ahler. Then using the metric $g$ on $M$
we construct the operator $D^*:\Omega^1\rightarrow
\Omega^0(\H)$, which is given by
\begin{equation}
D^*(\alpha)=d^*\alpha-d^*(I_1\alpha)i_1
-d^*(I_2\alpha)i_2-d^*(I_3\alpha)i_3.
\label{dstar}
\end{equation}
Now for a smooth real function $f$ on a K\"ahler manifold,
$d^*(Idf)=0$. Using this we can show that 
$D^*D(a)=\Delta a$, where $\Delta$ is the usual Laplacian.
Thus q-holomorphic functions on a hyperk\"ahler manifold are harmonic.
When $n=1$, $D$ is elliptic, and is the Dirac operator $D_+$. 
When $n>1$, $D$ is overdetermined elliptic. 

Here are two basic properties of q-holomorphic functions.

\begin{lem} Suppose that $a$ is q-holomorphic on $M$,
and that $q\in\H$. Then $qa$ is q-holomorphic.

Suppose that $i=r_1i_1+r_2i_2+r_3i_3\in\I$ satisfies 
$i^2=-1$, and that $I=r_1I_1+r_2I_2+r_3I_3$ is the 
corresponding complex structure on $M$. Suppose that 
$y+zi$ is a complex function on $M$ that is holomorphic
w.r.t.~$I$. Then $y+zi$ is q-holomorphic on $M$, regarding 
$y+zi=y+zr_1i_1+zr_2i_2+zr_3i_3$ as an element of~$\Omega^0(\H)$.
\label{holqhollem}
\end{lem}

\begin{proof} Let $q=q_0+q_1i_1+q_2i_2+q_3i_3\in\H$,
and define $Q=q_0+q_1I_1+q_2I_2+q_3I_3$, regarding $Q$
as an endomorphism of $T^*M$. Let $a\in\Omega^0(\H)$.
It is easy to verify that $D(qa)=Q\cdot D(a)$. Thus 
$D(qa)=0$ if $D(a)=0$, and $qa$ is q-holomorphic
whenever $a$ is q-holomorphic. This proves the first part.

If $y+zi$ is holomorphic w.r.t.~$I$, then $dy+I(dz)=0$
by the Cauchy-Riemann equations. But 
\begin{align*}
0&=dy+I(dz)=dy+(r_1I_1+r_2I_2+r_3I_3)(dz)\\
&=D(y+zr_1i_1+zr_2i_2+zr_3i_3),
\end{align*}
so that $y+zr_1i_1+zr_2i_2+zr_3i_3$ is q-holomorphic. This
completes the lemma.
\end{proof}

The operator $D$ of \eq{qdirac} was also studied by Baston
\cite{Bas}, upon quaternionic manifolds rather than hypercomplex 
manifolds. He calls $D$ the {\it Dirac-Fueter operator}, and uses 
the Penrose transform to interpret $D$ as a holomorphic object on 
the twistor space $Z$ of $M$. Baston \cite[p.~44-45]{Bas} shows that
$\Ker D$ on a quaternionic manifold $M$ can be identified with the
sheaf cohomology group $H^1(Z,{\cal O}_Z(-3))$, giving a twistor
interpretation of q-holomorphic functions. He also constructs an
exact complex of operators resolving~$D$,~\cite[p.~43-44]{Bas}. 

\subsection{Hypercomplex manifolds and H-algebras}
\label{h32}

Let $M$ be a hypercomplex manifold, and let $A$ be the vector space 
of q-holomorphic functions on $M$. In this section we will 
prove that $A$ is an H-algebra.

\begin{dfn} Let $M$ be a hypercomplex manifold. Define 
$A\subset\Omega^0(\H)$ to be the vector space of 
q-holomorphic functions on $M$. Let $a\in A$, and define 
$(q\cdot a)(m)=q(a(m))$ for $m\in M$. Then $q\cdot a\in A$ by 
Lemma \ref{holqhollem}, and this gives an $\H$-action on $A$, 
so $A$ is an $\H$-module. Define a subspace $A'$ in $A$ by
$A'=\bigl\{a\in A:a(m)\in\I$ for all $m\in M\bigr\}$. 
For each $m\in M$, define $\theta_m:A\ra\H$ by $\theta_m(a)=a(m)$. 
Then $\theta_m\in A^\t$, and if $a\in A'$ then $\theta_m(a)\in\I$,
so that~$\theta_m\in A^\d$.

Suppose $a\in A$, and $\alpha(a)=0$ for all $\alpha\in A^\d$. 
Since $\theta_m\in A^\d$, $a(m)=0$ for each $m\in M$, and so $a=0$. 
Thus $A$ is an \ah-module, by Definition \ref{ahmdef}.
Define the element $1\in A$ to be the constant function on $M$ 
with value 1. Then $1\notin A'$, but $\I\cdot 1\subset A'$.
We will also write $A_M$ for $A$, when we wish to specify the 
manifold~$M$.
\label{adef}
\end{dfn}

The following proposition gives us a greater understanding 
of the quaternionic tensor product.

\begin{prop} Let\/ $M$ and $N$ be hypercomplex manifolds, and let\/
$U,V$ be \ah-submodules of the \ah-modules $A_M,A_N$ of 
q-holomorphic functions on $M,N$ respectively. Define $W$ 
to be the vector space of smooth, $\H$-valued functions $w$ 
on $M\times N$, such that for each $m\in M$, the function 
$n\mapsto w(m,n)$ lies in $V$, and for each $n\in N$,
the function $m\mapsto w(m,n)$ lies in $U$. Then each such\/
$w$ is a q-holomorphic function on $M\t N$, and\/
$W$ is an \ah-submodule of\/ $A_{M\t N}$. Also, there
is a canonical injective \ah-morphism $\phi:U\oth V\ra W$.
If\/ $U,V$ are finite-dimensional, $\phi$ is an \ah-isomorphism.
\label{qholtensprop}
\end{prop}

\begin{proof} Suppose $w:M\t N\ra\H$ is a smooth
function, such that for each $m\in M$, the function 
$n\mapsto w(m,n)$ lies in $V$, and for each $n\in N$,
the function $m\mapsto w(m,n)$ lies in $U$. The
condition for $w$ to be q-holomorphic is $D(w)=0$.
But $D(w)=D_M(w)+D_N(w)$, where $D_M$ involves only 
derivatives in the $M$ directions, and $D_N$ only derivatives 
in the $N$ directions. 

Let $n\in N$. Then the function $m\mapsto w(m,n)$ is equal 
to some $u\in U$. Thus $D_M(w)(m,n)=D(u)(m)$. But the functions 
in $U$ are q-holomorphic, so $D(u)=0$. Therefore $D_M(w)=0$, and 
similarly $D_N(w)=0$. So $D(w)=0$, and $w$ is q-holomorphic, as 
we have to prove. It is clear that the space $W$ of such 
functions $w$ is closed under addition and multiplication 
by $\H$. Thus $W$ is an $\H$-submodule of $A_{M\t N}$, so 
$W$ is an \ah-submodule of~$A_{M\t N}$.

Now let $\epsilon\in U\oth V$. Then $\epsilon\in\H\ot
(U^\d)^*\ot(V^\d)^*$, so $\epsilon$ defines a linear
map $\epsilon:U^\d\ot V^\d\rightarrow\H$. Define a map
$w:M\t N\ra\H$ by $w(m,n)=\epsilon(\theta_m\ot\theta_n)$.
For $m\in M$, define $w_m:N\ra\H$ by $w_m(n)=w(m,n)$.
Since $\epsilon\in U\oth V$, $\epsilon\in(U^\d)^*
\ot\iota_V(V)$, so $w_m\in\iota(V)$, regarding $w_m$
as an element of $\H\ot(V^\d)^*$. Thus $w_m\in V$.
Similarly, defining $w_n(m)=w(m,n)$ for $n\in N$,
we find $w_n\in U$ for each~$n\in N$.

To show that $w\in W$, we only need to show that $w$
is smooth. In general, if $f$ is a function on $M\t N$,
such that for each $m\in M$, the function $n\mapsto f(m,n)$ 
is smooth, and for each $n\in N$, the function $m\mapsto f(m,n)$ 
is smooth, it does not follow that $f$ is smooth. However, 
because real tensor products involve only finite sums as in 
\S\ref{h11}, $\epsilon\in(U^\d)^*\ot\iota_V(V)$ implies that $w_m$ 
lies in some finite-dimensional subspace of $V$ for all $m\in M$,
and similarly $w_n$ lies in a finite-dimensional subspace
of $U$ for all $n\in N$. These imply that $w$ is smooth.
Thus~$w\in W$.

Define $\phi(\epsilon)=w$. In this way we define a map
$\phi:U\oth V\ra W$. It is easy to show that $\phi$
is an \ah-morphism. Also, if $w=0$ it is easily
seen that $\epsilon=0$, so $\phi$ is injective. Thus
$\phi$ is an injective \ah-morphism, as we have to prove.
Suppose $U,V$ are finite-dimensional, and let $w\in W$.
We must find $\epsilon$ in $U\oth V$ such 
that~$\phi(\epsilon)=w$. 

Choose bases of the form $\{\theta_{m_b}:b=1,\dots,k\}$ for $U^\d$ 
and $\{\theta_{n_c}:c=1,\dots,l\}$ for $V^\d$. It can be shown that 
such bases exist. Let $\epsilon:U^\d\ot V^\d\ra\H$ be the unique 
linear map satisfying $\epsilon(\theta_{m_b}\ot\theta_{n_c})=
w(m_b,n_c)$ for $b=1,\dots,k$, $c=1,\dots,l$. One may prove that 
$\epsilon\in U\oth V$, and $\phi(\epsilon)=w$. Thus $\phi$ is an 
injective and surjective \ah-morphism, and clearly is an 
\ah-isomorphism. This completes the proof.
\end{proof}

The following lemma is trivial, and the proof will be omitted.

\begin{lem} Suppose $M$ is a hypercomplex manifold, and $N$ is a
hypercomplex submanifold of $M$. If\/ $a$ is a q-holomorphic function
on $M$, then $a\vert_N$ is q-holomorphic on $N$. Let\/
$\rho:A_M\ra A_N$ be the restriction map. Then $\rho$ is
an \ah-morphism.
\label{restrictlem}
\end{lem}

Now we can define the multiplication map $\mu_A$ on $A$. 

\begin{dfn} Let $M$ be a hypercomplex manifold, and $A$ the \ah-module
of q-holomorphic functions on $M$. By Proposition 
\ref{qholtensprop} there is a canonical \ah-morphism 
$\phi:A\oth A\ra A_{M\t M}$. Now $M$ is embedded in $M\t M$ 
as the diagonal submanifold $\bigl\{(m,m):m\in M\bigr\}$,
and this is a hypercomplex submanifold of $M\t M$, isomorphic to $M$
as a hypercomplex manifold. Therefore Lemma \ref{restrictlem} gives
an \ah-morphism $\rho:A_{M\t M}\ra A$. Define an \ah-morphism
$\mu_A:A\oth A\ra A$ by~$\mu_A=\rho\circ\phi$.
\label{amudef}
\end{dfn}

Here is the main result of this section.

\begin{thm} Let\/ $M$ be a hypercomplex manifold. Then Definition
\ref{adef} defines an \ah-module $A$ and an element\/ 
$1\in A$, and Definition \ref{amudef} defines an
\ah-morphism $\mu_A:A\oth A\ra A$. With these definitions,
$A$ is an H-algebra in the sense of\/~\S\ref{h21}.
\label{qholhalgthm}
\end{thm}

\begin{proof} We must show that Axiom A is satisfied.
Parts $(i)$ and $(ii)$ are trivial. For part $(iii)$,
observe that the permutation map $A\oth A\ra A\oth A$
that swaps round the factors, is induced by the map
$M\t M\ra M\t M$ given by $(m_1,m_2)\mapsto(m_2,m_1)$.
Since the diagonal submanifold is invariant under this,
it follows that $\mu_A$ is invariant under permutation, 
and so~$\Lambh^2A\subset\Ker\mu_A$.

Let $\Delta^2_M$ be the `diagonal' submanifold in $M\t M$,
and let $\Delta^3_M$ be the `diagonal' submanifold in
$M\t M\t M$. We interpret part $(iv)$ as follows.
$A\oth A\oth A$ is a space of q-holomorphic functions
on $M\t M\t M$. The maps $\mu_A\oth\id$ and $\id\oth\mu_A$
are the maps restricting to $\Delta^2_M\t M$ and 
$M\t\Delta^2_M$ respectively. Thus $\mu_A\circ(\mu_A\oth\id)$ 
is the result of first restricting to $\Delta^2_M\t M$ and 
then to $\Delta^3_M$, and $\mu_A\circ(\id\oth\mu_A)$ is the 
result of first restricting to $M\t\Delta^2_M$ and then to 
$\Delta^3_M$. Clearly $\mu_A\circ(\mu_A\oth\id)=
\mu_A\circ(\id\oth\mu_A)$, proving part~$(iv)$.

Interestingly, the proof of part $(iv)$ does not use the 
associativity of quaternion multiplication. This raises
the possibility of generalizing the definition to give
`associative algebras over a nonassociative field'.
Part $(v)$ is given in Definition \ref{adef}. Finally,
part $(vi)$ follows easily from the fact that $1$ is
the identity in $\H$. Thus all of Axiom A of \S\ref{h21}
applies, and $A$ is an H-algebra.
\end{proof}

One problem with the H-algebra of all q-holomorphic 
functions on a hypercomplex manifold is that it is too large to
work with -- it is not in general finitely-generated,
for instance. Therefore, it is convenient to restrict 
to H-subalgebras of functions satisfying some condition.
The condition we shall use is that of {\it polynomial growth}.

\begin{dfn} Let $M$ be a hypercomplex manifold, and let $r:M\ra[0,\infty)$
be a given continuous function. Suppose $f\in\Omega^0(\H)$ on $M$,
and let $k\ge 0$ be an integer. We say that {\it $f$ has 
polynomial growth of order $k$}, written $f=O(r^k)$, if
there exist positive constants $C_1,C_2$ such that
$\md{f}\le C_1+C_2 r^k$ on~$M$.

The H-algebra of q-holomorphic functions on $M$ is $A$. For 
integers $k\ge 0$, define $P_k=\bigl\{a\in A:a=O(r^k)\bigr\}$,
and define $P=\bigcup_{k=0}^\infty P_k$. Then $P$ is a filtered 
\ah-module. It is easy to see that $P$ is an H-subalgebra of $A$, 
and satisfies Axiom AF of \S\ref{h23}. Thus, $P$ is a filtered 
H-algebra. We call $P$ the {\it filtered H-algebra of 
q-holomorphic functions of polynomial growth on $M$}, and we 
write $P_M$ for $P$ when we wish to specify the manifold~$M$.
\label{pgdef}
\end{dfn}

The main example we have in mind in making this definition, is
the case that $M$ is a complete, noncompact hyperk\"ahler manifold, and
$r:M\ra [0,\infty)$ is the distance function from some point
$m_0\in M$. Then $P_k$ is independent of the choice of base
point $m_0$. In good cases, such as those discussed in Chapter 4,
$P$ is an FGH-algebra. Again, in a good case, $P$ determines the 
hypercomplex structure of $M$ explicitly and uniquely. Thus we can define 
the hypercomplex structure of $M$ completely using only a finite-dimensional 
amount of algebraic data.

\subsection{Hyperk\"ahler manifolds and HP-algebras}
\label{h33}

Let $M$ be a hyperk\"ahler manifold, and $A$ the vector space of 
q-holomorphic functions on $M$. Since $M$ is hypercomplex, $A$ is an 
H-algebra by Theorem \ref{qholhalgthm}. In this section we will 
see that $A$ is also an HP-algebra in the sense of \S\ref{h22}. 
To save space, and because we have wandered from the main
subject of the paper, we shall omit the proofs of Proposition 
\ref{thetprop} and Theorem \ref{hkhpalg}. The proofs are 
elementary calculations, though not especially easy, and the 
author can supply them to the interested reader on request.

\begin{dfn} Let $M$ be a hyperk\"ahler manifold. Then $M\t M$ is also 
a hyperk\"ahler manifold. Let $\Delta^2_M=\bigl\{(m,m):m\in M\bigr\}$.
Then $\Delta^2_M$ is a hyperk\"ahler submanifold of $M\t M$. We shall write 
$M\t M=M^1\t M^2$, using the superscripts ${}^1$ and ${}^2$
to distinguish the two factors. Let $\nabla$ be the Levi-Civita 
connection on $M$. Define $\nabla^1,\nabla^2$ to be the lift of 
$\nabla$ to the first and second factors of $M$ in $M\t M$ 
respectively. Then $\nabla^1$ and $\nabla^2$ commute. Let 
$\nabla^{12}$ be the Levi-Civita connection on $M\t M$. 
Then~$\nabla^{12}=\nabla^1+\nabla^2$.

Let $x\in A_{M\t M}$. Then $\nabla^1\nabla^2x\in C^\infty
(\H\ot T^*M^1\ot T^*M^2)$ over $M\t M$. Restrict 
$\nabla^1\nabla^2x$ to $\Delta^2_M$. Then $\Delta^2_M\cong M$ 
and $T^*M^1\vert_{\Delta^2_M}\cong T^*M^2\vert_{\Delta^2_M}
\cong T^*M$. Thus $\nabla^1\nabla^2x\vert_{\Delta^2_M}\in 
C^\infty(\H\ot T^*M\ot T^*M)$ over $M$. Define a linear map 
$\Theta:A_{M\t M}\ra\Omega^0(\H)\ot\I$ by
\begin{equation}
\begin{split}
\Theta(x)=\bigl\{g^{ab}(I_1)_a^c\nabla_b^1\nabla_c^2x
\vert_{\Delta^2_M}\bigr\}\ot i_1
&+\bigl\{g^{ab}(I_2)_a^c\nabla_b^1\nabla_c^2x
\vert_{\Delta^2_M}\bigr\}\ot i_2\\
&+\bigl\{g^{ab}(I_3)_a^c\nabla_b^1\nabla_c^2x
\vert_{\Delta^2_M}\bigr\}\ot i_3,
\end{split}
\label{thetdef}
\end{equation}
using index notation for tensors on $M$ in the obvious way.
Here $g$ is the hyperk\"ahler metric on $M$, and $I_1,I_2,I_3$ the
complex structures.
\label{mthetdef}
\end{dfn}

Here are some properties of $\Theta$.

\begin{prop} This map satisfies $\Theta(x)\in A\ot\I$.
Also, $\Theta:A_{M\t M}\ra A\ot\I$ is an \ah-morphism,
and if\/ $\Theta(x)=x_1\ot i_1+x_2\ot i_2+x_3\ot i_3$ and\/
$m\in M$, then~$x_1(m)i_1+x_2(m)i_2+x_3(m)i_3=0\in\H$.
\label{thetprop}
\end{prop}

Now we can define the map $\xi_A$ of \S\ref{h22}.

\begin{dfn} Proposition \ref{qholtensprop} defines an
\ah-morphism $\phi:A\oth A\ra A_{M\t M}$. Definition
\ref{mthetdef} and Proposition \ref{thetprop} define an
\ah-morphism $\Theta:A_{M\t M}\ra A\ot\I$. Let $Y$ be the 
\ah-module of Definition \ref{ymoddef}. Then $(Y^\d)^*\cong\I$. 
Recall that $A\cong\iota_A(A)$, so we may identify
$A\ot\I\cong\iota_A(A)\ot(Y^\d)^*\subset\H\ot(A^\d)^*\ot(Y^\d)^*$.
Define $\xi_A:A\oth A\ra\iota_A(A)\ot(Y^\d)^*$ to be the
composition~$\xi_A=\Theta\circ\phi$.
\label{xidef}
\end{dfn}

Here is the main result of this section.

\begin{thm} This $\xi_A$ maps $A\oth A$ to $A\oth Y$. It is 
an \ah-morphism, and satisfies Axioms P1 and P2 of\/ \S\ref{h22}. 
Thus, by Theorem \ref{qholhalgthm} and Definition \ref{hpldef},
if\/ $M$ is a hyperk\"ahler manifold, then the vector space $A$ of 
q-holomorphic functions on $M$ is an HP-algebra.
\label{hkhpalg}
\end{thm}

Given a continuous function $r:M\ra[0,\infty)$, Definition
\ref{pgdef} defines the filtered H-algebra $P$ of q-holomorphic
functions on $M$ of polynomial growth. It is natural to ask
whether $P$ is a filtered HP-algebra. One must show that Axiom PF 
of \S\ref{h23} holds, which relates the Poisson bracket and the 
filtration. This is not automatic, but depends on the asymptotic
properties of the hyperk\"ahler structure and the function $r$ on $M$,
and must be verified for each case. These properties also 
determine whether $P$ is closed under $\xi_A$ at all.

\subsection{Reconstructing a hypercomplex manifold from its H-algebra}
\label{h34}

In \S\ref{h32} we saw that the vector space of q-holomorphic
functions on a hypercomplex manifold is an H-algebra, providing a
transform from geometric to algebraic objects. Now we shall
consider whether this transform can be reversed. It turns out
that under certain circumstances, an H-algebra does explicitly
determine a unique hypercomplex manifold. Therefore, it should be possible
to construct new hypercomplex manifolds by writing down their H-algebras.

Throughout this section, let $M$ be a hypercomplex manifold, let $A$
be the H-algebra of q-holomorphic functions on $M$, let $P$
be an H-subalgebra of $A$, and let $Q$ be an \ah-submodule of
$P$ that generates $P$ as an H-algebra, in the sense of \S\ref{h24}.
From \S\ref{h32}, if $m\in M$, then $m$ defines an $\H$-linear map 
$\theta_m:A\ra\H$ such that $\theta_m\in A^\d$. Now $\H$ is itself 
an H-algebra in the obvious way, and it is easy to see that in 
fact $\theta_m$ is an H-algebra morphism. Therefore 
$\theta_m\vert_P:P\ra\H$ is also an H-algebra morphism,
and~$\theta_m\vert_P\in P^\d$.

Conversely, suppose $\theta\in P^\d$, so that $\theta:P\ra\H$ is 
an $\H$-linear map. From Definition \ref{halgmordef} we calculate 
that $\theta$ is an H-algebra morphism if and only if $\theta$ 
satisfies the quadratic equation $\mu_P^\d(\theta)=\lambda_{P,P}
(\theta\ot\theta)$ in $(P\oth P)^\d$, where $\mu_P^\d:P^\d\ra
(P\oth P)^\d$ is the dual of the multiplication map $\mu_P$, and 
$\lambda_{P,P}:P^\d\ot P^\d\ra(P\oth P)^\d$ is defined in~\S\ref{h11}.

Suppose $\theta_1,\theta_2:P\ra\H$ are H-algebra morphisms, and 
that $\theta_1\vert_Q=\theta_2\vert_Q$. Because $Q$ generates
$P$, it is easy to see that $\theta_1=\theta_2$. Hence, H-algebra
morphisms from $P$ are determined by their restrictions to $Q$.
As $\theta_1\in P^\d$, we have $\theta_1\vert_Q\in Q^\d$.
The principal case we have in mind is that $P$ is a `polynomial
growth' H-subalgebra as in Definition \ref{pgdef}, and $Q$ is 
finite-dimensional. By restricting to $Q$ we can work in a
finite-dimensional situation. Motivated by the above, we make 
the following definition.

\begin{dfn} Let $P$ be an H-algebra, and $Q$ an \ah-submodule of 
$P$ that generates $P$. Define $M_{P,Q}$ by
\begin{equation*}
M_{P,Q}=\bigl\{\theta\vert_Q:\theta\in P^\d,\quad
\mu_P^\d(\theta)=\lambda_{P,P}(\theta\ot\theta)\bigr\},
\end{equation*}
so that $M_{P,Q}$ is a closed subset of $Q^\d$. Now suppose that
$M$ is a hypercomplex manifold, and $P$ an H-subalgebra of the H-algebra
$A$ of q-holomorphic functions on $M$. For each $m\in M$, 
$\theta_m:A\ra\H$ is an H-algebra morphism, and so 
$\theta_m\vert_Q$ lies in $M_{P,Q}$. Define a map 
$\pi_{P,Q}:M\ra M_{P,Q}$ by~$\pi_{P,Q}(m)=\theta_m\vert_Q$.
\end{dfn}

\begin{lem} Suppose $Q$ is finite-dimensional. Then 
$M_{P,Q}$ is an affine real algebraic subvariety of\/ $Q^\d$,
that is, it is the zeros of a finite collection of
polynomials on~$Q^\d$.
\label{algvarlem}
\end{lem}

\begin{proof} By Lemma \ref{genlem} there is an H-algebra
morphism $\phi_Q:F^Q\ra P$, and as $Q$ generates $P$, $\phi_Q$
is surjective. Let $I\subset F^Q$ be the kernel of $\phi_Q$.
As $F^Q$ is a free H-algebra, each element $x$ of $Q^\d$ 
defines a unique H-algebra morphism $\theta_x:F^Q\ra\H$, that
restricts to $x$ on $Q$. Clearly, if $x\in Q^\d$, then 
$x\in M_{P,Q}$ if and only if~$I\subset\Ker\theta_x$.

For each $y\in F^Q$, define a function $\psi_y:Q^\d\ra\H$
by $\psi_y(x)=\theta_x(y)$ for $x\in Q^\d$. It is easy to see
that if $y\in F^Q_k$, then $\psi_y$ is an $\H$-valued polynomial
on $Q^\d$ of degree at most $k$. But $x\in M_{P,Q}$ if and
only if $\psi_y(x)=0$ for each $y\in I$. Thus $M_{P,Q}$ is the 
zeros of a collection of polynomials on $Q^\d$. By Hilbert's 
Basis Theorem, we can choose a finite number of polynomials, 
which define~$M_{P,Q}$. 
\end{proof}

Our aim is to recover $M$ and its hypercomplex structure from $P$ and $Q$. 
If we are lucky, $\pi_{P,Q}$ will be (locally) a bijection, so 
that $M_{P,Q}$ gives us the manifold $M$, at least as a set. Then
we can try and use $P$ to define a hypercomplex structure on $M_{P,Q}$.
However, there are a number of ways in which this process
could fail.
\begin{itemize}
\item $M_{P,Q}$ might not be a submanifold of $Q^\d$.
\item $\pi_{P,Q}$ might not be (locally) injective.
\item $\pi_{P,Q}$ might not be (locally) surjective.
\item Even if $M_{P,Q}$ is a submanifold of $Q^\d$ and
$\pi_{P,Q}$ is a diffeomorphism, we may be unable to
define the hypercomplex structure on $M_{P,Q}$, because $P$ may
contain only partial information about the structure.
\end{itemize}

Unfortunately, all of these possibilities do occur, and examples 
will be given in \S\ref{h42}. Sometimes the hypercomplex manifold $M$ can 
be reconstructed, and sometimes not; for polynomial growth 
H-algebras $P$ and complete $M$, this seems to depend only on the 
asymptotic behaviour of $M$ at infinity. Next we shall explain how, 
in good cases, the hypercomplex structure of $M$ may be recovered from~$P$.

\begin{lem} Let\/ $M$ be a hypercomplex manifold of dimension $4k$, $A$ 
the H-algebra of q-holomorphic functions on $M$, $P$ an 
H-subalgebra of $A$, and\/ $Q$ an \ah-submodule of\/ $P$ 
generating $P$. Let\/ $m\in M$. Then the derivative of\/ 
$\pi_{P,Q}$ at\/ $m$ gives a linear map~$d_m\pi_{P,Q}:T_mM\ra Q^\d$.

Regard\/ $Q$ as a vector space of\/ $\H$-valued functions 
on $Q^\d$. Pulling these functions back to $T_mM$ using 
$d_m\pi_{P,Q}$ gives a linear map $(d_m\pi_{P,Q})^*:Q\ra 
T_m^*M\ot\H$. Define $V_m=\Im\bigl((d_m\pi_{P,Q})^*\bigr)$,
so that\/ $V_m$ is a linear subspace of\/ $T_m^*M\ot\H$.
Then $\dim V_m\le 12k$, and\/ $V_m$ determines the hypercomplex 
structure on $T_mM$ if and only if\/~$\dim V_m=12k$.
\label{recondlem}
\end{lem}

\begin{proof} Let $v\in V_m$, and write
$v=v_0\ot 1+v_1\ot i_1+v_2\ot i_2+v_3\ot i_3$. Now $v$ is the
first derivative at $m$ of some element of $Q$, which is a 
q-holomorphic function on $M$. By definition of q-holomorphic
function, it follows that $v_0+I_1v_1+I_2v_2+I_3v_3=0$, where
$I_1,I_2,I_3$ are the complex structures on $T_m^*M$. Let
$W_m\subset T_m^*M\ot\H$ be the subspace of elements 
$w_0\ot 1+w_1\ot i_1+w_2\ot i_2+w_3\ot i_3$ satisfying 
$w_0+I_1w_1+I_2w_2+I_3w_3=0$. Then $W_m$ has dimension $12k$,
and $V_m\subset W_m$, so $\dim V_m\le 12k$, as we have to prove.

If $\dim V_m=12k$, then $V_m=W_m$. But $W_m$ determines 
$I_1,I_2$ and $I_3$. (For instance, $w_1=I_1w_0$ if and only
if $w_0\ot 1+w_1\ot i_1\in W_m$, so $W_m$ determines $I_1$.)
Thus if $\dim V_m=12k$, then $V_m$ determines the hypercomplex 
structure on $T_mM$. Now $W_m$ is an $\H$-submodule of $T_m^*M\ot\H$, 
and it can be shown that if $W$ is an $\H$-submodule of $T^*M\ot\H$
such that $W\ne W_m$, $\dim W=12k$, and $W$ is close to $W_m$,
then $W$ determines a different hypercomplex structure. It is easy
to see that if $\dim V_m<12k$, then we may choose such a $W$ with 
$V_m\subset W$. Thus $V_m$ is consistent with two different 
hypercomplex structures. Therefore $V_m$ determines the hypercomplex 
structure on $T_mM$ if and only if~$\dim V_m=12k$.
\end{proof}

\begin{cor} In the situation of Lemma \ref{recondlem}, suppose 
that\/ $\dim V_m=12k$, and that\/ $\pi_{P,Q}$ is surjective near 
$m$. Then $P$ determines the hypercomplex structure of\/ $M$ near~$m$.
\label{recondcor}
\end{cor}

\begin{proof} As $\dim V_m=12k$, we deduce that $d_m\pi_{P,Q}$ 
is injective. Together with the surjectivity assumption, this 
implies that $\pi_{P,Q}$ is a diffeomorphism near $m$. Also,
$\dim V_n=12k$ for $n$ near $m$. The conclusion follows from the lemma.
\end{proof}

In the special case $\dim M=4$, we can give a more convenient 
condition for $P$ to determine the hypercomplex structure of~$M$.

\begin{prop} Let\/ $M$ be a hypercomplex manifold of dimension $4$, $A$ the 
H-algebra of q-holomorphic functions on $M$, $P$ an H-subalgebra of 
$A$, and\/ $Q$ an \ah-submodule of\/ $P$ generating $P$. Suppose 
that\/ $Q$ is stable, that\/ $M_{P,Q}$ is a submanifold of\/ $Q^\d$, 
and that\/ $\pi_{P,Q}$ is a (local) diffeomorphism. Then $P$ 
determines the (local) hypercomplex structure of\/~$M$.
\label{4dimprop}
\end{prop}

\begin{proof} Because of Corollary \ref{recondcor}, it is
sufficient to show that $\dim V_m=12$ for each $m\in M$. 
Suppose for a contradiction that $m\in M$ and $\dim V_m<12$. 
As $\pi_{P,Q}$ is a diffeomorphism, $d_m\pi_{P,Q}$ is injective. 
Let $R=\Ker\bigl((d_m\pi_{P,Q})^*\bigr)$, so that $R$ is a proper
\ah-submodule of $Q$, and $V_m\cong Q/R$. Since $d_m\pi_{P,Q}$ is 
injective and $\dim M=4$, we have~$\dim Q^\d=\dim R^\d+4$.

Now $V_m$ is an $\H$-module with $\dim V_m<12$, so $\dim V_m\le 8$.
But $\dim Q=\dim R+\dim V_m$, so $\dim Q\le\dim R+8$.
By Definition \ref{mindef}, using $\dim Q^\d=\dim R^\d+4$, we
see that the virtual dimension of $R$ is greater or equal to that 
of $Q$. Using the stability of $Q$, we then prove that for
each nonzero $q\in\I$, the image of $\id\oth\chi_q:Q\oth X_q\ra Q$
lies in $R$, and as $R$ is a proper \ah-submodule, this contradicts
the semistability of $Q$. Thus $\dim V_m=12$ for all $m\in M$,
and the proposition is complete.
\end{proof}

The author has not found a satisfactory analogue of this 
proposition for higher dimensions. One can also consider the
problem of reconstructing a hyperk\"ahler manifold from an HP-algebra.
It can be solved easily using a similar approach.

\subsection{Asymptotically conical hyperk\"ahler manifolds}
\label{h35}

In this section we shall define an interesting class of hyperk\"ahler
manifolds, and conjecture some theory about them.

\begin{dfn}
Let $N$ be a compact manifold of dimension $4n-1$, and
set $C=N\times(0,\infty)$. Let $t:C\ra(0,\infty)$ be the
projection to the second factor. Let $v$ be the vector
field $t\,\partial/\partial t$ on $C$. Suppose $C$ has a 
hyperk\"ahler structure, with metric $g$ and complex structures
$I_1,I_2,I_3$. We say that {\it $C$ is a hyperk\"ahler cone} if
$v$ is a Killing vector of $I_1,I_2,I_3$, and 
$g=t^2h+dt^2$, where $h$ is a Riemannian metric on $N$
that is independent of~$t$.
\label{hkconedef}
\end{dfn}

\begin{dfn}
Let $C$ be a hyperk\"ahler cone. Let $f$ be a q-holomorphic function
on $C$, and let $k\ge 0$ be an integer. We say {\it $f$
is homogeneous of degree $k$} if $f=t^kf_N$, where 
$f_N$ is an $\H$-valued function on $N$, independent of $t$.
Define $B^k$ to be the \ah-module of q-holomorphic functions
on $C$ that are homogeneous of degree $k$. Define 
$B=\bigoplus_{k=0}^\infty B^k$. Then $B$ is an \ah-submodule
of the H-algebra $A_C$ of q-holomorphic functions on $C$.
Clearly, $B$ is a graded H-algebra.
\label{bconedef}
\end{dfn}

Here are some basic properties of hyperk\"ahler cones. The proof is
left to the reader.

\begin{lem} Let\/ $C$ be a hyperk\"ahler cone. Then the vector fields 
$v$, $I_1v$, $I_2v$, $I_3v$ generate an action of the Lie algebra 
$\R\op{\go su}(2)$ on $C$, and exponentiating them gives an action 
of the Lie group $\R\t SU(2)$ on $C$. Let\/ $c\in C$, and let\/ 
${\cal O}_c$ be the orbit of\/ $c$ under $\R\t SU(2)$. Then 
${\cal O}_c$ is a hyperk\"ahler submanifold of\/ $C$, and is 
isomorphic as a hyperk\"ahler manifold to $\bigl(\H\setminus\{0\}\bigr)/
\Gamma$, for some finite subgroup~$\Gamma\subset SU(2)$.
\label{conelem}
\end{lem}

\begin{lem} Let\/ $C$ be a hyperk\"ahler cone. Then the graded H-algebra $B$
of Definition \ref{bconedef} is equal to the filtered H-algebra
$P_C$ of q-holomorphic functions of polynomial growth on~$C$.
\label{conepolylem}
\end{lem}

\begin{proof} Let $f$ be a q-holomorphic function on $C$, of
polynomial growth of degree $k$. Let $c\in C$. By Lemma
\ref{conelem}, ${\cal O}_c$ is isomorphic to 
$\bigl(\H\setminus\{0\}\bigr)/\Gamma$, so that the universal
cover $\tilde{\cal O}_c$ is isomorphic to $\H\setminus\{0\}$.
Restricting $f$ to ${\cal O}_c$ and lifting to $\tilde{\cal O}_c$,
the result is a q-holomorphic function on $\H\setminus\{0\}$, 
of polynomial growth of degree~$k$.

Now all such functions are in fact polynomials of degree $k$ 
on $\H$, by a classical result about harmonic functions of 
polynomial growth on $\R^n$. Therefore, we may decompose $f$ 
on $C$ into a sum of homogeneous polynomials on each ${\cal O}_c$. 
Clearly, these homogeneous pieces lie in $B^j$ for $j\le k$, so 
$f\in\bigoplus_{j=0}^kB^j$. We have shown that if $f\in P_C$, 
then $f\in B$, and $P_C\subset B$. But the inclusion $B\subset P_C$ 
is immediate, so $B=P_C$. Clearly, the filtrations on $B$ and
$P_C$ agree, and the lemma is complete.
\end{proof}

Next, we shall define asymptotically conical hyperk\"ahler manifolds.

\begin{dfn} Let $M$ be a complete hyperk\"ahler manifold of dimension $4n$, 
with metric $g^M$ and complex structures $I_1^M,I_2^M,I_3^M$. Let 
$C=N\t(0,\infty)$ be a hyperk\"ahler cone of dimension $4n$, with metric 
$g^C$ and complex structures $I_1^C,I_2^C,I_3^C$. Let $K$ be a 
compact subset of $M$. Then $M\setminus K$ and $N\t(1,\infty)$ are 
submanifolds of $M$ and $C$. Suppose that $\Phi:M\setminus K\ra 
N\t(1,\infty)$ is a diffeomorphism. Let $\nabla$ be the Levi-Civita 
connection on $C$, and $l$ be a positive integer.

We say that {\it $M$ is asymptotically conical}, or {\it AC} 
to order $l$, if 
\begin{gather}
\begin{gathered}
\Phi_*(g^M)=g^C+O(t^{-l}),\qquad
\nabla\bigl(\Phi_*(g^M)\bigr)=O(t^{-l-1}),\\
\nabla^2\bigl(\Phi_*(g^M)\bigr)=O(t^{-l-2}),\qquad\text{and}
\end{gathered}
\label{gacleq}\\
\begin{gathered}
\Phi_*(I_j^M)=I_j^C+O(t^{-l}),\qquad
\nabla\bigl(\Phi_*(I_j^M)\bigr)=O(t^{-l-1}),\\
\nabla^2\bigl(\Phi_*(I_j^M)\bigr)=O(t^{-l-2}),
\qquad\text{for $j=1,2,3$.} 
\end{gathered}
\label{iacleq}
\end{gather}
These equations should be interpreted as follows. Let 
$T$ be a tensor field on $N\t(1,\infty)$. Then $T=O(t^{-k})$ means that
$\md{T}\le \kappa t^{-k}$ on $N\t(1,\infty)$, where $\kappa$ is a 
positive constant, and $\md{\,.\,}$ is taken w.r.t.~the metric $g^C$.
We call $C$ the {\it asymptotic cone of\/~$M$}.
\label{acconedef}
\end{dfn}

Our aim in the remainder of this section is to explore the structure
of the H-algebra of q-holomorphic functions on an AC hyperk\"ahler manifold.
We shall now state -- but not prove, and I do not know a complete proof 
-- a powerful result relating the q-holomorphic functions on the
AC manifold and its asymptotic cone. An {\it incomplete} proof will 
be given shortly, that depends on conjectures to be stated below.

\begin{thm} Let\/ $M$ be a hyperk\"ahler manifold that is AC of order $l$,
with asymptotic cone $C$. Let\/ $B$ be the graded H-algebra of 
q-holomorphic functions on $C$, defined in Definition \ref{bconedef}. 
Let\/ $P$ be the filtered H-algebra of q-holomorphic functions on $M$ 
with polynomial growth. Then $B$ is an SGH-algebra, $P$ is an 
SFH-algebra, and\/ $B$ is isomorphic to $P$ to order $l$, in the sense 
of Definition \ref{asympalgdef}. This implies that $B$ is the associated 
graded H-algebra of\/ $P$, as in Proposition~\ref{asgrhaprop}.
\label{acconethm}
\end{thm}

We shall state two conjectures, and then prove the theorem
assuming these. Before making the conjectures, we shall define 
{\it tensor fields of polynomial growth} on AC hyperk\"ahler manifolds.

\begin{dfn} Let $M$ be an AC hyperk\"ahler manifold, with asymptotic cone
$C$. Let $T$ be a smooth tensor field or function on $M$, and $k$ 
be an integer. Then the expression $\md{T}=O(t^k)$ means that 
$\md{T}\le\kappa\Phi^*(t)^k$ on $M\setminus K$ for some positive 
constant $\kappa$, where $\md{\,.\,}$ is taken w.r.t.~the metric 
$g^M$, and $\Phi,K,t$ and $g^M$ are as in Definition~\ref{acconedef}.
\label{tokdef}
\end{dfn}

With this definition we can state the first conjecture.

\begin{conj} Let\/ $M$ be a hyperk\"ahler manifold that is AC of order 
$l\ge 1$, and let\/ $k\ge -1$ be an integer. Let\/ $x:M\ra\H$ be a smooth 
function, so that\/ $D(x)$ is a 1-form on $M$, where $D$ is the operator 
of\/ \S\ref{h31}. Suppose that\/ $\nabla^aD(x)=O(t^{k-a-1})$ for $a=0,1,2$,
where $\nabla$ is the Levi-Civita connection on $M$. Then there exists 
a smooth function $y:M\ra\H$ such that\/ $D(x)=D(y)$, and\/ $y=O(t^k)$. 
Moreover, if\/ $x$ takes values in $\I$, then $y$ can be chosen to take 
values in~$\I$.
\label{dasympconj}
\end{conj}

This conjecture is a result in analysis, and I believe that it
can be proved using existing mathematical ideas and techniques.
Much work has been done on a similar problem, that of studying
the {\it harmonic} functions of polynomial growth on an
{\it asymptotically flat\/} manifold. Some relevant examples are
\cite[Th.~9.2, p.~76]{LP}, \cite[Th.~1.17, p.~674]{Bar}, and 
in particular the proof of Theorem 3.1 in \cite[p.~678]{Bar}. 
See also \cite{LY}, in which Li and Yau study {\it holomorphic}
functions of subquadratic growth on an asymptotically flat K\"ahler 
manifold.

These results give strong relations between the harmonic functions
of polynomial growth on a Riemannian manifold $M$ asymptotic to 
$\R^n$, and the harmonic polynomials on $\R^n$. Conjecture 
\ref{dasympconj} is modelled on them. I believe that generalizing 
from asymptotically flat to asymptotically conical manifolds should 
be easy. The problems will come from dealing with the operator $D$ 
rather than $\Delta$. One possible tool to use here is Baston's
elliptic complex of operators \cite[p.~43-44]{Bas} resolving~$D$.

\begin{conj} Let $C$ be a hyperk\"ahler cone. Then the graded H-algebra
$B$ of Definition \ref{bconedef} is an SGH-algebra.
\label{coneconj}
\end{conj}

In fact, the author's calculations suggest that hyperk\"ahler cones $C$
always have $B^{2k}\cong\R^a\ot\Sh^kY$, where $a\ge 0$ is an 
integer depending on $k$ and $Y$ is the \ah-module
of \S\ref{h22}, and that most $C$ also have $B^{2k+1}=\{0\}$.
Conjecture \ref{coneconj} would follow immediately from this.
Assuming Conjectures \ref{dasympconj} and \ref{coneconj},
we will now prove Theorem~\ref{acconethm}.
\medskip

\noindent{\it Sketch proof of Theorem \ref{acconethm}.}
Let $j\ge 0$ be an integer. We shall construct a linear map 
$\phi_j:B_j/B_{j-l}\ra P_j/P_{j-l}$. Let $b\in B_j$. Then $b$ is a 
q-holomorphic function on $C$, with polynomial growth of degree $j$. 
With a partition of unity, one may construct a smooth, $\H$-valued 
function $x$ on $M$, such that $\Phi_*(x)=b$ on $N\t(2,\infty)\subset 
C$. Using equations \eq{gacleq} and \eq{iacleq} and the fact that 
$\nabla^ab=O(t^{j-a})$ on $C$ for large $t$ (this is easily proved), it 
can be shown that $\nabla^aD(x)=O(t^{j-l-a-1})$ for $a=0,1,2$ on~$M$.

Putting $k=j-l$, or $k=-1$ if $l>j+1$, we may apply Conjecture
\ref{dasympconj}. It shows that there exists a smooth, $\H$-valued
function $y$ on $M$ with $y=O(t^{j-l})$ or $y=O(t^{-1})$ respectively,
such that $D(x)=D(y)$ on $M$. Therefore, $D(x\!-\!y)=0$, so $x\!-\!y$ is
q-holomorphic on $M$. Clearly, $x\!-\!y$ has polynomial growth of order
$j$, so $x\!-\!y\in P_j$. Now $y$ may not be unique, but since 
$y=O(t^{j-l})$ or $O(t^{-1})$, any two solutions $y$ differ by
an element of $P_{j-l}$ or $P_{-1}=\{0\}$. Therefore, $x\!-\!y+P_{j-l}$ 
is a well-defined element of $P_j/P_{j-l}$, depending only on~$b$.

Define a map $\phi_j:B_j/B_{j-l}\ra P_j/P_{j-l}$ by 
$\phi_j(b+B_{j-l})=x-y+P_{j-l}$. Then $\phi_j$ is a well-defined
$\H$-linear map. Suppose that $b\in B_j'$. Then $b$ takes
values in $\I$. So $x$ takes values in $\I$, and by Conjecture
\ref{dasympconj} we may choose $y$ to take values in $\I$.
Thus $x\!-\!y\in P_j'$, and $\phi_j$ maps $B_j'/B_{j-l}'$ 
to $P_j'/P_{j-l}'$. Therefore $\phi_j$ is an \ah-morphism, 
provided $B_j/B_{j-l}$ and $P_j/P_{j-l}$ are \ah-modules. With a 
little more work, one shows that $\phi_j$ is an \ah-isomorphism.

It remains to verify the conditions of Definition \ref{asympalgdef}.
From the definition of $\phi_j$ it immediately follows that
$\phi_j$ takes $B_k/B_{j-l}$ to $P_k/P_{j-l}$ for $j-l\le k\le j$,
and also that $\phi_j=\phi_{j+1}$ on $B_j/B_{j-l+1}$. By Conjecture 
\ref{coneconj}, $B_j/B_{j-l}$ is a stable \ah-module, so $P_j/P_{j-l}$ 
is also a stable \ah-module. The equation 
$\mu^P_{jkl}\circ(\phi_j\oth\phi_k)=\phi_{j+k}\circ\mu^B_{jkl}$
comes naturally out of the construction of $\phi_j$. Thus $B$
is isomorphic to $P$ to order $l$, by definition. Finally, it 
follows easily that $B$ is the associated graded H-algebra of~$P$.
\hfill \rule{.5em}{.8em} \medskip

\section{Examples, applications and conclusions}

This chapter was difficult to write, because of the many examples, 
little bits of theory, and quaternionic versions of this and that which
begged to be included. For reasons of time and space I have been
ruthless, discussing a few topics only, and not in great depth. 
However, I think that one could easily fill another paper the length
of this one with interesting material.

Section \ref{h41} finds the HP-algebra of q-holomorphic functions of 
polynomial growth on $\H$, in a series of simple steps, as an example. 
In \S\ref{h42} we give examples of how the programme of \S\ref{h34}
(to recover a hypercomplex manifold from its H-algebra) may fail. Then 
\S\ref{h43} shows how to make HL-algebras and HP-algebras out of 
ordinary Lie algebras. This device is applied in \S\ref{h44}, which is 
about hyperk\"ahler manifolds with symmetries. The high point of 
\S\ref{h44} is a (conjectural) algebraic method to explicitly construct 
`coadjoint orbit' hyperk\"ahler manifolds, using HP-algebras.

In \S\ref{h45} we look at at the simplest nontrivial hyperk\"ahler 
manifold --- the Eguchi-Hanson space. It fits into our theory both 
as an AC hyperk\"ahler manifold, and as a `coadjoint orbit'. A careful 
investigation of the H-algebra and its deformations reveals a surprise: 
an unexpected family of singular hypercomplex structures with remarkable 
properties. Section \ref{h46} interprets self-dual connections, or 
`instantons', over a hypercomplex manifold, as modules over its 
H-algebra. Finally, \S\ref{h47} concludes the paper with some 
research problems.

\subsection{Q-holomorphic functions on $\H$}
\label{h41}

Let $\H$ have real coordinates $(x_0,\dots,x_3)$, so that
$(x_0,\dots,x_3)$ represents $x_0+x_1i_1+x_2i_2+x_3i_3$. Now $\H$ 
is naturally a hypercomplex manifold with complex structures given by
$I_1dx_2=dx_3$, $I_2dx_3=dx_1$, $I_3dx_1=dx_2$ and $I_jdx_0=dx_j$, 
for $j=1,2,3$. The study of q-holomorphic functions on $\H$ is called 
{\it quaternionic analysis}, and is surveyed in~\cite{Sud}.

\begin{ex}
\label{hlinex}
First we shall determine the \ah-module $U$ of all linear 
q-holomorphic functions on $\H$. Let $q_0,\dots,q_3\in\H$,
and define $u=q_0x_0+\cdots+q_3x_3$ as an $\H$-valued 
function on $\H$. A calculation shows that $u$ is 
q-holomorphic if and only if 
$q_0+q_1i_1+q_2i_2+q_3i_3=0$. It follows that $U\cong\H^3$.
Also, $U'$ is the vector subspace of $U$ with $q_j\in\I$ for
$j=0,\dots,3$. Let us identify $U$ with $\H^3$ explicitly by 
taking $(q_1,q_2,q_3)$ as quaternionic coordinates. Then
\begin{equation*}
U'=\bigl\{(q_1,q_2,q_3)\in\H^3:\text{$q_j\in\I$ for $j=1,2,3$
and $q_1i_1+q_2i_2+q_3i_3\in\I$}\bigr\}.
\end{equation*}
Thus $U'\cong\R^8$, and $\dim U=4j$, $\dim U'=2j+r$ with
$j=3$ and $r=2$, so the virtual dimension of $U$ is 2. This
is because $\H\cong\C^2$, so the complex dimension of $\H$ is 2.
It is easy to see that $U$ is a stable \ah-module.
\end{ex}

\begin{ex}
\label{hpolex}
Let $k\ge 0$ be an integer, and let $U^{(k)}$ be the \ah-module
of q-holomorphic functions on $\H$ that are homogeneous polynomials
of degree $k$. We shall determine $U^{(k)}$. Write $A$ for the 
H-algebra of q-holomorphic functions on $\H$, and $\mu_A:A\oth A
\ra A$ for the multiplication map. By Example \ref{hlinex}, 
$U^{(1)}=U\subset A$. Thus $\mu_A$ induces an \ah-morphism 
$\mu_A:U\oth U\ra A$, and composing $\mu_A$ $k\!-\!1$ times gives
an \ah-morphism $\mu_A^{k-1}:\bigoth^kU\ra A$. Clearly,
$\Im\mu_A^{k-1}\subset U^{(k)}$. Also, $\mu_A^{k-1}$ is symmetric 
in the $k$ factors of $U$, so it makes sense to restrict to~$\Sh^kU$. 

Thus we have constructed an \ah-morphism $\mu_A^{k-1}:\Sh^kU\ra 
U^{(k)}$. It is easy to show that $\mu_A^{k-1}$ is injective on
$\Sh^kU$. By Example \ref{hlinex}, $U$ is stable with
$j=3$ and $r=2$. Thus Proposition \ref{symmantiprop}
shows that $\dim\Sh^kU=2(k+1)(k+2)$. But Sudbery
\cite[Th.~7, p.~217]{Sud} shows that 
$\dim U^{(k)}=2(k+1)(k+2)$. It follows that $\mu_A^{k-1}$
is an isomorphism, and~$U^{(k)}\cong\Sh^kU$.
\end{ex}

The interpretation of Example \ref{hpolex} is simple. If $V$ 
is the linear polynomials on some vector space, then $S^kV$ 
is the homogeneous polynomials of degree $k$. Here we have a
quaternionic analogue of this, replacing $S^k$ by $\Sh^k$. We have 
found an elegant construction of the spaces $U^{(k)}$, 
important in quaternionic analysis, that gives insight into their 
algebraic structure and dimension.

\begin{ex}
\label{hhalgex}
Let us consider the filtered H-algebra $P$ of q-holomorphic 
functions of polynomial growth on $\H$, as in \S\ref{h32}. 
Clearly, the functions in $U^{(k)}$ have polynomial growth 
of order $k$, so that $U^{(k)}\subset P_k\subset P$. Thus
$\bigoplus_{j=0}^kU^{(j)}\subset P_k$, and 
$\bigoplus_{j=0}^\infty U^{(j)}\subset P$. Now it is a well-known
result in complex analysis that all holomorphic functions
on $\C$ of polynomial growth, are polynomials. The obvious
analogue of this is that all q-holomorphic functions on
$\H$ of polynomial growth are sums of elements of~$U^{(k)}$.

This is in fact true, and can be proved using the theory in
\cite{Sud}. Therefore $P_k=\bigoplus_{j=0}^kU^{(j)}$, and 
$P=\bigoplus_{j=0}^\infty U^{(j)}$. But from Example \ref{hpolex},
$U^{(j)}=\Sh^jU$. Thus $P=\bigoplus_{j=0}^\infty\Sh^jU$.
So, by Definition \ref{freealgdef}, $P$ is isomorphic to the 
free algebra $F^U$ generated by $U$, with its natural filtration. 
As $U$ is finite-dimensional, $P$ is finitely-generated, and in 
particular, $P$ is an FGH-algebra in the sense of~\S\ref{h24}. 

The full H-algebra $A$ of q-holomorphic functions on $\H$ is
obtained by completing $P$, by adding in convergent power series.
The analytic details are beyond the scope of this paper.
Note, however, that because $A$ contains all holomorphic functions
on $\C^2$, $A$ is certainly not finitely-generated, so that $P$
has a much simpler structure than $A$. We may generalize this 
example to the hypercomplex manifold $\H^n$. It is easy to see that the
H-algebra of q-holomorphic functions on $\H^n$ of polynomial
growth is $F^{nU}$, the free H-algebra generated by $n$ copies
of~$U$.
\end{ex}

Now $\H$ is a hyperk\"ahler manifold, so by Theorem \ref{hkhpalg}, $A$ 
and $P$ should be HP-algebras. We shall define the HP-algebra
structure on~$P$.

\begin{ex}
\label{hhpalgex}
We must construct an \ah-morphism $\xi_P:P\oth P\ra P\oth Y$.
From above $U=U^{(1)}\subset P$, so consider $\xi_P:U\oth U\ra
P\oth Y$. Since $\xi_P$ is antisymmetric, we may restrict
to $\Lambh^2U$. Now $U$ is stable and has $j=3$, $r=2$, so
by Proposition \ref{symmantiprop}, we have $\dim\Lambh^2U=8$ 
and $\dim(\Lambh^2U)'=5$. But these are the same dimensions
as those of the \ah-module $Y$ of Definition \ref{ymoddef}.
In fact there is a natural isomorphism $\Lambh^2U\cong Y$. Now 
$U^{(0)}\cong\H$, so that $U^{(0)}\oth Y\cong Y$. Thus we have 
\ah-isomorphisms~$\Lambh^2U^{(1)}\cong U^{(0)}\oth Y\cong Y$.

It is easy to show that the restriction of $\xi_P$ to
$\Lambh^2U^{(1)}$ gives exactly this isomorphism
$\Lambh^2U^{(1)}\cong U^{(0)}\oth Y$. Thus we have defined
$\xi_P$ on a generating subspace $U^{(1)}$ for $P$.
Using Axiom P2, we may extend $\xi_P$ {\it uniquely} to
all of $P$, because the action of $\xi_P$ on the generators
defines the whole action. Now $P$ is a filtered H-algebra,
and $\xi_P$ satisfies $\xi_P(P_j\oth P_k)\subset P_{j+k-2}$.
Thus Axiom PF of \S\ref{h23} holds, and $P$ is a filtered 
HP-algebra.
\end{ex}

\subsection{Hypercomplex manifolds undetermined by their H-algebras}
\label{h42}

In \S\ref{h34} we explained how, under good conditions, it 
is possible to reconstruct a hypercomplex manifold from an H-algebra
of q-holomorphic functions upon it. Here are three examples
where this cannot be done, illustrating different ways in
which the reconstruction can fail.

\begin{ex}
\label{hoverzex}
Since $\Z\subset\R\subset\H$, $\Z$ acts on $\H$ by translation,
and so $\H/\Z\cong\R^3\t{\cal S}^1$ is a hyperk\"ahler manifold. We 
shall determine the filtered H-algebra $P$ of q-holomorphic functions 
on $M=\H/\Z$ of polynomial growth. But $\H$ covers $M$, so
any q-holomorphic function of polynomial growth on $M$
lifts to a q-holomorphic function of polynomial growth on $\H$.
So by Example \ref{hhalgex}, $P$ is the H-subalgebra of $F^U$
that is invariant under~$\Z$.

Clearly, any polynomial invariant under $\Z$ is also invariant
under $\R\subset\H$. Referring to Example \ref{hhalgex}, the 
elements of $U$ invariant under $\R$ are those with $q_0=0$.
Therefore, the $\R$-invariant polynomials in $U$ are $\bigl\{
(q_1,q_2,q_3)\in\H^3:q_1i_1+q_2i_2+q_3i_3=0\bigr\}\subset U$.
But this is the \ah-module $Y$ of \S\ref{h22}. So, one may show 
that the H-algebra $P$ of q-holomorphic functions on $M$ of 
polynomial growth is the free H-algebra~$F^Y$.

Now consider reconstructing $M$ from $P$, as in \S\ref{h34}.
We have $P=F^Y$, so we put $Q=Y$. Then $Q^\d=\I=\R^3$. As 
$P=F^Q$, $M_{P,Q}$ is the whole of $Q^\d=\R^3$. But 
$M\cong\R^3\t{\cal S}^1$, and the map $\pi_{P,Q}:M\ra M_{P,Q}$
is simply the projection to the first factor $\R^3\t{\cal S}^1\ra 
\R^3$. Therefore, in this case, $M_{P,Q}$ is a manifold, and 
$\pi_{P,Q}$ is surjective, but not (even locally) injective. Thus, 
we cannot recover the manifold $M$ from~$P$.
\end{ex}

\begin{ex}
\label{mztex}
Much of Atiyah and Hitchin's book \cite{AH} is an in-depth study 
of a particular complete, noncompact hyperk\"ahler 4-manifold, the 
2-monopole moduli space $M^0_2$. As a manifold, $M^0_2$ is diffeomorphic 
to a complex line bundle over $\mathbb{RP}^2$, and the zero section 
gives a submanifold $\mathbb{RP}^2\subset M^0_2$. The isometry group 
of the metric $g$ on $M^0_2$ is $SO(3)$. This isometry group acts in a 
{\it nontrivial} way on the hyperk\"ahler structure. There is a vector 
space isomorphism $\langle I_1,I_2,I_3\rangle\cong\mathfrak{so}(3)$, 
that identifies the action of $SO(3)$ on $\langle I_1,I_2,I_3\rangle$ 
with the adjoint action on~$\mathfrak{so}(3)$.

Because $SO(3)$ does not fix the hyperk\"ahler structure, there are no
hyperk\"ahler moment maps (see \S\ref{h44}). However, there are some 
K\"ahler moment maps. Let $v_1,v_2,v_3\in\mathfrak{so}(3)$ be identified 
with $I_1,I_2,I_3$ under the isomorphism above. Then $v_1,v_2,v_3$ are 
Killing vectors of $g$ on $M$. Let $\alpha_1,\alpha_2,\alpha_3$ be 
the 1-forms on $M^0_2$ dual to $v_1,v_2,v_3$ under $g$. Then a brief 
calculation shows that the 1-form $I_j\alpha_k+I_k\alpha_j$ is closed, 
for $j,k=1,\dots,3$. Since $b^1(M^0_2)=0$, there exists a unique 
real function $f_{jk}$ with $df_{jk}=I_j\alpha_k+I_k\alpha_j$,
and such that~$\int_{\mathbb{RP}^2}f_{jk}dA=0$.

The $f_{jk}$ form a vector space $S^2\R^3\cong\R^6$ of real
functions on $M^0_2$. For $j,k=1,2,3$, let $q_{jk}\in\H$ with
$q_{jk}=q_{kj}$, and define $z=\Sigma_{j,k=1}^3q_{jk}f_{jk}$. 
When is $z$ a q-holomorphic function on $M^0_2$? Calculation
shows that $z$ is q-holomorphic if and only if $\Sigma_{j=1}^3
q_{jk}i_j=0$ for $k=1,2,3$. Therefore, the \ah-module of
q-holomorphic $z$ is
\begin{equation}
Z=\bigl\langle f_{22}-f_{33}+2i_1f_{23},
f_{33}-f_{11}+2i_2f_{31},
f_{11}-f_{22}+2i_3f_{12}\rangle.
\label{zdefeq}
\end{equation}
Note that only the trace-free part $S^2_0\mathfrak{so}(3)\cong\R^5$
appears here. From \eq{zdefeq} we find that $Z\cong\H^3$ and 
$Z'\cong\R^7$, and in fact there is a canonical isomorphism 
$Z\cong Y\oth Y$, where $Y$ is the \ah-module defined in~\S\ref{h22}.

We have found an \ah-module $Z\cong Y\oth Y$ of q-holomorphic
functions on $M^0_2$. By Lemma \ref{genlem} there is an
H-algebra morphism $\phi_Z:F^Z\ra A$, where $A$ is the H-algebra
of q-holomorphic functions on $M^0_2$. Now Atiyah and Hitchin
\cite{AH} define the metric on $M^0_2$ explicitly, using an
elliptic integral. Their construction uses transcendental functions,
not just algebraic functions. Because of this, it can be shown 
that the 5 functions $f_{jk}$ used in \eq{zdefeq} are algebraically 
independent. Therefore, there can be no polynomial relations 
between them, so that $\phi_Z$ is injective.

Thus $F^Z$ is an H-algebra of q-holomorphic functions on $M^0_2$.
Now $Z^\d\cong\R^5$, and since $F^Z$ is free, $M_{F^Z,Z}$ is the
whole of $Z^\d$. Thus $\pi_{F^Z,Z}$ maps $M^0_2$ to $\R^5$. It can
be shown that $\pi_{F^Z,Z}$ is generically injective, but because
of the dimensions $\pi_{F^Z,Z}$ cannot be surjective. Thus we
cannot recover the manifold $M^0_2$ from $F^Z$. I claim that
$F^Z$ is in fact the whole H-algebra of q-holomorphic functions 
of polynomial growth on~$M^0_2$.
\end{ex}

\begin{ex}
\label{strex}
Let $k\ge 2$, and let $M=\H^k$. Then by Example \ref{hhalgex},
the H-algebra of q-holomorphic functions of polynomial growth
on $M$ is $F^{kU}$. Now $U\cong\H^3$, so $kU\cong\H^{3k}$.
It can be shown that if $k\ge 2$, then the {\it generic} 
$\H$-submodule $Q\cong\H^{3k-1}$ of $kU$ is a stable \ah-module 
with $Q'\cong\R^{8k-4}$ and $Q^\d=\R^{4k}$. Let $Q$ be such an 
\ah-submodule, and let $P=F^Q\subset F^U$. Then $P$ is an
H-algebra of q-holomorphic functions on~$M=\H^k$.

Consider reconstructing $M$ and its hypercomplex structure from $P$, as
in \S\ref{h34}. Since $P$ is free, $M_{P,Q}$ is the whole
of $Q^\d$, which is $\R^{4k}$. But $M\cong\R^{4k}$, and in fact
$\pi_{P,Q}:M\ra M_{P,Q}$ is a diffeomorphism, the identity.
However, in Lemma \ref{recondlem} we have $V_m\cong Q$ for each 
$m\in M$, so that $\dim V_m=12k-4$. Thus, the lemma shows that the 
hypercomplex structure of $M$ cannot be recovered from $P$. This means 
that $P$ gives full information about the manifold $M$, but only partial 
information about the hypercomplex structure.
\end{ex}

\subsection{HL-algebras and HP-algebras}
\label{h43}

Here is a simple construction of HL-algebras.

\begin{ex}
\label{liealgex}
Let $\g$ be a Lie algebra. Then the Lie bracket $[\,,\,]$ 
on $\g$ gives a linear map $\lambda:\g\ot\g\ra\g$, such
that $\lambda(x\ot y)=[x,y]$ for $x,y\in\g$. Let $Y$ be
the \ah-module defined in \S\ref{h22}, and define 
$A_\g$ to be the \ah-module $\g\ot Y$. Then $A_\g\oth A_\g
\cong (\g\ot\g)\ot(Y\oth Y)$ and $A_\g\oth Y\cong\g\ot(Y\oth Y)$.
Define a linear map $\xi_{A_\g}:A_\g\oth A_\g\ra A_\g\oth Y$ by
$\xi_{A_\g}=\lambda\ot\id$, as a map from $(\g\ot\g)\ot(Y\oth Y)$ 
to $\g\ot(Y\oth Y)$. It is easy to show that $\xi_{A_\g}$
satisfies Axiom P1 of \S\ref{h22}, using the Jacobi identity
for $\g$ to prove part $(iii)$. Therefore, $A_\g$ is an 
HL-algebra, by Definition~\ref{hpldef}.
\end{ex}

\begin{ex}
\label{findex}
Example \ref{hhpalgex} constructed a filtered HP-algebra $P$. 
Because $\xi_P(P_j\oth P_k)\subset P_{j+k-2}$, each of $P_0,P_1$ 
and $P_2$ are closed under $\xi_P$. Therefore, restricting $\xi_P$ 
to $P_0,P_1$ and $P_2$ gives them the structure of HL-algebras.
Now $P_0\cong\H$ and $\xi_P$ is zero on $P_0$, but
$P_1\cong\H\op U$, $P_2\cong\H\op U\op\Sh^2U$, which
are both finite-dimensional, and $\xi_P$ is nontrivial on
both. It is easy to see that none of these is of the
form $\g\ot Y$. Thus, there exist nontrivial, 
finite-dimensional HL-algebras that do not arise from
Example~\ref{liealgex}.
\end{ex}

The next example is an aside about Lie and Poisson algebras.

\begin{ex}
\label{gpoalgex}
Let $\g$ be a Lie algebra. The {\it symmetric algebra} 
$S(\g)=\bigoplus_{k=0}^\infty S^k\g$ of $\g$ is a free, 
commutative algebra generated by $\g$. Define a bracket 
$\{\,,\,\}:S(\g)\t S(\g)\ra S(\g)$ as follows. Let 
$k,l\ge 0$ be integers. If $x,y\in\g$, then $x^k\in S^k\g$ 
and $y^l\in S^l\g$, and $S^k\g,S^l\g$ are generated by 
such elements. When $k=0$ or $l=0$ define $\{\,,\,\}=0$ on
$S^k\g\times S^l\g$, and when $k,l>0$ define $\{x^k,y^l\}=
kl\,\sigma\bigl(x^{k-1}\ot[x,y]\ot y^{l-1}\bigr)$. Here 
$[\,,\,]$ is the Lie bracket on $\g$ and $\sigma:S^{k-1}\g
\ot\g\ot S^{l-1}\g\ra S^{k+l-1}\g$ is the symmetrization 
operator, a projection. 

This definition extends uniquely to give a bilinear operator 
$\{\,,\,\}:S^k\g\t S^l\g\ra S^{k+l-1}\g$, so we have found 
a bilinear bracket $\{\,,\,\}$ on $S(\g)\t S(\g)$. It is
shown in \cite[\S 1.4]{BV} that $\{\,,\,\}$ is a Poisson
bracket on $S(\g)$, so that $S(\g)$ is a Poisson algebra,
the Poisson algebra of the Lie algebra~$\g$.
\end{ex}

Using Example \ref{gpoalgex} as a model, here is a 
construction of HP-algebras.

\begin{ex}
\label{unienvex}
Let $A$ be an HL-algebra. Then \S\ref{h24} defines the
free H-algebra $F^A$ generated by $A$. To make $F^A$ into
an HP-algebra, we must give a Poisson bracket $\xi_{F^A}$
on $F^A$. Let $k,l$ be positive integers, and define a map 
$\xi_{k,l}:\Sh^kA\oth\Sh^lA\ra\Sh^{k+l-1}A$ by $\xi_{k,l}=
kl\,\sigma_\H\circ(\id\oth\xi_A\oth\id)\circ\iota$, using
the sequence of maps
\begin{align*}
&\Sh^kA\oth \Sh^lA{\buildrel\iota\over\ra}
\Sh^{k-1}A\oth A\oth A\oth\Sh^{l-1}A\,\,
{\buildrel\id\smalloth\xi_A\smalloth\id\over\longra}\\
&\Sh^{k-1}A\oth A\oth Y\oth\Sh^{l-1}A
{\buildrel\sigma_\H\over\longra}
\Sh^{k+l-1}A\oth Y.
\end{align*}

Here $\iota$ is the inclusion, and $\sigma_\H$ the 
symmetrization operator of \S\ref{h11}. For $k=0$ or $l=0$, 
let $\xi_{k,l}=0$. Define $\xi_{F^A}:F^A\oth F^A\ra F^A\oth Y$ to 
be the unique linear map such that the restriction of $\xi_{F^A}$ 
to $\Sh^kA\oth\Sh^lA$ is $\xi_{k,l}$, for all $k,l\ge 0$. Then 
$\xi_{F^A}$ is a well-defined \ah-morphism. A calculation 
following those in \cite[\S 1.4]{BV} shows that $\xi_{F^A}$ 
satisfies Axioms P1 and P2 of \S\ref{h22}. Thus, by Definition 
\ref{hpldef}, $F^A$ is an HP-algebra. Moreover, $F^A$ is a
filtered H-algebra with the natural filtration, and it is easy 
to show that Axiom PF holds. So $F^A$ is a filtered HP-algebra.
\end{ex}

Combining Examples \ref{liealgex} and \ref{unienvex},
we see that if $\g$ is a Lie algebra, then $F^{A_\g}$ is
a filtered HP-algebra. 

\subsection{Hyperk\"ahler manifolds with symmetries}
\label{h44}

Let $M$ be a hyperk\"ahler manifold, and suppose $v$ is a Killing 
vector of the hyperk\"ahler structure on $M$. A {\it hyperk\"ahler 
moment map for $v$} is a triple $(f_1,f_2,f_3)$ of smooth real 
functions on $M$ such that $\alpha=I_1df_1=I_2df_2=I_3df_3$, where
$\alpha$ is the 1-form dual to $v$ under the metric $g$.
Moment maps always exist if $b^1(M)=0$, and are unique up
to additive constants.

More generally, let $M$ be a hyperk\"ahler manifold, let $G$ be a Lie
group with Lie algebra $\g$, and suppose $\Phi:G\ra\Aut(M)$
is a homomorphism from $G$ to the group of automorphisms
of the hyperk\"ahler structure on $M$. Let $\phi:\g\ra\Vect(M)$ be 
the induced map from $\g$ to the Killing vectors.
Then a {\it hyperk\"ahler moment map for the action $\Phi$ of\/ $G$} is 
a triple $(f_1,f_2,f_3)$ of smooth functions from $M$ to $\g^*$,
such that for each $x\in\g$, $(x\!\cdot\!f_1,x\!\cdot\!f_2,
x\!\cdot\!f_3)$ is a hyperk\"ahler moment map for the vector field 
$\phi(x)$, and in addition, $(f_1,f_2,f_3)$ is equivariant 
under the action $\Phi$ of $G$ on $M$ and the coadjoint action 
of $G$ on~$\g^*$.

Moment maps are a familiar part of symplectic geometry,
and hyperk\"ahler moment maps were introduced by Hitchin et 
al.~as part of a quotient construction for hyperk\"ahler manifolds
\cite{HKLR}, \cite[p.~118-122]{Sal}. Hyperk\"ahler moment maps
will exist under quite mild conditions on $M$ and $G$,
for instance if $b^1(M)=0$ and $G$ is compact. We shall
use them to construct q-holomorphic functions on hyperk\"ahler
manifolds with symmetries.

\begin{ex}
\label{vecmomex}
Let $M$ be a hyperk\"ahler manifold and $v$ a nonzero Killing vector 
of the hyperk\"ahler structure on $M$. Suppose $(f_1,f_2,f_3)$ is a 
hyperk\"ahler moment map for $v$. We shall make q-holomorphic functions 
on $M$ out of the real functions $f_j$. Let $q_1,q_2,q_3\in\H$, 
and consider the $\H$-valued function $y=q_1f_1+q_2f_2+q_3f_3$ 
on $M$. By construction the $f_j$ satisfy $I_1df_1=I_2df_2=
I_3df_3$. Using this equation, the fact that $v$ is nonzero,
and the definition of q-holomorphic in \S\ref{h31}, it is easy 
to show that $y$ is q-holomorphic if and only if $q_1i_1+q_2i_2
+q_3i_3=0$. Thus we have constructed an \ah-module 
$Y\cong\bigl\{(q_1,q_2,q_3)\in\H^3:q_1i_1+q_2i_2+q_3i_3=0\bigr\}$ 
of q-holomorphic functions on $M$. It is isomorphic to the 
\ah-module $Y$ of Definition~\ref{ymoddef}.
\end{ex}

\begin{ex}
\label{momex}
Now let $M$ be a hyperk\"ahler manifold, $G$ a Lie group, 
$\Phi:G\ra\Aut(M)$ an action of $G$ on $M$ preserving the hyperk\"ahler 
structure, and $\phi:\g\ra\Vect(M)$ the induced map. Suppose 
$\phi$ is injective, and that $(f_1,f_2,f_3)$ is a hyperk\"ahler moment 
map for $\Phi$. By Example \ref{vecmomex}, each nonzero $x\in\g$ 
gives us an \ah-module $Y$ of q-holomorphic functions on $M$. 
Clearly, these fit together to form a canonical \ah-module 
$\g\ot Y$ of q-holomorphic functions on~$M$.

We have already met $\g\ot Y=A_\g$ as an HL-algebra in Example
\ref{liealgex}, where we defined a Poisson bracket $\xi_{A_\g}$
on it. In the context of this example, $\g\ot Y$ is an
\ah-submodule of the HP-algebra $A_M$ of q-holomorphic functions
on $M$, which derives its own Poisson bracket $\xi_{A_M}$ from
the hyperk\"ahler structure of $M$. A computation shows that $\g\ot Y$ 
is closed under $\xi_{A_M}$, and that~$\xi_{A_M}=\xi_{A_\g}$.
\end{ex}

Example \ref{momex} shows that the HL-algebras of Example 
\ref{liealgex} are related to the HP-algebras of hyperk\"ahler manifolds 
with symmetry groups. In the following lemma we extend this
to the associated HP-algebras defined by Example~\ref{unienvex}.

\begin{lem} Let\/ $M$ be a hyperk\"ahler manifold, and\/ $A_M$ the
HP-algebra of q-holomorphic functions on $M$. Let\/ $G$
be a Lie group, $\Phi:G\ra\Aut(M)$ be an action of $G$ on $M$
preserving the hyperk\"ahler structure, and $\phi:\g\ra\Vect(M)$ be the 
induced map. Suppose $(f_1,f_2,f_3)$ is a hyperk\"ahler moment map for 
$\Phi$. Then there is a canonical HP-algebra morphism 
$\Phi_*:F^{A_\g}\ra A_M$, where $F^{A_\g}$ is defined by 
Examples \ref{liealgex} and~\ref{unienvex}.
\label{hpmorlem}
\end{lem}

\begin{proof} In Example \ref{momex} we constructed an 
\ah-submodule $A_\g=\g\ot Y$ of $A_M$. By Lemma \ref{genlem},
there is a unique H-algebra morphism $\Phi_*=\phi_{A_\g}:
F^{A_\g}\ra A_M$. To complete the proof we must show that
$\Phi_*$ is an HP-algebra morphism, where the HP-algebra
structure on $F^{A_\g}$ is defined by Example \ref{unienvex},
so we must show that $\Phi_*$ identifies the Poisson brackets
on $F^{A_\g}$ and~$A_M$.

In Example \ref{unienvex} we remarked that $\xi_{A_M}=\xi_{A_\g}$
on $A_\g$. Thus $\Phi_*$ identifies the Poisson brackets
on this subspace. Because $\Phi_*$ is an H-algebra morphism
and $A_\g$ generates $F^{A_\g}$, we can deduce from Axiom P2 
that if $\xi_{A_\g}$ and the pullback of $\xi_{A_M}$ agree
on $A_\g$, they must agree on the whole of $F^{A_\g}$.
Thus $\Phi_*$ identifies the Poisson brackets of $F^{A_\g}$
and $A_M$, and $\Phi_*$ is an HP-algebra morphism.
\end{proof}

In the situation of the lemma, $\Phi_*(F^{A_\g})$ is an
HP-algebra containing information about $M$ and $G$. This 
suggests that to understand hyperk\"ahler manifolds with symmetries
better, it may be helpful to study HP-algebras that are 
images of $F^{A_\g}$. The next two examples construct such
images.

\begin{ex} 
\label{krcoadex}
Let $G$ be a Lie group with Lie algebra $\g$. Then $A_\g=\g\ot Y$ 
and $Y^\d=\I$ by definition, so $A_\g^\d=\g^*\ot\I$. As in the 
proof of Lemma \ref{algvarlem}, each element $f$ of $F^{A_\g}$ 
induces an $\H$-valued polynomial $\psi_f$ on $A_\g^\d=\g^*\ot\I$. 
Now $G$ acts on $\g^*$ by the coadjoint action, so $G$ acts on 
$\g^*\ot\I$, with trivial action on $\I$. Let $\Omega\subset
\g^*\ot\I$ be an orbit of $G$, that is contained in no proper 
vector subspace of $\g^*\ot\I$. Define an \ah-submodule 
$I^\Omega$ in $F^{A_\g}$ by 
\begin{equation}
I^\Omega=\bigl\{f\in F^{A_\g}:
\text{$\psi_f\equiv 0$ on $\Omega$}\bigr\}.
\label{iomegeq}
\end{equation}
Clearly, $\mu_{F^{A_\g}}\bigl(I^\Omega\oth F^{A_\g}\bigr)\subset 
I^\Omega$, so $I^\Omega$ is an {\it ideal}.

Define a subset $M^\Omega\subset\g^*\ot\I$ by
\begin{equation}
M^\Omega=\bigl\{x\in \g^*\ot\I:
\text{$\psi_f(x)=0$ for all $f\in I^\Omega$}\bigr\}.
\label{momegeq}
\end{equation}
Then $M^\Omega$ is an affine real algebraic variety in
$\g^*\ot\I$, as in Lemma \ref{algvarlem}. Also, $\Omega\subset 
M^\Omega$ by \eq{iomegeq}, and $M^\Omega$ is invariant under
the action of $G$ on $\g^*\ot\I$. Define $A^\Omega=\bigl\{\psi_f
\vert_{M^\Omega}:f\in F^{A_\g}\bigr\}$. Thus $A^\Omega$ is a vector 
space of $\H$-valued functions on~$M^\Omega$. 

Now $I^\Omega$, $F^{A_\g}$ and $A^\Omega$ fit into a natural, 
\ah-exact sequence 
\begin{equation}
0\ra I^\Omega{\buildrel\iota\over\longra} 
F^{A_\g}{\buildrel\rho\over\longra}A^\Omega\ra 0,
\label{omegseqeq}
\end{equation}
where $\iota$ is the inclusion map, and $\rho$ is the restriction
map from $\g^*\ot\I$ to $M^\Omega$. Since $I^\Omega$ is an
ideal, intuitively $A^\Omega$ should be an H-algebra. 

In \S\ref{h24} we saw that the best sort of ideal is a {\it stable 
filtered ideal}. Therefore, let us assume that $I^\Omega$ is a
stable filtered ideal. Then Lemma \ref{ideallem} shows that
$A^\Omega$ is an SFH-algebra, and $\rho$ a filtered H-algebra 
morphism. Also, because we suppose that $\Omega$ is not contained 
in any proper subspace of $\g^*\ot\I$, and $\Omega\subset M^\Omega$,
we see that $\rho\vert_{\g\ot Y}$ is injective, and $\g\ot Y$
is an \ah-submodule of $A^\Omega$ that generates $A^\Omega$.
It can be shown that $M_{A^\Omega,\g\ot Y}=M^\Omega$. Observe 
that $A^\Omega$ is an FGH-algebra in the sense of~\S\ref{h24}.
\end{ex}

\begin{ex}
\label{hpcoadex}
In the situation of the previous example, we shall show that 
$A^\Omega$ is an HP-algebra. The linear map $\xi_{F^{A_\g}}:
F^{A_\g}\oth A_\g\ra F^{A_\g}\oth Y$ is easy to understand, 
because $A_\g=\g\ot Y$, so we may write the map as 
$\xi_{F^{A_\g}}:\g\ot(F^{A_\g}\oth Y)\ra F^{A_\g}\oth Y$, 
which just gives the Lie algebra action of $\g$ on $F^{A_\g}$. 
Since $I^\Omega$ is $G$-invariant, it follows that $\xi_{F^{A_\g}}$ 
maps $I^\Omega\oth A_\g$ to $I^\Omega\oth Y$. But $A_\g$ generates 
$F^{A_\g}$, and so using Axiom P2 and the fact that 
$\mu_{F^{A_\g}}\bigl(I^\Omega\oth F^{A_\g}\bigr)\subset I^\Omega$, 
it follows that~$\xi_{F^{A_\g}}\bigl(I^\Omega\oth F^{A_\g}\bigr)
\subset I^\Omega\oth Y$.

Using the assumption in Example \ref{krcoadex}, this inclusion
is just what is needed to prove that the Poisson bracket 
$\xi_{F^{A_\g}}$ can be pushed down to $A^\Omega$ using $\rho$,
inducing a Poisson bracket $\xi_{A^\Omega}$ on $A^\Omega$, so that 
$A^\Omega$ is an HP-algebra. As $\rho$ is a filtered H-algebra 
morphism, $A^\Omega$ is a filtered HP-algebra, and $\rho$ a filtered 
HP-algebra morphism.
\end{ex}

Examples \ref{krcoadex} and \ref{hpcoadex} construct a large
family of HP-algebras $A^\Omega$ associated to Lie groups. As in 
\S\ref{h34}, we can try to use $A^\Omega$ to construct a hyperk\"ahler
structure on $M^\Omega$. This suggests that associated to
each Lie group $G$, there is a natural family of hyperk\"ahler manifolds.
Now Kronheimer \cite{Kr3,Kr4}, Biquard \cite{Bi} and 
Kovalev \cite{Ko} have also constructed hyperk\"ahler manifolds associated 
to Lie groups, from a completely different point of view.

Let $G$ be a compact Lie group with Lie algebra $\g$, and let 
the complexification of $G$ be $G^c$ with Lie algebra $\g^c$. 
Kronheimer found that certain moduli spaces of singular 
$G$-instantons on $\R^4$ are hyperk\"ahler manifolds. These moduli 
spaces can be identified with coadjoint orbits of $G^c$ in $(\g^c)^*$,
and have hyperk\"ahler metrics invariant under $G$. Kronheimer's 
construction worked only for certain special coadjoint orbits, and 
more general cases were handled by Biquard and Kovalev.

Although these metrics look very algebraic, their construction 
is in fact analytic, and the algebraic description of these metrics 
is not well understood. I propose that Examples \ref{krcoadex} and
\ref{hpcoadex} provide this algebraic description. This was proved
in \cite[\S 11-12]{Joy2} for Kronheimer's metrics \cite{Kr3,Kr4}, 
which are the simplest case, but I have not yet proved it for the 
metrics of Biquard and Kovalev. Here is a conjecture about this.

\begin{conj} 
\label{krconj}
We conjecture the following relations between Examples 
\ref{krcoadex} and \ref{hpcoadex}, and Kronheimer,
Biquard and Kovalev's `coadjoint orbit' metrics.
\begin{itemize}
\item The assumption made in Example \ref{krcoadex} always holds.
\item Using the techniques of \S\ref{h34}, the HP-algebra
$A^\Omega$ determines a hyperk\"ahler structure on a dense open set 
of~$M^\Omega$.
\item These hyperk\"ahler structures on subsets of $M^\Omega$ include
all those of \cite{Kr3}, \cite{Kr4}, \cite{Bi} and \cite{Ko} 
as special cases. However, generically the structures on 
$M^\Omega$ are new, and do not coincide with those of 
Kronheimer, Biquard and Kovalev.
\item For some $\Omega$, $M^\Omega$ is a cone in $\g^*\ot\I$,
and $A^\Omega$ is an SGH-algebra. Then $M^\Omega$ is a hyperk\"ahler cone,
as in \S\ref{h35}. These are the {\it nilpotent orbits} 
of~\cite{Kr4}.
\item For each $\Omega$ there is an $\tilde\Omega$, such that
$M^{\tilde\Omega}$ is a `nilpotent orbit', and the associated
graded H-algebra of $A^\Omega$ is $A^{\tilde\Omega}$. Then
the (singular) hyperk\"ahler manifold $M^\Omega$ has an AC end, as in
Definition \ref{acconedef}, with asymptotic cone~$M^{\tilde\Omega}$.
\end{itemize}
\end{conj}

Further study of these spaces from the HP-algebra point of view
will probably lead to a much clearer understanding of the algebra
and geometry underlying Kronheimer's metrics. We will look at an
example in detail in the next section.

\subsection{The Eguchi-Hanson space}
\label{h45}

The Eguchi-Hanson space $M$ is a noncompact hyperk\"ahler manifold 
of dimension 4, with a metric written down by Eguchi and Hanson 
\cite{EH}. It is of interest to us as an example for two reasons. 
Firstly, it is the simplest interesting example of a `coadjoint orbit' 
metric, as it fits into the theory of \S\ref{h44} with $G=SU(2)$ or 
$SO(3)$. Secondly, it is the simplest hyperk\"ahler {\it asymptotically 
locally Euclidean} or {\it ALE space}. An ALE space is an AC hyperk\"ahler 
manifold of dimension 4, with asymptotic cone $\H/\Gamma$, for some finite 
subgroup~$\Gamma\subset SU(2)$. 

The ALE spaces for cyclic $\Gamma$ were described explicitly by Gibbons 
and Hawking \cite{GH} and Hitchin \cite{Hi}, and a complete construction 
and classification of ALE spaces was given by Kronheimer \cite{Kr1}, 
\cite{Kr2}. The Eguchi-Hanson space has asymptotic cone $\H/\{\pm1\}$,
and as a manifold it is diffeomorphic to the total space of
$T^*\mathbb{CP}^1$. In this section we will treat the Eguchi-Hanson space 
from our H-algebra point of view. It gives us our first explicit example 
of the theories of \S\ref{h35} and \S\ref{h44}. On the way, we will
determine the deformations of the H-algebra of the Eguchi-Hanson space. 
The result and its implications are quite surprising.

Let $M$ be the Eguchi-Hanson space. Then the hyperk\"ahler cone of $M$, 
in the sense of \S\ref{h35}, is $\H/\{\pm1\}$. Our first step is to
find the SGH-algebra $B$ of q-holomorphic functions on $\H/\{\pm1\}$
with polynomial growth.

\begin{ex}
\label{acconex}
By Example \ref{hhalgex}, the graded H-algebra of q-holomorphic
functions on $\H$ is $F^U$, where $U$ is defined in the example.
Define $\sigma:\H\ra\H$ by $\sigma(q)=-q$. Then $\sigma$ acts on 
$F^U$, and as the functions in $U$ are linear, $\sigma$ acts as $-1$ 
on $U$. Therefore $\sigma$ acts as $(-1)^j$ on $\Sh^jU$. Now the 
SGH-algebra $B$ of q-holomorphic functions on $\H/\{\pm1\}$ is just 
the $\sigma$-invariant part of $F^U$. Therefore 
$B=\bigoplus_{j=0}^\infty\Sh^{2j}U$. For compatibility with \S\ref{h35} 
we adopt the grading $B^{2j}=\Sh^{2j}U$ and $B^{2j+1}=\{0\}$,
and the corresponding filtration.

A calculation using Proposition \ref{symmantiprop} and
the definitions of $Y$ and $U$ in \S\S\ref{h22} and \ref{h41} shows
that $\dim(\Sh^jY)=4k$, $\dim(\Sh^jY)'=2k\!+\!s$ with $k=j\!+\!1$ and
$s=1$, and $\dim(\Sh^{2j}U)=4l$, $\dim(\Sh^{2j}U)'=2l+t$ with 
$l=(2j\!+\!1)(j\!+\!1)$ and $t=2j\!+\!1$. Thus $l=(2j\!+\!1)k$ and 
$t=(2j\!+\!1)s$. In fact, by carefully investigating the geometry it 
can be shown that~$\Sh^{2j}U\cong\R^{2j+1}\ot\Sh^jY$. 

Here is one way to see this. It is well-known that $SO(4)$ acts 
irreducibly on $\R^4$, but that $\R^4\ot\C$ splits as $\C^2\ot\C^2$. 
In the same way, $U$ is irreducible, but $U\ot\C$ may be written as 
$\C^2\ot W$, where $W$ is a {\it complex} \ah-module with 
$W\oth W=\Sh^2W\cong Y\ot\C$ as complex \ah-modules. Thus we may 
write $B^{2j}\cong\R^{2j+1}\ot\Sh^jY$, where $\R^{2j+1}$ is the real 
part of $\C^{2j+1}=S^j(\C^2)$, and~$B=\bigoplus_{j=0}^\infty
\R^{2j+1}\ot\Sh^jY$.
\end{ex}

Next, we present $B$ as an FGH-algebra generated by $Q\subset B$,
and determine the algebraic variety~$M_{B,Q}$.

\begin{ex}
\label{bfghex}
To write $B$ as an FGH-algebra, the quotient of a free H-algebra 
$F^Q$ by an ideal $I\subset F^Q$, one must choose a subspace $Q$ of 
$B$ that generates $B$. Define $Q=B^2=\R^3\ot Y$, where $\R^3$
is equipped with a Euclidean metric and inner product. Then $\phi_Q:
F^Q\ra B$ is defined in Lemma \ref{genlem}, and we set $I=\Ker\phi_Q$.
Since $B^{2j}=\R^{2j+1}\ot\Sh^jY$ and $\Sh^jQ=\R^{(j+1)(j+2)/2}
\ot\Sh^jY$, it can be shown that $\phi_Q:\Sh^jQ\ra B^{2j}$ is 
surjective for all $j\ge 0$, and has kernel $I^{2j}\subset\Sh^jQ$ 
with~$I^{2j}\cong\R^{j(j-1)/2}\ot\Sh^jY$. 

Now $\Sh^jY$ is stable by Proposition \ref{symmantiprop}, as $Y$ is 
stable. Thus $I$ is a stable filtered ideal in the sense of Definition 
\ref{idealdef}. Therefore $B$ is generated by $Q$, and $B$ is an 
FGH-algebra by Definition \ref{fghdef}. Set $J=I^4\cong\Sh^2Y$. Another 
calculation shows that $I$ is generated by $J$ in the sense of 
Definition \ref{idealdef}. Moreover, if $h$ is the Euclidean metric 
on $\R^3$, then the \ah-module $J\subset\Sh^2Q$ is~$J=\langle h\rangle
\ot\Sh^2Y\subset S^2\R^3\ot\Sh^2Y=\Sh^2Q$.

Thus $B$ is the quotient of the free SFH-algebra $F^Q$, where 
$Q=\R^3\ot Y$, by the stable filtered ideal $I$ in $F^Q$, where $I$ is 
generated by the \ah-module $J=\Sh^2Y\subset\Sh^2Q$. This is an explicit
description of $B$ as an FGH-algebra. To finish the example, we will
find and describe the subset $M_{B,Q}$ of $Q^\d$, defined in \S\ref{h34}.
As $Q=\R^3\ot Y$, $Q^\d=(\R^3)^*\ot\I$, as $Y^\d\cong\I$. Identify 
$(\R^3)^*$ with $\R^3$ using the metric. Then a point $\gamma$
in $Q^\d$ may be written $\gamma=\Sigma_{k=1}^3v_k\ot i_k$, where
$v_1,v_2,v_3$ are vectors in~$\R^3$.

The proof of Lemma \ref{algvarlem} showed that each element $y$ 
of $F^Q$ defines an $\H$-valued polynomial $\psi_y$ on $Q^\d$. In our
case, as $J$ generates $I$, $M_{B,Q}$ is the zeros of the polynomials
$\psi_y$ for $y\in J$. Moreover, as $J\subset\Sh^2Q$, these 
polynomials are homogeneous quadratics on $Q^\d$. A computation using
the definitions of $Y$, $Q$ and $J$ shows that $\gamma\in M_{B,Q}$ 
if and only if
\begin{equation}
v_1\cdot v_1=v_2\cdot v_2=v_3\cdot v_3, \quad
v_1\cdot v_2=0, \quad v_2\cdot v_3=0, \quad v_3\cdot v_1=0,
\label{vdoteq}
\end{equation}
where `$\cdot$' is the inner product on $\R^3$. So $v_1,v_2,v_3$ 
must be orthogonal in $\R^3$, and of equal length.
\end{ex}

Now we look at $M_{B,Q}$ more closely, and interpret it as
a `coadjoint orbit'.

\begin{ex}
\label{vlgex}
The previous example identified the algebraic variety $M_{B,Q}$ 
with the set of triples $(v_1,v_2,v_3)$ of vectors in $\R^3$ that are 
orthogonal and of equal length. For each $r>0$, define $S_{r,+}$ to be 
the set of orthogonal triples $(v_1,v_2,v_3)$ that form a positively 
oriented basis of $\R^3$, and for which $\md{v_1}=\md{v_2}=\md{v_3}=r$.
Let $S_{r,-}$ be defined in the same way, but with negative orientation.
Then $S_{r,+}$ and $S_{r,-}$ are subsets of~$M_{B,Q}$.

Now $SO(3)$ acts on $\R^3$ preserving the Euclidean metric, and the 
action preserves the equations \eq{vdoteq}. Thus $M_{B,Q}$ 
decomposes into orbits of the $SO(3)$- action. It is easy to show that
the orbits are the single point $(0,0,0)$ and the sets $S_{r,\pm}$
for $r>0$. Moreover, $SO(3)$ acts freely on the set $S_{r,\pm}$, so
that $S_{r,\pm}\cong SO(3)\cong \mathbb{RP}^3$. Therefore $M_{B,Q}$
is the disjoint union of $\{0\}$ and 2 copies of~$(0,\infty)\t\mathbb{RP}^3$.

However, $\H/\{\pm1\}$ is the disjoint union of $\{0\}$ and 1 copy of 
$(0,\infty)\t\mathbb{RP}^3$. Therefore, $M_{B,Q}$ is actually the union 
of {\it two distinct copies} of $\H/\{\pm1\}$, which meet at $0$, and
the map $\pi_{B,Q}:\H/\{\pm1\}\ra M_{B,Q}$ is an isomorphism with
one of these copies. This is something of a surprise. The 2 copies are 
not separate algebraic components of the variety $M_{B,Q}$, as the
complexification of $M_{B,Q}$ in $Q^\d\ot\C$ has only one component
containing both copies.

Let us look at $B$ from the point of view of \S\ref{h44}.
Above we gave an action of $SO(3)$ on $M_{B,Q}$. This action
induces an isomorphism between the $\R^3$ in the equation
$Q=\R^3\ot Y$, with the Lie algebra $\mathfrak{so}(3)$ of $SO(3)$. 
Thus $Q\cong A_{\mathfrak{so}(3)}$, where $A_{\mathfrak{so}(3)}$ is
defined in Example \ref{liealgex}. So $B$ is the quotient
of $F^{A_{\mathfrak{so}(3)}}$ by a stable filtered ideal. Therefore,
$B$ fits into framework of Example~\ref{krcoadex}.

In fact, $S_{r\pm}$ are orbits in $\mathfrak{so}(3)^*\ot\I$, and if
we put $\Omega=S_{r,+}$ or $S_{r,-}$ for any $r>0$, then calculation
shows that $A^\Omega=B$, so that $M^\Omega=M_{B,Q}$. This is the 
simplest example of an HP-algebra of this form, and it is one
of the `nilpotent orbits' mentioned in Conjecture~\ref{krconj}.
\end{ex}

Our goal, remember, is to find the filtered H-algebra $P$ of 
q-holomorphic functions of polynomial growth on the Eguchi-Hanson 
space $M$. From the definition \cite{EH} of the Eguchi-Hanson metric 
one finds that $M$ is AC to order 4, with asymptotic cone $\H/\{\pm1\}$.
Therefore, Theorem \ref{acconethm} applies to show that $P$ is an
SFH-algebra and $B$ and $P$ are isomorphic to order 4, in the
sense of Definition \ref{asympalgdef}. In the next example we will 
determine all such SFH-algebras.

\begin{ex}
\label{detpex}
Suppose $P$ is an SFH-algebra isomorphic to order 4 with the H-algebra
$B$ of the previous examples. As $B_{-2},P_{-2}$ are zero, $\phi_2:B_2\ra 
P_2$ is an \ah-isomorphism. Thus $P_2\cong\H\op Q$. Since $Q$ generates 
$B$, Proposition \ref{asgrgenprop} shows that $P_2$ generates $P$. But 
the $\H$ in $P_2$ is the multiples of 1, so it does not generate anything. 
Therefore $P$ is generated by $Q$, and $P$ is the quotient of $F^Q$ by a 
stable filtered ideal~$I^P$.

Define $J^P=I^P_4$, so that $J^P$ is an \ah-submodule of $\H\op 
Q\op\Sh^2Q$. By a similar argument to Proposition \ref{asgrgenprop}, 
we find that because $J$ generates $I$, $J^P$ generates $I^P$. Thus
$P$ is determined by the \ah-submodule $J^P$. Using the maps 
$\phi_2,\phi_4$ to compare $I^P$ and $I$, we see that $I^P_2=\{0\}$, 
that $J^P\subset\H\op J$. Let $\pi_1,\pi_2$ be the projections
from $\H\op J$ to the $\H$ and $J$. Then $\pi_2:J^P\ra J$ is an 
\ah-isomorphism. Let $\iota:J\ra J^P$ be its inverse, so that
$\iota:J\ra\H\op J$ is an \ah-morphism with~$J^P=\iota(J)$. 

As $\pi_2\circ\iota$ is the identity on $J$, $J^P$ is determined
by $\pi_1\circ\iota$, i.e.~by an \ah-morphism from $J$ to $\H$.
Let $\lambda\in J^\d$. Then $\lambda\in J^\t$, so $\lambda:J\ra\H$
is an $\H$-linear map. By definition of $J^\d$, $\lambda$ maps $J'$
to $\I$ which is $\H'$. Thus $\lambda$ is an \ah-morphism from 
$J$ to $\H$, and it is easy to see that all such \ah-morphisms are
elements of~$J^\d$.

Therefore, we have shown that all SFH-algebras $P$ isomorphic with
$B$ to order 4 (in the sense of Definition \ref{asympalgdef}) are 
constructed in the following manner. We are given $Q=\R^3\ot Y$, and 
a fixed \ah-module $J\subset\Sh^2Q$, that has $J\cong\Sh^2Y$, so that 
$J^\d\cong\R^5$. Choose an element $\lambda\in J^\d$. Then 
$\lambda:J\ra\H$ is an \ah-morphism. Define $J^\lambda$ to be the 
image of the \ah-morphism $\lambda\op\id:J\ra\H\op J$. 

As $\H=\Sh^0Q$ and $J\subset\Sh^2Q$, we may regard $\H\op J$ as an 
\ah-submodule of $F^Q$. Since $J^\lambda$ is an \ah-submodule of 
$\H\op J$, it is an \ah-submodule of $F^Q$. Define $I^\lambda$ by 
$I^\lambda=\mu_{F^Q}(J^\lambda\oth F^Q)$. Then $I^\lambda$ is a stable 
filtered ideal of $F^Q$. Define $P^\lambda$ to be the quotient of $F^Q$ 
by $I^\lambda$, following Lemma \ref{ideallem}. By Definition 
\ref{fghdef}, $P^\lambda$ is an FGH-algebra. 
\end{ex}

Next we identify the variety~$M_{P^\lambda\!,Q}$.

\begin{ex}
\label{plhcex}
Let $\lambda\in J^\d$ be given. Then the previous example defines 
an FGH-algebra $P^\lambda$, generated by $Q\subset P^\lambda$.
We shall describe the real algebraic variety $M_{P^\lambda\!,Q}\subset
Q^\d$. Our treatment follows Example \ref{bfghex} closely, and uses 
the same notation. Now $J^\lambda$ may be interpreted as a vector 
space of $\H$-valued polynomials on $Q^\d$, and $M_{P^\lambda\!,Q}$ is 
the zeros of these polynomials. 

There is a natural identification between $(\Sh^2Y)^\d$ and
$S^2_0\R^3$, which is the space of trace-free $3\t 3$ symmetric 
matrices. Using this identification we may write $\lambda$ in 
coordinates as $\lambda=(a_{kl})$, where $(a_{kl})$ is a $3\t 3$ 
matrix with $a_{jk}=a_{kj}$ and $\Sigma_j\,a_{jj}=0$. Using
the vectors $v_1,v_2,v_3$ as coordinates on $Q^\d$ as before,
the equations defining $M_{P^\lambda\!,Q}$ turn out to be

\begin{equation}
\begin{gathered}
v_1\cdot v_1-a_{11}=v_2\cdot v_2-a_{22}=v_3\cdot v_3-a_{33}, \\
v_1\cdot v_2=a_{12}, \quad v_2\cdot v_3=a_{23}, \quad 
v_3\cdot v_1=a_{31}.
\end{gathered}
\label{avdoteq}
\end{equation}
These are 5 real equations, because $(J^\lambda)^\d\cong\R^5$. When 
$\lambda=0$, $a_{kl}=0$ and we recover the equations~\eq{vdoteq}.

As in Example \ref{vlgex}, the equations \eq{avdoteq} are
invariant under the action of $SO(3)$ on $\R^3$. Thus $M_{P^\lambda\!,Q}$
is invariant under the $SO(3)$- action, and is a union of $SO(3)$ orbits.
Choose any point $m=(v_1,v_2,v_3)$ in $M_{P^\lambda\!,Q}$ such that
$v_1,v_2,v_3$ are linearly independent (this is true for a generic 
point), and let $\Omega$ be the orbit of $m$. Then it can be shown 
that $A^\Omega=P^\lambda$ and $M^\Omega=M_{P^\lambda\!,Q}$, as in
Example \ref{vlgex}. Thus $P^\lambda$ is one of the algebras $A^\Omega$ 
of Example \ref{krcoadex}, with Lie group $G=SO(3)$. Therefore,
by Example \ref{hpcoadex}, $P^\lambda$ is a filtered HP-algebra.
\end{ex}

Finally, we interpret $M_{P^\lambda\!,Q}$, and describe its
hyperk\"ahler structure.

\begin{ex}
\label{ehex}
In Example \ref{plhcex}, the element $\lambda\in J^\d$ determines
a matrix $(a_{jk})$. Now $SO(3)$ acts on $\I$, inducing automorphisms 
of $\H$. But automorphisms of $\H$ act on HP-algebras in an natural
way. Thus $SO(3)$ acts on the category of HP-algebras. This action
of $SO(3)$ on the HP-algebras $P^\lambda$ turns out to be conjugation
of the matrix $(a_{jk})$ by elements of $SO(3)$. So, up to
automorphisms of $\H$, we may suppose that $P^\lambda$ is defined
by a matrix $(a_{jk})$ that is diagonal, and for which~$a_{11}\ge 
a_{22}\ge a_{33}$.

Let the matrix $(a_{jk})$ be of this form. Then the equations
\eq{avdoteq} can be rewritten

\begin{equation}
v_1\cdot v_1-a=v_2\cdot v_2-b=v_3\cdot v_3, \quad
v_1\cdot v_2=0, \quad v_2\cdot v_3=0, \quad v_3\cdot v_1=0,
\label{bvdoteq}
\end{equation}
where $a,b$ are constants with $a\ge b\ge 0$. We shall consider  
the following three cases separately. Case 1 is $a=b=0$, case 2 
is $a>0$, $b=0$, and case 3 is $a\ge b>0$. In case 1, we have
$\lambda=0$, so $P^\lambda=B$ and $M_{P^\lambda\!,Q}$ is 2 copies
of $\H/\{\pm1\}$ meeting at 0, as in Example \ref{vlgex}.

In case 2, there is one special orbit of $SO(3)$, the orbit
$\bigl\{(v_1,0,0):v_1\cdot v_1=a\bigr\}$. Clearly this is a 
2-sphere ${\cal S}^2$. It is easy to see that $M_{P^\lambda\!,Q}$ is 
the union of ${\cal S}^2$ with 2 copies of $(0,\infty)\t\mathbb{RP}^3$.
Also, $M_{P^\lambda\!,Q}$ is singular at ${\cal S}^2$ but nonsingular
elsewhere. A more careful investigation shows that $M_{P^\lambda\!,Q}$
is actually the union of 2 nonsingular, embedded submanifolds
$M_1,M_2$ of $Q^\d$, that meet in a common ${\cal S}^2$.
In fact, both $M_1$ and $M_2$ are the Eguchi-Hanson space.
This case is the HP-algebra of the Eguchi-Hanson space $M$,
and the map $\pi_{P^\lambda\!,Q}:M\ra M_{P^\lambda\!,Q}$ is
a diffeomorphism from $M$ to $M_1$, say.

Thus, we have found the HP-algebra $P^\lambda$ of q-holomorphic
functions of polynomial growth on the Eguchi-Hanson space $M$.
When we attempt to reconstruct $M$ from $P^\lambda$ we find
that $M_{P^\lambda\!,Q}$ contains not one but two distinct copies
of $M$, that intersect in an~${\cal S}^2$.

In case 3, $SO(3)$ acts freely on $M_{P^\lambda\!,Q}$, and
$M_{P^\lambda\!,Q}$ is a {\it nonsingular} submanifold of $Q^\d$
diffeomorphic to $\R\t\mathbb{RP}^3$. This 4-manifold has two ends,
each modelled on that of $\H/\{\pm1\}$. In cases 1 and 2, we saw that
$M_{P^\lambda\!,Q}$ was the singular union of 2 copies of $\H/\{\pm1\}$
or the Eguchi-Hanson space; in case 3 this singular union is resolved
into one nonsingular 4-manifold.

Let us apply the programme of \S\ref{h34}, to construct a hypercomplex
structure on $M_{P^\lambda\!,Q}$. This can be done, and as $Q$ is 
stable, Proposition \ref{4dimprop} shows that the hypercomplex structure is
determined by $P^\lambda$ wherever it exists. At first sight, therefore,
it seems that $M_{P^\lambda\!,Q}$ carries a nonsingular hypercomplex 
structure with two ends, both asymptotic to $\H/\{\pm1\}$. However, 
explicit calculation reveals that although $M_{P^\lambda\!,Q}$ is 
nonsingular, its hypercomplex structure has a singularity on the 
hypersurface $v_3=0$ in $M_{P^\lambda\!,Q}$, which is diffeomorphic 
to~$\mathbb{RP}^3$.

What is the nature of this hypersurface singularity? Consider the
involution $(v_1,v_2,v_3)\mapsto(-v_1,-v_2,-v_3)$ of $M_{P^\lambda\!,Q}$.
This preserves the hypercomplex structure, but it is orientation-reversing 
on $M_{P^\lambda\!,Q}$. It also preserves the hypersurface $v_3=0$. Now a 
hypercomplex structure has its own natural orientation. Therefore, the 
hypercomplex structure changes orientation over its singular 
hypersurface $v_3=0$.

This singular hypercomplex manifold has a remarkable property. Each 
element of $P^\lambda$ is a q-holomorphic function on $M_{P^\lambda\!,Q}$, 
which is smooth, as $M_{P^\lambda\!,Q}$ is a submanifold of $Q^\d$. Thus, 
the hypercomplex manifold has a full complement of q-holomorphic functions 
that {\it extend smoothly over the singularity}. I shall christen
singularities with this property {\it invisible singularities},
because you cannot see them using q-holomorphic functions.
I think this phenomenon is worth further study.

As $P^\lambda$ is an HP-algebra, $M_{P^\lambda\!,Q}$ is actually 
hyperk\"ahler. In Example \ref{plhcex} we showed that this hyperk\"ahler 
structure is $SO(3)$-invariant. Now in \cite{BGPP}, Belinskii et 
al.~explicitly determine all hyperk\"ahler metrics with an $SO(3)$-action 
of this form, by solving an ODE. Thus, the metric on $M_{P^\lambda\!,Q}$ 
is given in \cite{BGPP}. Belinskii et al.~show that the singularity is a 
curvature singularity -- that is, the Riemann curvature becomes infinite 
upon it. 

The metrics can also be seen from the twistor point of view. In 
\cite{Hi}, Hitchin constructs the twistor spaces of some ALE spaces, 
including the Eguchi-Hanson space. He uses a polynomial 
$z^k+a_1z^{k-1}+\cdots+a_k$, introduced on \cite[p.~467]{Hi}. On
p.~468 he assumes this polynomial has a certain sort of factorization, 
to avoid singularities. If this assumption is dropped, then in the case 
$k\!=\!2$, Hitchin's construction yields the twistor spaces of our 
singular hyperk\"ahler manifolds~$M_{P^\lambda\!,Q}$.
\end{ex}

\subsection{Self-dual vector bundles and H-algebra modules}
\label{h46}

First we define q-holomorphic sections of a vector bundle.

\begin{dfn} Let $M$ be a hypercomplex manifold. Let $E$ be a vector 
bundle over $M$, and $\nabla_E$ a connection on $E$. Define 
an operator $D_E:\H\ot C^\infty(E)\ra C^\infty(E\ot T^*M)$ by
\begin{align*}
&D_E\bigl(1\ot e_0+i_1\ot e_1+i_2\ot e_2+i_3\ot e_3\bigr)=\\
&\nabla_Ee_0+I_1(\nabla_Ee_1)+I_2(\nabla_Ee_2)+I_3(\nabla_Ee_3),
\end{align*}
where $I_1,I_2,I_3$ act on the $T^*M$ factor of $E\ot T^*M$,
and $e_0,\dots,e_3\in C^\infty(E)$. We call an element $e$ of 
$\H\ot C^\infty(E)$ a {\it q-holomorphic section} if~$D_E(e)=0$.

Define $Q_{M,E}$ to be the vector space of q-holomorphic sections 
in $\H\ot C^\infty(E)$. Then $Q_{M,E}$ is closed under the 
$\H$-action $p\cdot(q\ot e)=(pq)\cdot e$ on $\H\ot C^\infty(E)$, so 
$Q_{M,E}$ is an $\H$-module. Define a real vector subspace $Q_{M,E}'$ 
by $Q_{M,E}'=Q_{M,E}\cap\I\ot C^\infty(E)$. As in Definition
\ref{adef}, one can show that $Q_{M,E}$ is an \ah-module.
\label{udef}
\end{dfn}

The point of this definition is the following theorem. The proof
follows that of Theorem \ref{qholhalgthm} in \S\ref{h32} very 
closely, so we leave it as an exercise.

\begin{thm} Let\/ $M$ be a hypercomplex manifold, and\/ $A_M$ the H-algebra
of q-holomorphic functions on $M$. Let\/ $E$ be a vector bundle
over $M$, with connection $\nabla_E$. Then the \ah-module $Q_{M,E}$
of Definition \ref{udef} is a module over the H-algebra $A_M$, in a
natural way.
\label{umodthm}
\end{thm}

Recall that modules over H-algebras were defined in Definition
\ref{halgdef}. To see the link between Theorems \ref{qholhalgthm}
and \ref{umodthm}, put $E=\R$ with the flat connection. Then a
q-holomorphic section of $E$ is a q-holomorphic function, so
$Q_{M,E}=A_M$, and Theorem \ref{umodthm} states that $A_M$
is a module over itself. But this follows trivially from 
Theorem~\ref{qholhalgthm}.

Now the equation $D_E(e)=0$ is in general overdetermined,
and for a {\it generic} connection $\nabla_E$ we find that
$Q_{M,E}=\{0\}$ for $\dim M>4$, and $Q_{M,E}'=\{0\}$ for 
$\dim M=4$. In these cases, Theorem \ref{umodthm} is trivial.
To make the situation interesting, we need $\nabla_E$ to
satisfy a curvature condition (integrability condition), 
ensuring that $D_E(e)=0$ has many solutions locally.
We will give this condition. The following proposition 
is a collection of results from~\cite[\S 2]{MS}.

\begin{prop} Let $M$ be a hypercomplex manifold of dimension $4n$. 
Then $I_1,I_2,I_3$ act as maps $T^*M\ra T^*M$, so we may
consider the map $\delta:\Lambda^2T^*M\ra\Lambda^2T^*M$
defined by $\delta=I_1\ot I_1+I_2\ot I_2+I_3\ot I_3$.
Then $\delta^2=2\delta+3$, so the eigenvalues of\/
$\delta$ are $3$ and $-1$. This induces a splitting
$\Lambda^2T^*M=\Lambda_+\op\Lambda_-$, where $\Lambda_+$ 
is the eigenspace of\/ $\delta$ with eigenvalue $3$, and
$\Lambda_-$ is the eigenspace of\/ $\delta$ with 
eigenvalue~$-1$. 

The fibre dimensions are $\dim\Lambda_+=2n^2+n$ and
$\dim\Lambda_-=6n^2-3n$. Also, $\Lambda_+$ is the subbundle 
of\/ $\Lambda^2T^*M$ of 2-forms that are of type $(1,1)$ 
w.r.t.~every complex structure $r_1I_1+r_2I_2+r_3I_3$, for 
$r_1,r_2,r_3\in\R$ and\/ $r_1^2+r_2^2+r_3^2=1$. When $n=1$,
the hypercomplex structure induces a conformal structure, and
the splitting $\Lambda^2T^*M=\Lambda_+\op\Lambda_-$
is the usual splitting into self-dual and anti-self-dual 2-forms.
\label{selfdualprop}
\end{prop}

Following Mamone Capria and Salamon \cite[p.~520]{MS}, we make 
the following definition.

\begin{dfn} Let $M$ be a hypercomplex manifold, let $E$ be a vector 
bundle over $M$, and let $\nabla_E$ be a connection on $E$. 
Let $F_E$ be the curvature of $\nabla_E$, so that 
$F_E\in C^\infty(E\ot E^*\ot\Lambda^2T^*M)$. We say that
$\nabla_E$ is a {\it self-dual connection} if 
$F_E\in C^\infty(E\ot E^*\ot\Lambda_+)$, that is, if the
component of $F_E$ in $E\ot E^*\ot \Lambda_-$ is zero.
\label{selfdualdef}
\end{dfn}

We use the notation $\Lambda_+$, $\Lambda_-$ and self-dual
to stress the analogy with the four-dimensional case, where
this notation is already standard. In four dimensions,
self-dual connections (also called {\it instantons}) are
a very important tool in differential topology, and have
been much studied. 

The point of the definition is this. If $\nabla_E$ is self-dual,
then $F_E$ is of type $(1,1)$ w.r.t.~the complex structure
$r_1I_1+r_2I_2+r_3I_3$ by Proposition \ref{selfdualprop},
and so $\nabla_E$ makes $\C\ot E$ into a holomorphic bundle 
w.r.t.~$r_1I_1+r_2I_2+r_3I_3$. Therefore $\C\ot E$ has many local
holomorphic sections w.r.t.~$r_1I_1+r_2I_2+r_3I_3$. But these
are just special sections $e$ of $\H\ot E$ (or $\I\ot E$) 
satisfying~$D_E(e)=0$. 

We deduce that locally, if $F_E$ is self-dual, then the equation 
$D_E(e)=0$ admits many solutions. Thus, it is clear that vector 
bundles with self-dual connections are the appropriate quaternionic 
analogue of {\it holomorphic vector bundles} in complex geometry.

Suppose that $M$ has a notion of q-holomorphic functions
of polynomial growth, as in Definition \ref{pgdef}, and that
$E$ is equipped with a metric. Sections of $E$ with polynomial 
growth can then be defined in the obvious way. Let $P_{M,E}\subset 
Q_{M,E}$ be the filtered \ah-submodule of q-holomorphic sections
of $E$ of polynomial growth, and let $P_M$ be the filtered H-algebra 
of q-holomorphic functions on $M$ with polynomial growth. Then
$P_{M,E}$ is a filtered module over the H-algebra~$P_M$.

The ideas of \S\ref{h34} can also be applied to the problem of 
reconstructing a self-dual connection from its module $P_{M,E}$. 
Under suitable conditions, $P_{M,E}$ entirely determines the
bundle $E$ and its connection $\nabla_E$. This suggests that
H-algebra techniques could be used to construct explicit
self-dual connections over hypercomplex manifolds of interest. 

Now the {\it ADHM construction} \cite{At} is an explicit 
construction of self-dual connections over the hypercomplex manifold $\H$. 
By Example \ref{hhalgex}, each such connection yields a module 
$P_{M,E}$ over the H-algebra $F^U$. The author hopes to study 
these modules in a later paper, and hence to provide a new 
treatment and proof of the ADHM construction. The approach
generalizes to instantons over other hypercomplex manifolds such as
the hyperk\"ahler ALE spaces, for which an ADHM-type construction is
given in~\cite{KrN}.

\subsection{Directions for future research}
\label{h47}

We have developed an extensive and detailed comparison between 
vector spaces, tensor products and linear maps, and their quaternionic 
analogues. As a result, many pieces of algebra that use vector
spaces, tensor products and linear maps as their building blocks
have a quaternionic version. We shall discuss the quaternionic version 
of algebraic geometry.

Complex algebraic geometry is the study of complex manifolds
using algebras of holomorphic functions upon them. In the same
way, let `quaternionic algebraic geometry' be the study of hypercomplex 
manifolds using H-algebras of q-holomorphic functions upon them. I 
believe that quaternionic algebraic geometry and its generalizations 
may be interesting enough to develop a small field of algebraic geometry 
devoted to them. 

I have tried to take the first steps in this direction in 
Chapters 2-4. These methods seem to have no application to
compact hypercomplex manifolds, unfortunately, so instead one should
study noncompact hypercomplex manifolds satisfying some restriction,
such as AC hypercomplex manifolds. The best category of H-algebras to 
use appears to be FGH-algebras. Here are a number of questions 
I think are worth further study. Almost all are problems in quaternionic 
algebraic geometry, and should clarify what I mean by it.
\medskip

\noindent{\bf Research problems about hypercomplex manifolds}
\begin{itemize}
\item Study the theory of `coadjoint orbit' hyperk\"ahler manifolds and
HP-algebras, begun in \S\ref{h44}, in much greater depth.
\item In Example \ref{ehex} we saw that hypercomplex manifolds derived 
from FGH-algebras can have `invisible singularities' with interesting 
properties. Study these singularities. Develop a theory of the 
singularities possible in `algebraic' hypercomplex manifolds.
\item Develop a {\it deformation theory} for FGH-algebras. That is,
given a fixed FGH-algebra, describe the family of `nearby' 
FGH-algebras. One expects to find the usual machinery of versal and
universal deformations, infinitesimal deformations and obstructions,
cohomology groups. However, the author's calculations suggest that
the H-algebra setting makes the theory rather complex and difficult.
The application is to deformations of `algebraic' hypercomplex manifolds.
\item Use this deformation theory to construct hyperk\"ahler deformations 
of singular hyperk\"ahler manifolds. Can you find new explicit examples 
of complete, nonsingular hyperk\"ahler manifolds of dimension at least~8?
\item In particular, try and understand the deformations of 
$\H^n/\Gamma$, for $\Gamma$ a finite subgroup of $Sp(n)$. In 
\cite{Joy3} the author used analysis to construct a special class
of hyperk\"ahler metrics on crepant resolutions of $\C^{2n}/\Gamma$.
They are called {\it Quasi-ALE metrics}, and satisfy complicated
asymptotic conditions at infinity.

It seems likely that these Quasi-ALE metrics are natural solutions 
to the deformation problem for the FGH-algebra of $\H^n/\Gamma$, and
thus that one could use hypercomplex algebraic geometry to study them,
and even construct them explicitly.
\item Study the H-algebras of ALE spaces, mentioned in \S\ref{h45}.
Use H-algebras to give a second proof of Kronheimer's classification
of ALE spaces \cite{Kr1}, \cite{Kr2}. 
\item Rewrite the ADHM construction \cite{At} for self-dual
connections on $\H$ in H-algebra language, as suggested in \S\ref{h46}.
Interpreted this way it becomes a beautiful algebraic construction
for modules of the H-algebra $F^U$ of Example~\ref{hhalgex}.
\item In a similar way, rewrite the ADHM construction on ALE spaces
\cite{KrN} in terms of H-algebras. Consider the possibility of a
`general ADHM construction' for FGH-algebras and algebraic hypercomplex 
manifolds, that captures the algebraic essence of the constructions
of \cite{At} and~\cite{KrN}.
\end{itemize}
\medskip

\noindent{\bf Other research problems}
\begin{itemize}
\item Understand and classify finite-dimensional HL-algebras,
by analogy with Lie algebras.
\item Apply quaternionic algebra in other areas of mathematics, to produce
quaternionic analogues of existing pieces of mathematics. These need have
no connection at all with hypercomplex manifolds. For instance, one can
look at a quaternionic version of the Quantum Yang-Baxter Equation, and
try and produce quaternionic knot invariants.
\item Generalize the quaternionic algebra idea by replacing $\H$ by some
noncommutative ring, or more general algebraic object. Can
you find any interesting algebraic applications for this?
One particularly interesting case is when $\H$ is replaced by a 
division ring over an algebraic number field, since these rings
have a strong analogy with~$\H$.
\item We have seen that one may associate an H-algebra to a hypercomplex 
manifold. However, the converse is false, and one cannot in general
associate a hypercomplex manifold to an H-algebra. Following \S\ref{h34},
given an H-algebra $P$ generated by finite-dimensional $Q$, one 
constructs a real algebraic variety $M_{P,Q}\subset Q^\d$. What
geometric structure does $P$ induce on $M_{P,Q}$, for general $P$\,?
Study these structures.
\item In \cite{Joy1} the author constructed nontrivial examples of 
manifolds with $n$ anticommuting complex structures, for arbitrary
$n>1$. Such structures have also been found by Barberis et 
al.~\cite{BMM}, who call them {\it Clifford structures}. The 
programme of this paper may be generalized to Clifford structures. 
One simply replaces $\H$ by the Clifford algebra $C_n$ in the definition 
of H-algebra, and then to a manifold with a Clifford structure one
may associate a `$C_n$-algebra', the analogue of H-algebra.
Use this theory to study and find new examples of manifolds with
Clifford structures.
\end{itemize}

\end{document}